\documentclass[12pt,fleqn]{report}
\usepackage[matrix,arrow,graph,curve]{xy}
\usepackage{amssymb}
\newtheorem{thm}{Theorem}[section]

\newtheorem{lem}{Lemma}[section]

\begin{document}
\thispagestyle{empty}
\begin{center}
\bf \LARGE{TRIANGULATIONS OF SPHERES AND DISCS}
\end{center}
\vspace{25mm}

\begin{center}
\LARGE {PANCHADCHARAM ELANGO\\}
\vspace{40mm}
\Large {INSTITUTE OF MATHEMATICAL SCIENCS \\ FACULTY OF SCIENCE 
\\ UNIVERSITY OF MALAYA \\}
\vspace{40mm}

\Large{DISSERTATION PRESENTED FOR THE \\ DEGREE OF MASTER OF SCIENCE \\
UNIVERSITY OF MALAYA \\ KUALA LUMPUR \\ (2002)}
\end{center}

\newpage
\pagenumbering{roman}
\setcounter{page}{1}

\centerline{\Large \bf{Acknowledgement}} 
\bigskip 

I would like to express my sincere thanks to my former supervisor, Associate Professor Dr. Thomas Bier, for his unlimited assistance  
and encouragement throughout the course of my study 
and for helping me to complete my research work on time. 
\\

I would also like to express my sincere thanks to my present supervisor, 
Associate Professor Dr. Angelina Chin Yan Mui for taking on the task
of supervising me after
my former supervisor went on leave. Her guidance and valuable suggestions
helped me a lot to complete my thesis.  
\\

Many thanks also to the Head and 
all staff of the Institute of Mathematical Sciences, University of Malaya 
for their assistance in numerous ways.  \\

Finally, I wish to express my special thanks to my loving mother and dearest wife for their encouragement to successfully complete my MSc studies. \\

\vspace{0.8cm} 

\noindent
Panchadcharam Elango, 2002.

\newpage
\centerline{\Large \bf{Abstract}}
\bigskip

The main objective of this research is to find the different types
of elliptic 
triangulations for planar discs and spheres. 
We begin in Chapter 1 with the mandatory introduction.  
In the second chapter we define and study the notion of a 
{\it patch}, that is, a triangulation of a planar disc. 
By introducing 
a suitable notion of degree for the vertex points, 
we focus on those patches with points of degrees $\le 6$.
Such patches are called {\it elliptic}.   
We show that the elliptic patches with precisely three points of degree 
four, denoted by $(0,3,0)$, can be classified. 
The number of vertex points 
of these patches are calculated, and we also describe 
their triangulation structures.  
From the classification of patches of the type
$(0,3,0)$, we describe and find the number of vertex points
for three other elliptic patches of types $ (0,2,2), (0,1,4), (0,0,6)$.  
We also describe an enlargement method for constructing patches
(which we call the {\it generic construction} method)
and apply this method to derive formulas for certain patches. 
In the third chapter we describe some configurations for constructing
elliptic spherical triangulations. These are the  
mutant configuration, the productive configuration and 
the self-reproductive
configuration. We also describe the face-fullering and edge-fullering
methods as well as the glueing of patches method
for constructing 
triangulations and patches. We show that there are only
19 possible types of elliptic triangulations for spheres
and determine the existence (as well as nonexistence) 
of all these types except
for a small number of cases.

\newpage
\centerline{\Large \bf{Abstrak}}
\bigskip

Tujuan utama penyelidikan ini adalah untuk mendapatkan jenis triangulasi
eliptik yang berbeza bagi cakera satahan dan sfera.   
Kita mula dalam Bab 1 dengan pengenalan. 
Dalam Bab 2, kita menakrif dan mengkaji konsep {\it tampal}, 
iaitu, triangulasi
bagi cakera satahan.
Dengan memperkenalkan takrif yang sesuai bagi darjah sesuatu
bucu, kita  menumpu kepada tampal yang mempunyai bucu berdarjah
$\le 6$. Tampal seperti ini dikatakan {\it eliptik}.    
Kita tunjukkan bahawa tampal eliptik dengan tiga bucu berdarjah
empat, ditulis sebagai $(0,3,0)$, boleh diklasifikasikan. 
Bilangan bucu untuk tampal seperti ini didapati, dan kita juga 
memperihalkan struktur triangulasinya.
Daripada klasifikasi tampal jenis
$(0,3,0)$, kita memperihalkan dan cari bilangan bucu bagi tiga lagi
tampal eliptik jenis $(0,2,2), (0,1,4), (0,0,6)$.  
Kita juga memberikan satu kaedah `pembesaran' untuk membina tampal
(yang dipanggil {\it kaedah membina generik}) dan gunakan kaedah ini
untuk mendapatkan formula bagi beberapa tampal. 
Dalam Bab 3 kita memberikan beberapa konfigurasi untuk membina  
triangulasi sfera eliptik. Ini termasuk konfigurasi ``mutant'', 
konfigurasi produktif dan konfigurasi swa-produktif semula.  
Kita juga menerangkan kaedah ``face-fullering'' dan ``edge-fullering''   
serta kaedah ``mengglu tampal'' untuk membina triangulasi dan tampal.
Kita tunjukkan bahawa hanya terdapat 19 kemungkinan untuk 
triangulasi eliptik bagi sfera dan tentukan kewujudan
(serta ketakwujudan) untuk semua jenis ini kecuali bagi sebilangan kecil  
kes.  

\newpage

\tableofcontents
\listoffigures
\listoftables
\include{intro}
\include{ch1}
\include{ch2}
\include{ch3}

\pagenumbering{arabic}
\setcounter{page}{1}

\chapter{Introduction}
A triangulation of a surface is a set $V$ of points, a set $E$ of $2$-subsets of 
$V$ and a set $T$ of $3$-subsets of $V$ satisfying the following conditions:
\begin{description}
\item(a) The union of elements of $T$ is connected;
\item(b) Every $2$-subset in $E$ is contained in precisely 
two $3$-subsets in $T$;
\item(c) The union of elements of $T$ containing a fixed point $x$ 
has a boundary
consisting of a single cycle, that is, they are of the form 
$$xx_1x_2, xx_2x_3, \dots , xx_{n-1}x_n, xx_nx_1.$$
\end{description}
A point in $V$ is called a {\it vertex}, a $2$-subset in $E$ is called an
{\it edge} and a $3$-subset in $T$ is called a {\it triangle}.
In this thesis we shall investigate triangulations of spheres
and planar discs.   

We begin by considering triangulations of planar discs in Chapter $2$. We call
such a triangulation a {\it patch} and show that a patch can be extended to a 
triangulation of a sphere. The {\it degree} of a vertex $x$ in a patch 
$P$ can be defined
as follows:
If the point $x$ is inside the patch $P$ (that is, does not lie on the boundary
of the patch),   
then the number of triangles $ \sigma_1,\sigma_2, \dots ,\sigma_d$ 
which contain $x$ is denoted by  $ d=d(x)$ and is called 
the degree of $x$.    
If the point $x$ is on the boundary of the patch $P$
and is incident with exactly $d'(x)$ edges 
(including the two boundary edges), then we define 
the degree of $x$ as the integer $d = d'(x) + 2.$  
Note that in this case $x$ lies in $d'(x)-1$ triangles.
A patch $P$ is said to be {\it elliptic}    
if the degrees of all the points in $P$   
are not greater than $6$, that is, $d(x) \leq 6$ for all $x \in P$.   
For a patch $P,$ we let $\alpha_d = \alpha_d(P)$ denote the number of points
of $P$ of degree $d$. By using Euler's equation, we get   
$$3 \alpha_3 + 2 \alpha_4 + \alpha_5    - \alpha_7 - 2\alpha_8 -\dots
 - (m-6)\alpha_m = 6$$
which gives us the equation 
$$3 \alpha_3 + 2 \alpha_4 + \alpha_5 = 6$$
in the case when $P$ is an elliptic patch. This latter equation
has only seven possible nonnegative solutions, that is, 
$$(\alpha_3, \alpha_4, \alpha_5)=
(2,0,0); (0,3,0); ( 1,1,1); (1,0,3); (0,2,2); (0,1,4); (0,0,6).$$
We call the $3$-tuple $(\alpha_3, \alpha_4, \alpha_5)$ the {\it type} 
of the patch $P$.   
 
The definition of patch first appeared in the literature in a paper by 
W. T. Tutte \cite{tutte}. 
B. Gr\"unbaum and T.S. Motzkin \cite{GM} implicitly used
patches of type $(2,0,0)$ to construct and classify triangulations
of spheres of the type $(4,0,0)$.   

In Section $2.1,$ we discuss patches of the type $(2,0,0)$ and give a 
general formula for the number of points of degree $6$.  
In Section $2.2,$ we discuss patches of the type $(0,3,0)$
with some boundary conditions. We derive general formulas 
for the number of points of degree $6$ of this type of patches
and describe the structures 
of the patches with those boundary conditions.  
In Section $2.3,$ we construct certain patches of types 
$(0,2,2),(0,1,4)$ and $(0,0,6)$. 
We then obtain further properties of patches of these types in Section 2.4
by using the classification of patches of the type 
$(0,3,0)$ obtained in Section 2.2. 
In Section $2.5$ we describe an enlargement
process for constructing patches, which we call {\it generic patches}. 
We then use
this construction method to construct patches of the type $(1,1,1)$ 
in Section $2.6$.  
In Section $2.7,$ the final section in Chapter $2,$ we give the construction of some
well-known patches and some patches which are used in Chapter $3$ (including 
some
of the type $(1,0,3)$). 

In Chapter $3,$ we consider triangulations of spheres. 
The {\it degree} of a vertex $x$ in a triangulation $T$ of a sphere is
the number of triangles $ \sigma_1,\sigma_2,\dots,\sigma_d$ 
which contain $x$ and is denoted by  $ d=d(x).$ 
A triangulation $T$ is said to be {\it elliptic} 
if it does not contain any point with degree
greater than $6$, that is, $d(x) \leq 6$ for every $x \in T.$
Again we use Euler's equation to get  
$$3 \alpha_3 + 2 \alpha_4 + \alpha_5 - \alpha_7 - 2\alpha_8 
 -\dots - (m-6)\alpha_m = 12,$$
which reduces to
$$3 \alpha_3 + 2 \alpha_4 + \alpha_5 = 12$$
in the elliptic case. There are 
$19$ nonnegative solutions $(a_3,a_4,a_5)$
for this equation. As in Chapter $2,$ we call $(a_3,a_4,a_5)$ the type
of the triangulation $T.$ It has been shown by Eberhard \cite{Eber} that 
for each of the solution $(a_3,a_4,a_5)$, there
exist a triangulation $T$ and a nonnegative integer     
$N=a_6 \geq 0$ with the property
$$ (\alpha_3(T),\alpha_4(T), \alpha_5(T), \alpha_6(T)) = (a_3,a_4,a_5,a_6).$$ 

Our main contribution in Chapter $3$ is to find, for each of the $19$ types
of triangulations, all possible values of $N=a_6.$
A summary list of these 19 types 
and what is known about the possible values of $N$
are given in the appendix of this thesis. 
We also describe in Chapter $3$ various methods to construct elliptic
spherical triangulations such as the mutant, productive and self-reproductive
configurations, the fullering constructions and the glueing of patches method.

We remark here that some non-existence results on triangulations have been obtained
by Gr\"unbaum \cite{GR}. Eberhard \cite{Eber} and Br\"uckner \cite{bruck} 
have determined the minimum values of $N$ such that the triangulations 
of type $(a_3,a_4,a_5,N)$ exist for each of the $19$ 
nonnegative solutions
$(a_3,a_4,a_5)$. 

Other references that deal with related aspects of triangulations are 
Barnette \cite{Bar3}, K\"uhnel and Lassmann \cite{KuLa},   
\cite{KuLa2} and Negami \cite{sei}.

\newpage
\chapter{Triangulations of Discs (Patches)}
Let $P$ be a triangulation of the unit circle 
\begin{equation} \label{U1}
D^2 = \{ (s,t)\in \mathbb R^2 ~:~s^2 + t^2 \leq 1 ~\}
\end{equation}
with boundary 
\begin{equation} \label{U2}
S^1 = \{ (s,t)\in \mathbb R^2 ~:~s^2 + t^2 = 1 ~\}. \end{equation}
$P$ will be called a {\it patch}.
For any point $x\in P$, 
we define the {\it degree $d$} of $x$ as follows: 
If the point $x$ is in the interior $Int(P)$ of the patch $P$  
(that is, $x$ does not lie on the boundary of $P$),  
then the number of triangles $ \sigma_1,\sigma_2, \dots ,\sigma_d$ 
which contain $x$ is denoted by  $d=d(x)$ and is called 
the degree of $x$.    
If the point $x$ is not inside the patch but on the boundary of the patch $P$,
that is, $x \notin Int(P)$ and $ x\in P \cap S^1 $ 
(denoted as $ x \in \partial P$),
and $x$ is incident with exactly $d'(x)$ edges 
(including the two boundary edges),
we define the degree of $x$ as the integer $d= d'(x) + 2$.    
Note that in this case $x$ lies in $d'(x)-1$ triangles.
We can easily see that $ d \geq 3$ if $x \in Int(P)$ and that 
$ d \geq 4$ if $x \in \partial P.$ From now on the degree of a point
$x \in P$ will be written as $d(x)$ or just $d$.      

A patch $P$ is said to be {\it elliptic} 
if the degree of all the points in $P$ 
is not greater than $6$, that is, $d(x) \leq 6$ for all $x \in P.$ Hereafter
we use the word patch to refer to an elliptic patch unless otherwise stated. 

Let $ \alpha_d $ be the number of points (lying on the boundary or not) 
of degree $d$ in a patch $P$. That is,  
\begin{equation} \label{U3}
\alpha_d(P) = \alpha_d = | \{ x \in P~:~ d(x) = d \} |. \end{equation}
The collection of the $\alpha_d$ is called the parameters of $P.$ 
In this case we obtain the equation 
\begin{equation} \label{U4}
3 \alpha_3 + 2 \alpha_4 + \alpha_5    - \alpha_7 - 2\alpha_8 -... - (m-6)\alpha_m = 6 ~.
\end{equation}
This is a consequence of the fact that the given triangulation can 
be extended to a triangulation of a sphere $T$ by adding  
a patch $Q$.    
That is, $T=P\cup Q$ where 
$Q$ is another patch with $b+1$ points and boundary length $b$ such that   
$b$ points lie on the boundary of $Q$ and the remaining   
one point is in the interior of $Q$ with degree $b$. 
If $\gamma_i$ denotes the corresponding parameter for this 
triangulation $T$, we have $~\gamma_5=\alpha_5+b,~~~\gamma_b=\alpha_b+1,~~~
\gamma_d=\alpha_d {\mbox ~~for~~} d \not =5,b.$
Hence by {\it Euler's} equation
\begin{eqnarray*}
\hspace{2.9cm} 
\sum_{d \geq 3}~(6-d)\gamma_d  =  12, 
\end{eqnarray*}
we get  
\begin{eqnarray*}
\hspace{2.9cm} 
\sum_{d \geq 3} (6-d) \alpha_d  + (6-b)\cdot 1 + b  =  12,  
\end{eqnarray*}
that is, 
\begin{eqnarray*}
\hspace{2.9cm} \sum_{d \geq 3} (6-d) \alpha_d  =  6 
\end{eqnarray*}
which gives us (\ref{U4}). 
In the particular case of an elliptic triangulation,   
$d(x) \leq 6$ for all $x \in P$ and
there are only seven possible  
nonnegative solutions for (\ref{U4}): 
\begin{eqnarray} \label{TU3}
(\alpha_3, \alpha_4, \alpha_5) &=&  
(2,0,0); (0,3,0); ( 1,1,1); (1,0,3); (0,2,2); \nonumber \\
&&\qquad  (0,1,4); (0,0,6).
\end{eqnarray}
The $3$-tuple 
$(\alpha_3, \alpha_4, \alpha_5)$   
in (\ref{TU3}) 
will be called the {\it type} of the patch. \\

Let $P$ be a patch with boundary length $b$, that is, $b$ boundary
edges. If we introduce one new point opposite 
(on the outside) to each edge of the boundary of $P$ 
and connect these new points with each other
and with the vertex points of the corresponding 
edges of $P$, we will get a strip of $2b$ new
triangles around the patch $P.$ This strip of triangles is called
a {\it belt} of $P.$

We note that if there exists a patch with 
$b$ boundary points, with a given set of (nonnegative) parameters like 
$ (a_3,a_4,a_5,a_6)$, then we can also get another patch 
with parameters $(a_3,a_4,a_5,a_6+b)$ by enlarging the first patch 
with a belt of points of degree six. 

For each given number $b \geq 3$ we may then ask for the list 
of $4$-tuples $(a_3, a_4, a_5, a_6)$ 
such that $ 3a_3+2a_4+ a_5 = 6$, $0 \leq a_6 \leq b$ and such that there 
exists a patch $P$ with precisely $b$ boundary points 
and an integer $m \geq 0$ satisfying  
\begin{equation} \label{U6}
( \alpha_3(P) , \alpha_4 (P), \alpha_5(P), \alpha_6(P)) =
 ( a_3,a_4,a_5,a_6+b \cdot m). ~ 
\end{equation}
If such a $4$-tuple of nonnegative integers exists, we shall 
denote it by $(a_3,a_4,a_5$, $a_6)_b.$ 
\\

It is easy to see that the removal of all boundary points, and of all 
edges and triangles of the boundary points from a patch $P$ 
gives another patch $P'$, or a $1$-dimensional complex (that is, a 
graph), or a $0$-dimensional complex (that is, just some points).   
The types of $P$ and $P'$ are not necessarily the same. \\
 
In what follows we use the usual terminology concerning patches. 
In particular, 
we write that $\beta_k=m $ if there are precisely $m$ points   
of degree $k < 6$ on the boundary. 
We also say that a point $x$ of degree $k$ is {\it almost on the boundary }
if $x$ is not on the boundary, but is contained in a triangle 
which has an edge of the boundary. In other words, a point is 
almost on the boundary iff it lies opposite to a boundary edge. \vskip3mm

We define the face numbers $f_i(P)$ ($i=1,2,3$) 
for a given patch $P$ as follows: 
Let $f_1(P)$ be the number of points, $f_2(P)$ the number of edges
and $f_3(P)$ the number of triangles for the given patch $P.$ 
Let the parameters of $P$ be $(a_3,a_4,a_5,a_6)_b$ 
(that is, $P$ has boundary length $b$).  
Then we can get the following relations involving  
the face numbers of $P$ and its parameters.   

\begin{equation} \label{face1}
f_1(P)=a_3+a_4+a_5+a_6,
\end{equation}
\begin{equation} \label{face2}
f_2(P)=3f_1(P)-(3+b),
\end{equation}
\begin{equation} \label{face3}
f_3(P)=2f_1(P)-(2+b).
\end{equation}
\medskip

\begin{lem}  \label{L1}
Let $P$ be a patch with a boundary consisting of points of degree 
$6$ only. Then the removal of all the boundary points, and of all the 
edges and triangles containing boundary points gives a new patch $P'$ 
of the same type as $P,$ or a graph, or just some points. 
If a patch $P'$ is obtained, then the face numbers of $P'$ satisfy   
\begin{equation} \label{EL1}
f_1(P') = f_1(P) - b,~~ f_2(P') = f_2(P) - 3b,~~ f_3(P') = f_3(P)- 2b.
\end{equation}  
\end{lem}
\smallskip

\noindent
We say that $P'$ is obtained from $P$ by {\it peeling a belt} of length $b.$ \\

\section{Patches of Type $(2,0,0)$} \label{patch200}
We first discuss a case in which the peeling of a belt may result in a graph. 
This may happen in the case of type $(2,0,0),$ where a patch consists of 
a boundary of length $b=2k$ and of $ m= (k-1) + 2kr $ points of degree $6,$ 
and two points of degree $3.$ Here $k-1$ points constitute a graph 
that forms a path of length $k$ between the two points of degree $3.$
The other $2kr $ vertex points of degree $6$ are arranged in 
$r$ belts of equal length around this path. It has been shown implicitly 
in \cite{GM} as a by-product of the classification 
of spherical triangulations of type $(4,0,0,m)$ that each patch of 
type $(2,0,0)$ has this form. In other words, 
for each patch of type $(2,0,0)$ there exist two integers $k \geq 2,r \geq 0 $ 
such that the number of vertex points of degree $6$ of the patch is 
$ m = 2kr+k-1,$ and the length of the boundary of the patch is $ b = 2k.$ 
Hence the classification of patches of type $(2,0,0)$ can essentially 
be deduced from the results in \cite{GM}. \\
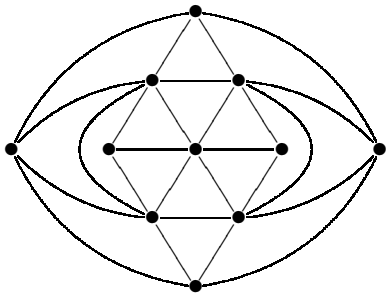
\begin{figure}[htb]
\hspace {4cm} 
\xymatrix @M=0ex@R=4ex@C=2ex{
&& & & & {\bullet}\ar@{-}[dr]="c" \ar@{-}[dl]="a" \ar@/^1pc/@{-}[ddrrrrr] \ar@/_1pc/@{-}[ddlllll]
	& & & && \\
&& & & {\bullet} \ar@/_/@{-}[dllll] \ar@{-}[dr] \ar@{-}[dl] \ar@{-}[rr]&  & {\bullet} \ar@/^/@{-}[drrrr] 
	\ar@{-}[dr]\ar@{-}[dl]& & && \\
{\bullet} &  && {\bullet}\ar@{-}[dr] ="b"\ar@{-}[rr]& & {\bullet} \ar@{-}[dr]\ar@{-}[dl]\ar@{-}[rr]& & 
	{\bullet}\ar@{-}[dl]="d"& &&{\bullet}\\
& && & {\bullet}\ar@/^/@{-}[ullll]\ar@{-}[dr]\ar@{-}[rr]&  & {\bullet} \ar@/_/@{-}[urrrr]\ar@{-}[dl]& & && \\
& && & & {\bullet} \ar@/_1pc/@{-}[uurrrrr]\ar@/^1pc/@{-}[uulllll]& & && &
\ar@/_2.3pc/ @{-}"a";"b" \ar@/^2.3pc/ @{-}"c";"d"}
\caption{Patch of type $(2,0,0,m)$ with $m=9,r=2,k=2$.}
\label{TyXX}
\end{figure}

\section{Patches of Type $(0,3,0)$} \label{pat030}
 
In the following we assume that the patches under consideration have more than 
$6$ vertex points, and have $ \beta_4>0.$  

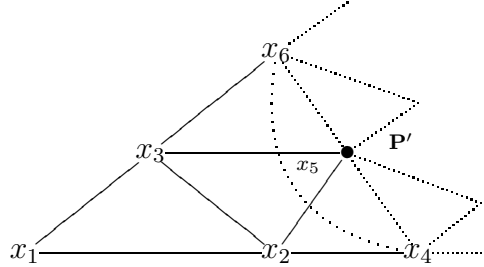
\begin{figure}[htb]
\hspace {3.5cm} 
\xymatrix @M=0ex@R=3ex@C=3.5ex{
& & & & & &  &\\
& & & & {x_6} \ar@{-}@{.}[ur] \ar@{-}@{.}[drr] \ar@{-}@{.}[ddr]& & &  \\
& & & & & & & \\
&  &{x_3} \ar@{-}[uurr]="a" \ar@{-}[rrr]_(.8){x_5} \ar@{-}[ddrr] \ar@{-}[ddll]& & &{\bullet} 
	\ar@{-}@{.}[drr]\ar@{-}@{.}[ddr]="b" \ar@{-}@{.}[ur]_-{\bf{P'}}& & \\
& & & & & &  & \\
{x_1} \ar@{-}[rrrr] & & & &{x_2} \ar@{-}[uur] \ar@{-}[rr] & & {x_4} \ar@{-}@{.}[ur] 
\ar@{-}@{.}[r]& &\ar@/_2pc/ @{-}@{.}"a";"b"}
\caption{Cutting the corner at $x_1$.}
\label{cutcorn12}
\end{figure}

If $x_1$ is a vertex on the boundary of degree $4,$ 
then let $ x_1x_2x_3 $ be the triangle containing $x_1$.    
Let $x_4,\,x_5,\,x_6$ be three points of (graphical)
distance $2$ from the vertex $x_1$ such that 
$ x_2x_4x_5, x_2x_3x_5, x_3x_5x_6 $ form triangles and assume that   
$ x_4x_6 $ do not form an edge. 
We may construct a new patch $P'$ by replacing the four 
triangles $ x_1x_2x_3,x_2x_4x_5, x_2x_3x_5, x_3x_5x_6 $
by a single triangle $ x_4x_5x_6$ (see Figure \ref{cutcorn12}). The patch $P'$ satisfies 
\begin{equation} \label{EP}
b' = b-3,~~ \beta_4' = \beta_4-1. 
\end{equation}
We will call the replacement of $P$ by $P'$ as {\it cutting the corner at } $x_1.$ 

If the patch is sufficiently large, it is clear that we can cut all (at most 3) 
corners until we obtain a new patch $\tilde{P}$ with $\beta_4 = 0.$ 
This means that all boundary 
points of $\tilde{P}$ have degree $6$.   
By peeling the belt around $\tilde{P} $ 
we obtain a patch $ P_0$ that is also of type $ (0,3,0)$.

In particular, if $\beta_4 = 2, $ then there exists precisely one point $q$ of degree 
$4$ that is not on the boundary of $P$.   
From the fact that $P$ has finitely many vertices, 
it is clear that there exists a minimal positive integer $k>0$ 
such that after cutting corners and peeling the resulting belts 
$k$ times there will result a patch $\tilde{P}$ which contains 
$q$ on the boundary. In this case we say that $q$ 
has distance $k$ from the boundary. 
For the patches of type $(0,3,0)$ with $ \beta_4 = 2$, we write  
$ [0,0,k]$ if the  point of degree $4$  that is not on the boundary 
has distance $k$ from the boundary. 

In the case when $\beta_4 =1,$ there exist precisely two points $q,r$ of degree 
$4$ that are not on the boundary of $P.$ From the fact that $P$ has finitely many vertices, 
it is clear that there exists a minimal positive integer $l>0$ 
such that after cutting corners and peeling the resulting belts 
$l$ times there will result a patch $\tilde{P}$ which contains 
(without loss of generality) 
$q$ on the boundary, and both $q,r$ not on the boundary 
after $l-1$ steps. In this case we say that $q$ has distance $l$ 
from the boundary. If $r$ is also on the boundary after 
removing the $l$ belts, 
we say that $r$ also has distance $l$ from the boundary. 
If $r$ is not on the boundary after removing $l$ belts, we use  
the previous case $ \beta_4 =2$ to define the distance $k$ of the point $r$ from the boundary. 

For the patches of type $(0,3,0)$ with $\beta_4=1$   
and $b$ boundary points,    
we write $P[0,l,k]_b$ or just $[0,l,k]$ 
if the two points of degree $4$ that are not on the boundary    
have distance $l$ and $k$ from the boundary respectively. We will prove
the following:

\begin{thm} \label{T1}
\begin{description}
\item{(i)} If there exists a patch with parameters $(0,3,0,N)_b,$ with $b$ boundary points, 
then the integer $b$ is a multiple of three, that is, $ b= 3h$ for some
integer $h.$  
\item{(ii)} There exists a patch with parameters $(0,3,0,N)_{3h} $ with 
$ \beta_4 =3$ for all values of $ h \geq 1$ iff 
\begin{equation} \label{E0}
 N = {h+2 \choose 2}-3.  \end{equation} 
The distance between any two of the three points of degree $4$ on the boundary is equal to $h.$  
\item{(iii)} There exists a patch with parameters $(0,3,0,N)_{3h},$ 
$ \beta_4 = 2$ and of the form 
$[0,0,k]$ with $h=k+l$ for all values of $l>0$ iff 
\begin{equation} \label{EET}
 N = {h+k +2 \choose 2 } - { 2k+1 \choose 2} + { k \choose 2} - 3 . 
\end{equation} 
The lengths of the two boundary parts between the two points of 
degree $4$ on the boundary are $h+l$ and $2h-l.$ 
\item{(iv)} There exists a patch with parameters $(0,3,0,N)_{3h}$ 
with $ \beta_4 = 1$  iff 
\begin{equation} \label{E2}
 N = {h+k+l +2 \choose 2 } - 3 { k+1 \choose 2} - 3{ l+1 \choose 2} - 3 . \end{equation} 
This patch is of the form $[0,l,k]$ where $ 0 < l \leq k < h.$

\end{description}
\end{thm} 
First we show the existence of patches of the types mentioned in parts
$(ii),(iii)$ and $(iv)$ of Theorem \ref{T1}.

\subsection{Construction of Patches $(0,3,0)$ with $\beta_4 = 1,2,3.$ }
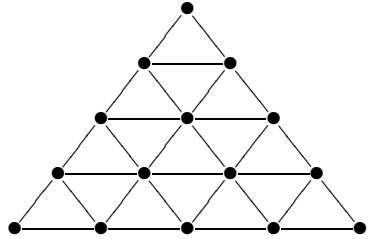
\begin{figure}[htb]
\hspace {3.5cm} 
\xymatrix @M=0ex@R=3ex@C=2ex{
& & & & {\bullet}\ar@{-}[dr]\ar@{-}[dl]& & & & \\
& & &{\bullet} \ar@{-}[dr]\ar@{-}[dl]\ar@{-}[rr]& & {\bullet}\ar@{-}[dr]\ar@{-}[dl]& & & &\\
& &{\bullet} \ar@{-}[dr]\ar@{-}[dl]\ar@{-}[rr]& & {\bullet}\ar@{-}[dr]\ar@{-}[dl]\ar@{-}[rr]& 
	&{\bullet}\ar@{-}[dr]\ar@{-}[dl] & & \\
&{\bullet} \ar@{-}[dr]\ar@{-}[dl]\ar@{-}[rr]& &{\bullet}\ar@{-}[dr]\ar@{-}[dl] \ar@{-}[rr]& 
	&{\bullet} \ar@{-}[dr]\ar@{-}[dl]\ar@{-}[rr]& & {\bullet}\ar@{-}[dr]\ar@{-}[dl]& \\
{\bullet}\ar@{-}[rr]& & {\bullet}\ar@{-}[rr]& & {\bullet}\ar@{-}[rr]& &{\bullet} \ar@{-}[rr]& &{\bullet}}
\caption{Patch $P_4[0,0,0]$ of type $(0,3,0,N)$ with $N=12,\beta_4=3$.}
\label{cutcornq1}
\end{figure}

We first demonstrate that for each $N$ of the form 
$$ N = {h+2 \choose 2}-3,  $$  
there exists a patch with parameters $ (0,3,0,N)_b$ 
with $b=3h,\,\beta_4 = 3$.  \vskip1mm
This construction is obvious as a tessellation of an equilateral triangle 
by $h$ layers of smaller equilateral triangles. This patch will be denoted by 
$P_h= P_h[0,0,0]$ or $P[0,0,0]_{3h}$.    
For an example with $N=12$ (that is, of 
$P_4=P_4[0,0,0]$), see Figure \ref{cutcornq1}.

\begin{figure}[htb]
\hspace {3.5cm} 
\xymatrix @M=0ex{
& & & & \\
& {\bullet}\ar@{-}[rr] & & {\bullet} & \\
{\bullet}\ar@{-}[rrrr]_-{lower boundary}\ar@{-}[ur]="a"^-{left edges} & & & & {\bullet}
	\ar@{-}[ul]="b"_-{right edges} \ar@/^2.0pc/ @{-}@{--}^(0.8){upper boundary}"a";"b"}
\caption{The Patch $Tr_g(P_k)$.}
\label{idlow}
\end{figure}
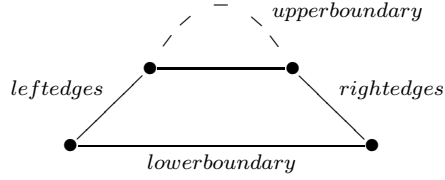

For each $N$ of the form 
$$ N = {h+k+2 \choose 2}-{ 2k+1 \choose 2} + { k \choose 2}- 3\quad
(0< k< h),$$  
we show that
there exists a patch with parameters $(0,3,0,N)_b$ of the form $[0,0,k]$ 
with $ b=3h,~\beta_4 = 2.$  
This patch will be denoted by $P_{h-k}[0,0,k]$ or just 
$[0,0,k]$.   
We can construct this patch by using a truncation 
of the triangular tessellation 
$P_h$ above. For $0 < g \le k$, we may regard  
$ P_g[0,0,0] $ as a subtriangulation of $ P_k[0,0,0]$    
by including it in such a way that one of the points of degree $4$ 
in $P_g$ and in $P_k$ are identical. We then define the 
{\it truncated tessellation } 
$ Tr_g(P_k) $ as the remaining part of the triangulation, that is,   
\begin{equation} \label{Trdef}
Tr_g(P_k) = P_k[0,0,0] \backslash P_g[0,0,0] . \end{equation} 
The boundary edges of $P_g[0,0,0]$ which are not boundary edges 
of $P_k[0,0,0] $ will be called the {\it upper boundary } of 
$ Tr_g(P_k).$ The part of the boundary of  $P_k$ 
consisting of the $k$ edges on one straight line 
all not in $P_g$ is called the {\it lower boundary } of $Tr_g(P_k).$ 
The two remaining parts of the boundary of 
$Tr_g(P_k)$ are called the {\it left} and {\it right}   
edges of $Tr_g(P_k)$ (see Figure \ref{idlow}).    

We now construct $P_l[0,0,k]$ where $l=h-k\ge 2$
as the union of four portions. 
Let us denote the first portion, which consists of two disjoint 
but identical parts, and which are to be called {\it sediments,} by 
$S_\alpha = Tr_{l-2}(P_{l-1})$ for $\alpha = R, L .$
The second portion is the base $B$ that is defined as   
$B= Tr_{2l}(P_{2l+k})$.   
The third portion, called the {\it core,} 
consists of only two triangles $ \triangle _1 \cup \triangle_2$
bordering a common edge. 
The last portion is called the {\it top,} and is defined to be 
a $P_{l-2}$.  
Then the construction is given by the union 
$$P_l[0, 0, k] = B \cup S_R \cup S_L \cup (\triangle _1 \cup \triangle_2)
\cup P_{l-2}.$$
Note that if $l=2$, then $S_R=S_L=P_1$ and 
$P_l[0,0,k]=B\cup S_R\cup S_L\cup (\triangle_1\cup\triangle_2)$.   

The upper boundary of the base $B$ is identified 
with the two lower boundaries of the sediments, 
but separated in the center by  
two boundary edges of $\triangle_1.$ 
Two of the three boundary edges of 
the top $P_{l-2}$ are identified with the upper 
boundaries of the sediments. Finally, the two remaining edges of 
the core are identified with the right-side edge of $S_L$ 
and with the left-side edge of $S_R$.       

We denote the subtriangulation consisting of the union of the last four 
pieces, that is, 
$S_R \cup S_L \cup (\triangle _1 \cup \triangle_2)\cup P_{l-2}$ 
by  
\begin{equation} \label{EqQs}
Q_l  = S_R \cup S_L \cup (\triangle _1 \cup \triangle_2)\cup P_{l-2}.
\end{equation}

In the case $l=h-k=1$, we construct $P_1[0,0,k]$ by taking the base
as $B=Tr_2(P_{2+k})$ and the top as a triangle $\triangle$ which  
connects the $3$ points on the upper boundary of $B$ so that 
$P_1[0,0,k]=B\cup \triangle$.  

Now let $l\ge 1$.
The {\it lower boundary } of $ P_l[0,0,k] $ is defined to be the lower 
boundary of its base $B.$  The rest of the boundary of $ P_l[0,0,k] $ 
is called the 
{\it upper boundary} of $P_l[0,0,k]$.   
\begin{figure}[htb]
\hspace {3cm} 
\xymatrix @M=0ex@R=3ex@C=2ex{
& & {\bullet}\ar@{-}[rr] \ar@{-}[dr]\ar@{-}[dl]& &{\bullet}\ar@{-}[dr] \ar@{-}[dl] \ar@{-}@{.}[rrr]& & 
	&{\bullet}\ar@{-}[dr] \ar@{-}@{.}[rrr]
	\ar@{-}[dl]& & &{\bullet}\ar@{-}[rr] \ar@{-}[dr] \ar@{-}[dl]& &{\bullet} \ar@{-}[dr] \ar@{-}[dl]& &  \\
& {\bullet} \ar@{-}[rr] & & {\bullet} \ar@{-}[rr] & &{\bullet}\ar@{-}@{.}[r] \ar@{-}@{.}[dd]& 
	{\bullet} \ar@{-}[rr] \ar@{-}[dr] \ar@{-}@{.}[ddl]	& &{\bullet}\ar@{-}[dl] \ar@{-}@{.}[r] 
		\ar@{-}@{.}[ddr]&{\bullet} \ar@{-}@{.}[dd]	\ar@{-}[rr] & & 	{\bullet}\ar@{-}[rr] & &{\bullet} &  \\
& & & & & & &{\bullet}\ar@{-}@{.}[d] & & & & & & &  \\
& {\bullet}\ar@{-}[rr] \ar@{-}[dr] \ar@{-}[dl]& & {\bullet}\ar@{-}[rr] \ar@{-}[dr] 
	\ar@{-}[dl]& &{\bullet} \ar@{-}[rr] \ar@{-}[dr] \ar@{-}[dl]& & {\bullet}\ar@{-}[rr] \ar@{-}[dr] 
	\ar@{-}[dl] & &{\bullet} \ar@{-}[rr] \ar@{-}[dr] \ar@{-}[dl]& & {\bullet}\ar@{-}[rr] \ar@{-}[dr] 
	\ar@{-}[dl] & &{\bullet} \ar@{-}[dr] \ar@{-}[dl] & \\
{\bullet}\ar@{-}[rr] & &{\bullet}\ar@{-}[rr] & & {\bullet} \ar@{-}[rr]& &{\bullet} \ar@{-}[rr]
	 && {\bullet}\ar@{-}[rr] & &{\bullet}\ar@{-}[rr] &  &{\bullet}\ar@{-}[rr] & &{\bullet} }
\end{figure}

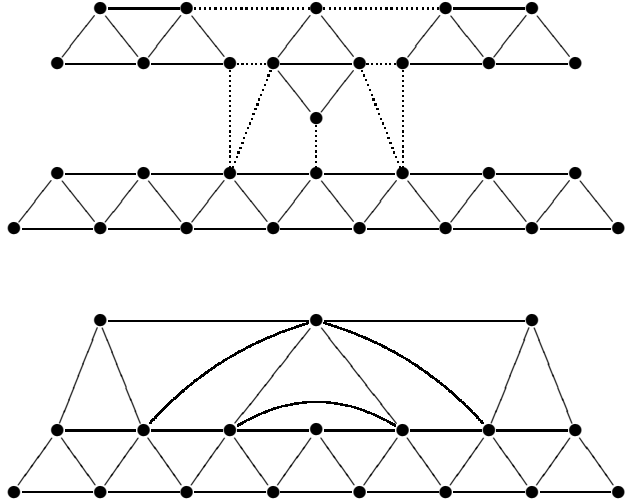
\begin{figure}[htb]
\hspace {3cm} 
\xymatrix @M=0ex@R=3.5ex@C=2ex{
& & & & & & & & & & & & & &  \\
& &{\bullet}\ar@{-}[rrrrr]\ar@{-}[ddr] \ar@{-}[ddl]& & & & &{\bullet} \ar@{-}[ddrr] \ar@{-}[ddll]
	\ar@/^0.5pc/@{-}[ddrrrr] \ar@/_0.5pc/@{-}[ddllll] \ar@{-}[rrrrr] 
	& & & & & {\bullet}\ar@{-}[ddr] \ar@{-}[ddl]& &  \\
& & & & & & & & & & & & & &  \\
& {\bullet}\ar@{-}[rr] \ar@{-}[dr] \ar@{-}[dl]& & {\bullet}\ar@{-}[rr] \ar@{-}[dr] 
	\ar@{-}[dl]& &{\bullet} \ar@{-}[rr] \ar@{-}[dr] \ar@{-}[dl] \ar@/^0.9pc/@{-}[rrrr] 
	& & {\bullet}\ar@{-}[rr] \ar@{-}[dr] 
	\ar@{-}[dl] & &{\bullet} \ar@{-}[rr] \ar@{-}[dr] \ar@{-}[dl]& & {\bullet}\ar@{-}[rr] \ar@{-}[dr] 
	\ar@{-}[dl] & &{\bullet} \ar@{-}[dr] \ar@{-}[dl] & \\
{\bullet}\ar@{-}[rr] & &{\bullet}\ar@{-}[rr] & & {\bullet} \ar@{-}[rr]& &{\bullet} \ar@{-}[rr]
	 && {\bullet}\ar@{-}[rr] & &{\bullet}\ar@{-}[rr] &  &{\bullet}\ar@{-}[rr] & &{\bullet} }
\caption{Identification of Base, Sediments and Core ($k=1$)}
\label{MmM2}
\end{figure}

It is not difficult to see that the length of the boundary of 
$P_l[0,0,k]$ is 
$$(2l+k)+2k+2+(l-2)=3(k+l)=3h.$$
We also note that the length of the upper boundary of  $ P_l[0,0,k] $ is 
$ h_{up} = 2k+l =2h-l$, and the length of the lower boundary of 
$ P_l[0,0,k] $ is $h_{low} = 2l+k = h+l$.   

It can be remarked that in this construction,
the integer $k$ has a geometric interpretation 
as the distance of the unique point of degree $4$ inside the patch 
$ P_l[0,0,k]$ to the lower boundary of the patch.    


In order to verify the expression for the number of points of degree six 
\begin{equation} \label{EMhk}
 N = N (h,k) = {h+k +2 \choose 2 } - { 2k+1 \choose 2} + { k \choose 2} - 3, 
\end{equation} 
we use induction on both $k$ and $h$.    
If $k=1$ and $h=2$,    
then we will have the patch $\{123, 145, 125, 256, 236, 367 \}$
with $N=4$. The right-hand side of 
(\ref{EMhk}) gives us 
${5 \choose 2}-{3\choose 2}
+{1\choose 2}-3=4$, which agrees with $N$.
Thus (\ref{EMhk}) (and so too (\ref{EET})) is true
for $(h,k)=(2,1)$.  

For the inductive step, we note that 
$P_{l-1}[0,0,k-1]$ can be embedded in 
$P_l[0,0,k]$ and that  
\begin{equation} \label{Egeo}
N(h,k) = N(h-2,k-1) + 3h,
\end{equation}
since the boundary of $P_l[0,0,k]$ consists of $3h$ points. 
On the other hand, we can compute using the right hand side of (\ref{EMhk}) 
that  
\begin{eqnarray}
 &  & {h+k +2 \choose 2 } - { 2k+1 \choose 2} + { k \choose 2} - 3 \\[2mm] 
\nonumber
 & = & {h+k-1 \choose 2 } - { 2k-1 \choose 2} + { k-1 \choose 2 } - 3 \\[2mm] 
\nonumber
 &  &  \qquad + (3h+3k)  ~~~ - (4k-1) ~~~+ (k-1) \\[3mm] \nonumber
 &= &  {h+k-1 \choose 2 } - { 2k-1 \choose 2} + { k-1 \choose 2 } - 3 + 3h. 
\end{eqnarray} 
This verifies equation (\ref{EMhk}). 

\begin{figure}[htb]
\hspace {3.5cm} 
\xymatrix @M=0ex@R=2.5ex@C=2ex{
& & & &{\bullet}\ar@{-}[rrrrrr] \ar@{-}[drrr] \ar@{-}@{.}[ddll]& & & & & &{\bullet} \ar@{-}[dlll]
	\ar@{-}@{.}[ddrr]& 
	& & &  \\
& & & & & & &{\bullet} \ar@{-}@{.}[d]& & & & & & &  \\
& &{\bullet}\ar@{-}[rrrrr]\ar@{-}[ddr] \ar@{-}[ddl]& & & & &{\bullet} \ar@{-}[ddrr] \ar@{-}[ddll]
	\ar@/^0.5pc/@{-}[ddrrrr] \ar@/_0.5pc/@{-}[ddllll] \ar@{-}[rrrrr] 
	& & & & & {\bullet}\ar@{-}[ddr] \ar@{-}[ddl]& &  \\
& & & & & & & & & & & & & &  \\
& {\bullet}\ar@{-}[rr] \ar@{-}[dr] \ar@{-}[dl]& & {\bullet}\ar@{-}[rr] \ar@{-}[dr] 
	\ar@{-}[dl]& &{\bullet} \ar@{-}[rr] \ar@{-}[dr] \ar@{-}[dl] \ar@/^0.9pc/@{-}[rrrr] 
	& & {\bullet}\ar@{-}[rr] \ar@{-}[dr] 
	\ar@{-}[dl] & &{\bullet} \ar@{-}[rr] \ar@{-}[dr] \ar@{-}[dl]& & {\bullet}\ar@{-}[rr] \ar@{-}[dr] 
	\ar@{-}[dl] & &{\bullet} \ar@{-}[dr] \ar@{-}[dl] & \\
{\bullet}\ar@{-}[rr] & &{\bullet}\ar@{-}[rr] & & {\bullet} \ar@{-}[rr]& &{\bullet} \ar@{-}[rr]
	 && {\bullet}\ar@{-}[rr] & &{\bullet}\ar@{-}[rr] &  &{\bullet}\ar@{-}[rr] & &{\bullet} }
\end{figure}

\vspace{0.8cm}

\begin{figure}[htb]
\hspace {3.5cm} 
\xymatrix @M=0ex@R=2.5ex@C=2ex{
& &{\bullet}\ar@{-}[rrrrr]\ar@{-}[ddr] \ar@{-}[ddl] \ar@/^1.8pc/@{-}[rrrrrrrrrr] 
	& & & & &{\bullet} \ar@{-}[ddrr] \ar@{-}[ddll]
	\ar@/^0.5pc/@{-}[ddrrrr] \ar@/_0.5pc/@{-}[ddllll] \ar@{-}[rrrrr] 
	& & & & & {\bullet}\ar@{-}[ddr] \ar@{-}[ddl]& &  \\
& & & & & & & & & & & & & &  \\
& {\bullet}\ar@{-}[rr] \ar@{-}[dr] \ar@{-}[dl] & & {\bullet}\ar@{-}[rr] \ar@{-}[dr] 
	\ar@{-}[dl]& &{\bullet} \ar@{-}[rr] \ar@{-}[dr] \ar@{-}[dl] \ar@/^0.9pc/@{-}[rrrr] 
	& & {\bullet}\ar@{-}[rr] \ar@{-}[dr] 
	\ar@{-}[dl] & &{\bullet} \ar@{-}[rr] \ar@{-}[dr] \ar@{-}[dl]& & {\bullet}\ar@{-}[rr] \ar@{-}[dr] 
	\ar@{-}[dl] & &{\bullet} \ar@{-}[dr] \ar@{-}[dl] & \\
{\bullet}\ar@{-}[rr] & &{\bullet}\ar@{-}[rr] & & {\bullet} \ar@{-}[rr]& &{\bullet} \ar@{-}[rr]
	 && {\bullet}\ar@{-}[rr] & &{\bullet}\ar@{-}[rr] &  &{\bullet}\ar@{-}[rr] & &{\bullet} }
\caption{Identification of Top, Sediments and Core ($k=1,l=3$)}
\label{TydddD}
\end{figure}
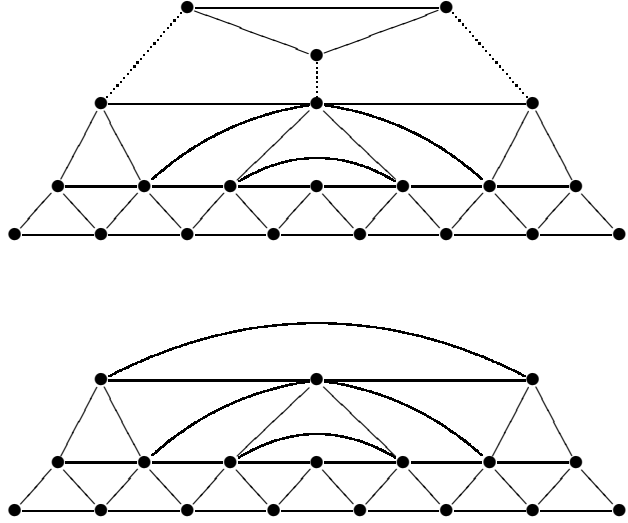

\vskip8mm
Assume that $ h > k \geq l > 0.$ Then for each $N$ of the form 
$$ N = {k+h+l+2 \choose 2}-3{ k+1 \choose 2} -3{ l+1 \choose 2}- 3,$$
it can be shown that there exists a patch with parameters 
$(0,3,0,N)_b$ with 
$ b=3h,\,\beta_4 = 1$ which is of the form $[0,l,k].$ 
For the construction we first assume that   
\begin{equation} \label{condhkl}
 l> \max \{ 2h-k-l, h+k-2l \} . \end{equation}
In this case we assemble the patch $ P[0,l,k]_{3h} $ out of the following 
five parts: 
$$ C=P[0,0,k-l]_{3(h-l)}; \qquad A=P[0,0,0]_{3(2h-k-l)}; $$
$$ B=P[0,0,0]_{3(h+k-2l)}; \qquad  Q_{l}; \qquad R[0,0,0]_{6l}. $$ 
Here $ R[0,0,0]_{6l} $ is given by cutting off two smaller triangles 
$ P[0,0,0]_{3(2l+k-2h)} $ and $P[0,0,0]_{3(3l-k-h)}$ from the large triangle 
$ P[0,0,0]_{6l} $ at two corners of 
$P[0,0,0]_{6l}$.
That is,  
\begin{equation}  \label{Rdef}
R[0,0,0]_{6l} := P[0,0,0]_{6l} \setminus ( P[0,0,0]_{3(2l+k-2h)}  
\cup P[0,0,0]_{3(3l-k-h)}).\qquad 
\end{equation}
The boundary of the triangulation $ R[0,0,0]_{6l} $ 
consists of altogether five edges, 
three of which are derived from 
the original boundary of the triangle $P[0,0,0]_{6l}$,    
and the remaining two are new boundary edges.
The part $Q_l$ is as defined in (\ref{EqQs}).    

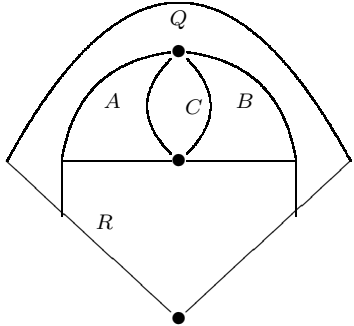
\begin{figure}[htb]
\hspace {4.5cm} 
\xymatrix @M=0ex@R=3.5ex@C=4ex{
& & & & & &\\
& & &{\bullet}\ar@/^1pc/@{-}[dd]_-{C} \ar@/_1pc/@{-}[dd] \ar@/_1pc/@{-}[ddll]^-{A} 
	\ar@/^1pc/@{-}[ddrr]_-{B}& & &\\
& & & & & &\\
& \ar@{-}[rr] \ar@{-}[d]& &{\bullet}\ar@{-}[rr] & &\ar@{-}[d] &\\
& & & & & &\\
& & & & & &\\
& & &{\bullet}\ar@{-}[uuulll]="b" _-{R}\ar@{-}[uuurrr]="a"& & & \ar@/_5pc/ @{-}"a";"b"^-{Q}}
\caption{Identification of five parts}
\label{ids5pc}
\end{figure}

The left-hand side of the boundary of $C=P[0,0,k-l]_{3(h-l)}$ 
is identified with 
one side of the triangle  $A=P[0,0,0]_{3(2h-k-l)}$; 
the right-hand side of the 
boundary of $C$ is identified with 
one side of the triangle $B=P[0,0,0]_{3(h+k-2l)}$.
The vertex of degree $4$,   
that is,  
one of the meeting points of the left and right-hand side boundaries 
of $C$ will be identified with the 
central point of degree $4$ of $Q_l$. This point is also identified with 
one of the corners of the triangles $A$ and $B$.  
Parts of the boundaries of $A$ and $B$   
are identified with the inner parts of the boundaries of $ Q_l$.
This leaves   
three parts of the boundary of $Q_l$ without identification, 
the top part which 
will be part of the boundary of the patch, 
and the two remaining inner parts $ Q \cap R$ 
(refer to Figure \ref{ids5pc}).
The central piece of the boundary of $ R[0,0,0]_{6l} $ 
which is opposite to the corner point of   
degree four is identified with the union 
of the remaining two boundaries of the 
triangles $A$ and $B$, namely,
 $ R \cap(A \cup B).$ 
The two new boundary parts 
of $R[0,0,0]_{6l}$ are identified with the 
remaining inner parts of the boundary  
of $Q_l$ as $ Q \cap R$ (refer to Figure \ref{ids5pc}) 
and are easily seen to have the corresponding lengths 
$2l+k-2h$ and $3l-k-h$.   

By using parts (ii) and (iii) of Theorem \ref{T1}, 
the number of points of degree $6$ of the 
parts $ C, A, B, Q, R$  can be seen to be as follows: 

\begin{eqnarray} \label{NC}
\label{NA}
N_A &=&  { 2h - k - l + 2 \choose 2 }-3 \\  \label{NB} 
N_B &=&  { h+k-2l +2 \choose 2 }-3 \\ \label{NC}
N_C &=&  { (h-l) +(k-l) + 2 \choose 2 } - { 2k-2l+1 \choose 2 } 
 + { k-l \choose 2}-3 \\ \label{NQ}
N_Q &=&  { 2l + 2 \choose 2 } - { 2l + 1 \choose 2 } + { l \choose 2 } -3 \\ 
\label{NR}
N_R &=&  { 2l+2 \choose 2 } - { 2l + k - 2h + 1 \choose 2 } 
 - { 3l-k-h+1 \choose 2 } -5
\end{eqnarray}
Note that the number of interior points of degree $6$ in $A, B$ are
$N_A-3(2h-k-l-1),\,
N_B-3(h+k-2l-1)$, respectively.  
The number of interior points of degree $6$ in $C$ is
$N_C-3(h-l-1)$ if $k=l$ and 
$N_C-(3h-3l-2)$ if $k>l$. 
The number of points of degree $6$ in $R$ except for those lying on the
common boundary with $Q,\,A,\,B$ is $N_R-2l+3$ whereas the number of points
of degree $6$ in $Q$ except for those lying on the common boundary
with $R,\,A,\,B$ is $N_Q-2(l-1)$. The number of points of degree $6$ 
lying on each of the common boundaries of $R,\,A,\,B,\,C,\,Q$
are as follows:
\begin{eqnarray*}
&& N_{(R\cap Q)_{l}} = 2l+k-2h+1\quad\mbox{(the left part of $R\cap Q$)},\\
&& N_{(R\cap Q)_{r}} = 3l-k-h+1\quad\mbox{(the right part of $R\cap Q$)},\\
&& N_{(A\cap R)\setminus (R\cap Q)_l} = 2h-k-l,\\
&& N_{(B\cap R)\setminus (A\cap R)\setminus (R\cap Q)_{r}} 
    = h+k-2l-1,\\
&& N_{(A\cap Q)\setminus (A\cap R)} 
    = 2h-k-l-1,\\
&& N_{(B\cap Q)\setminus (A\cap Q)\setminus (B\cap R)} 
    = h+k-2l-1,\\
&& N_{(A\cap C)\setminus (A\cap Q)\setminus (A\cap R)} 
 =\left\{ \begin{array}{ll}
   2h-k-l-2, &\mbox{if $k=l$}\\
   2h-k-l-1, &\mbox{if $k>l$}
\end{array}, \right.  
\\
&& N_{(B\cap C)\setminus (B\cap Q)\setminus (B\cap R)} 
   = h+k-2l-1. 
\end{eqnarray*}
Therefore the total number of points of degree $6$ 
lying on the common boundaries is
\begin{eqnarray*}
\bar N=
\left\{ \begin{array} {ll}
 6h-4l-4, &\mbox{if $k=l$}\\
 6h-4l-3, &\mbox{if $k>l$}
\end{array}
\right. 
.
\end{eqnarray*}
We also note that 
\begin{eqnarray*}
2 N_A &=& (2h-k-l+2)(2h-k-l+1)-6 \\
	 &=& 4h^2+k^2+l^2-4hk-4hl+2kl+6h-3k-3l-4, \\
2 N_B &=& (h+k-2l+2)(h+k-2l+1)-6 \\
      &=& h^2+k^2+4l^2+2hk-4hl-4kl+3h+3k-6l-4, \\
2 N_C &=& (h+k-2l+2)(h+k-2l+1)-(2k-2l+1)(2k-2l) \\
		& & +(k-l)(k-l-1)-6\\
   &=& h^2-2k^2+l^2+2hk-4hl+2kl+3h-3l-4, \\
2 N_Q &=& (2l+2)(2l+1)-2l(2l+1)+l(l-1)-6 \\
	&=& l^2+3l-4, \\
2N_R &=& (2l+2)(2l+1)-(2l+k-2h+1)(2l+k-2h)\\
	& & -(3l-k-h+1)(3l-k-h)-10\\
	&=& -5h^2-2k^2-9l^2+2hk+14hl+2kl+3h+l-8.  
\end{eqnarray*}
It follows that for $k=l$, the number of points of degree $6$ 
in this construction is
\begin{eqnarray*}
N &=& (N_A-3(2h-k-l-1))+(N_B-3(h+k-2l-1))\\
 &&\quad +(N_C-3(h-l-1))
  +(N_R-2l+3)\\
 &&\qquad +(N_Q-2(l-1))+\bar N\\
 &=& \frac 12(h^2-2k^2-2l^2+2hk+2hl+2kl+3h-4).
\end{eqnarray*}
If $k>l$, we also get the same $N$. 
Since 
\begin{eqnarray*}
& &{h+k+l +2 \choose 2 } - 3 { k+1 \choose 2} - 3{ l+1 \choose 2} - 3\\
 	&=&\frac{1}{2}[(h^2+hk+hl+h+hk+k^2+kl+k+hl+kl+l^2+l\\
& &  +2h+2k+2l+2)-(3k^2+3k)-(3l^2+3l)]-3  \\
	&=&\frac{1}{2}(h^2-2k^2-2l^2+2hk+2hl+2kl+3h-4),  
\end{eqnarray*}
we have the equality  
\begin{equation} \label{VERIFY}
N={h+k+l +2 \choose 2 } 
- 3 { k+1 \choose 2} - 3{ l+1 \choose 2} - 3 . 
\end{equation} 

The various cases of smaller $l,$ that is, those cases where (\ref{condhkl}) 
does not hold can be discussed in an analogous way.

\subsubsection{Some preliminaries for the proof of Theorem \ref{T1}}

We remark that for any patch $P$ with an interior point $q \in P$,   
the distance of $q$ to the boundary of $P$ can be defined as the minimal 
number of triangles that are intersected by a topological path connecting  
$q$ to the boundary of $P.$ 

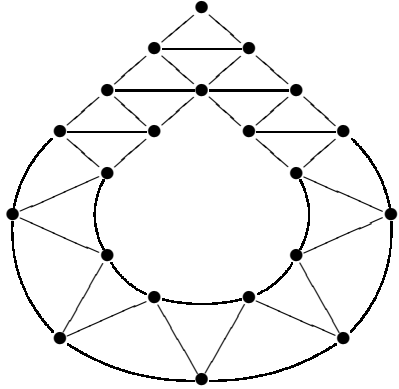
\begin{figure}[htb]
\hspace {3.5cm} 
\xymatrix @M=0ex@R=2ex@C=2.3ex{
& & & &{\bullet}\ar@{-}[dr] \ar@{-}[dl]& & & & \\
& & & {\bullet}\ar@{-}[rr]\ar@{-}[dr] \ar@{-}[dl]& &{\bullet}\ar@{-}[dr] \ar@{-}[dl] & & & \\
& & {\bullet}\ar@{-}[rr] \ar@{-}[dr] \ar@{-}[dl]="a"& &{\bullet} \ar@{-}[rr] \ar@{-}[dr] \ar@{-}[dl]& 
 	&{\bullet} \ar@{-}[dr]="b" \ar@{-}[dl]&  &\\
& {\bullet}\ar@{-}[rr] \ar@{-}[dr] \ar@/_0.3pc/@{-}[ddl]& &{\bullet} \ar@{-}[dl] & &{\bullet}\ar@{-}[rr] 
	\ar@{-}[dr] & & {\bullet}\ar@{-}[dl] \ar@/^0.3pc/@{-}[ddr] &\\
& & {\bullet} \ar@{-}[dll] \ar@/_0.4pc/@{-}[dd]& & & & {\bullet}\ar@{-}[drr] \ar@/^0.4pc/@{-}[dd]& & \\
{\bullet}\ar@{-}[drr] \ar@/_0.5pc/@{-}[dddr]& & & & & & & &{\bullet} \ar@{-}[dll]
	\ar@/^0.5pc/@{-}[dddl]\\
& &{\bullet} \ar@/_0.2pc/@{-}[dr] \ar@{-}[ddl]& & & & {\bullet}\ar@/^0.2pc/@{-}[dl] \ar@{-}[ddr]& & \\
& & &{\bullet} \ar@/_0.2pc/@{-}[rr] \ar@{-}[dll] \ar@{-}[ddr]& &{\bullet} \ar@{-}[ddl] \ar@{-}[drr]& & & \\
&{\bullet}\ar@/_0.4pc/@{-}[drrr] & && & & &{\bullet}\ar@/^0.4pc/@{-}[dlll] & \\
& & && {\bullet} & & & & }
\caption{A belt around patch $[o,l,k]$.}
\label{belth}
\end{figure}

If $P$ is a patch such that all its boundary points have degree $6,$ 
with the exception of one point which has degree $4,$ then  
we may consider the union of all triangles which are incident with the 
boundary of $P.$ This union of triangles, as depicted in Figure \ref{belth} 
will be referred to as a {\it belt } around the patch. 

The following is an immediate consequence of our definitions.

\begin{lem}  \label{Lbelt}
\begin{description} 
\item (i) If the boundary of the patch $P$ and the interior of all the 
triangles of a belt around $P$ are removed, the remaining topological space 
is another (triangulated) patch $P'$. 
\item (ii) If $q$ is an interior point of a patch $P$ 
which has distance $k$ to the 
boundary, then the distance of $q$ to the boundary of $P'$ is $k-1$. 
\end{description} 
\end{lem}

\subsection{Proof of Theorem \ref{T1}}

\subsubsection{Proof of part (i):}
This is shown by induction on the number $m= f_1(P)$ 
of points of the patch $P$.     
The theorem is true for $m=3,6$ since the patches with   
$m=3,6$ points where three of the points have degree $4$
are easily seen to be unique. 

Assume that $P$ is a patch of type $(0,3,0,m-3)_b$ with $m>6$ points. 
From the description 
of the process of peeling the boundary, we may assume that the boundary of the 
patch $P$ contains a vertex point $q$ of degree $4$. Hence the point 
$q$ has precisely two adjacent points $x,y$ 
which are both on the boundary of $P$.   
It is clear that $ d(x)=4 $ or $d(y) = 4$ 
implies that $m=3,$ which contradicts our assumption that
$m>6.$ Since the type of $P$ is $(0,3,0),$ it follows that $ d(x)=d(y)=6.$ 
Obviously $xy$ is an interior edge of $P,$ which is contained in the triangle 
$ qxy,$ and so it must be contained in a second triangle $xyz$. Since  
$d(x)=d(y)=6$, it follows that the edges $xz, yz$ are interior 
edges. There exist two further points $ x',y'$ such that $xx', yy'$ 
are boundary 
edges, which are then contained in the triangles $xx'z, yy'z$ respectively. 

If the degree of $z$ is $ d(z) = 4$, 
then clearly $ x'y'$ must be an edge on the boundary, and 
$x'y'z$ must be a triangle. This implies that $d(x')=d(y')=5$; 
a contradiction.  
Since $P$ has parameters $(0,3,0,m-3),$ it follows that $d(z)=6.$ 
This implies that 
$x'$ and $y'$ do not form an edge, and hence $x'y'z$ is not a triangle. 
\begin{figure}[htb]
\hspace {3.5cm} 
\xymatrix @M=0ex@R=3ex@C=3.5ex{
& & & & & &  &\\
& & & & {y'} \ar@{-}@{.}[ur] \ar@{-}@{.}[drr] \ar@{-}@{.}[ddr]& & &  \\
& & & & & & & \\
&  &y \ar@{-}[uurr]="a" \ar@{-}[rrr]_(.8){z} \ar@{-}[ddrr] \ar@{-}[ddll]& & &{\bullet} 
	\ar@{-}@{.}[drr]\ar@{-}@{.}[ddr]="b" \ar@{-}@{.}[ur]_-{\bf{P'}}& & \\
& & & & & &  & \\
q \ar@{-}[rrrr] & & & &x \ar@{-}[uur] \ar@{-}[rr] & & {x'} \ar@{-}@{.}[ur] 
\ar@{-}@{.}[r]& &\ar@/_2pc/ @{-}@{.}"a";"b"}
\caption{Cutting the corner at $q$.}
\label{pat030}
\end{figure}
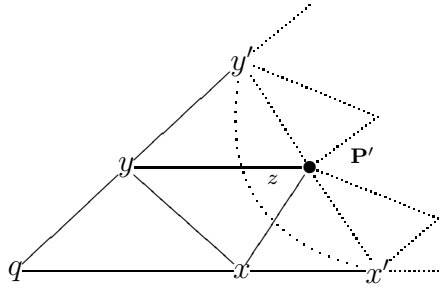

By cutting the corner at $q,$ as explained in Section \ref{pat030}
(see Figure \ref{cutcornq2}), we 
replace the four triangles $ qxy, xx'z,yy'z,xyz $ by a single triangle 
$ x'y'z, $ thereby removing the three points $ q,x,y $ and 
obtaining a patch $P'$ which has $ d(z) = 4.$ Note that 
in $P'$ the points $x',y'$ are boundary points, but $z$ is an interior point. 
Observe that the number of points, and the length of the boundaries of 
$P$ and $P'$ are related by the equations:
\begin{equation}  \label{EPdash}
b' = b-3, \qquad m' = m-3. 
\end{equation}

Since the number of vertex points of $P'$ is $ m'= m-3<m$,    
it follows from the hypothesis of the 
induction that the length of the boundary $b'$ 
of the patch $P'$ is a multiple 
of $3,$ that is, 
$ b'=3h'$ for some positive integer $h'$.
By (\ref{EPdash}) we see that the number 
$ b = 3 + 3h' = 3(h'+1) $      
is also a multiple of $3.$ This proves the first part of the theorem.

The existence of patches of the types mentioned in parts (ii), (iii) and
(iv) of Theorem \ref{T1} have already been shown in 
Section 2.2.1. We are thus left with the task of showing the necessity part. 
\medskip
 
\vspace{0.5cm}
\subsubsection{Proof of part (ii):} 
We prove the necessity part of part (ii) by induction 
on the boundary length $b=3h$.
The statement  (ii) is clearly true for $h=1,$ since the only patch 
with boundary length 
$3$ is a triangle. 

Assume that the statement (ii) is true for any $h' < h.$ 
Let $P$ be a patch with parameters $(0,3,0,N)_{3h}$ with $\beta_4 =3,$ so that the boundary length is $3h>1.$ 

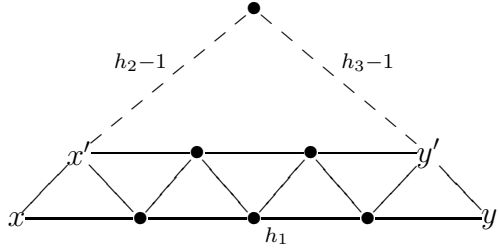
\begin{figure}[htb]
\hspace {3.5cm} 
\xymatrix @M=0ex@R=3ex@C=3ex{
&&&&{\bullet}\ar@{-}@{--}[dddrrr]^-{h_3-1} \ar@{-}@{--}[dddlll]_-{h_2-1}&&&&\\
&&&&&&&&\\
&&&&&&&&\\
& {x'} \ar@{-}[rr] \ar@{-}[dl] \ar@{-}[dr]& & {\bullet} \ar@{-}[rr] \ar@{-}[rr] \ar@{-}[dl] \ar@{-}[dr]
& & {\bullet} \ar@{-}[rr] \ar@{-}[rr] \ar@{-}[dl] \ar@{-}[dr]& &{y'}  \ar@{-}[dl] \ar@{-}[dr] &\\
{x} \ar@{-}[rr] & & {\bullet} \ar@{-}[rr]_(1.2){h_1}& & {\bullet} \ar@{-}[rr]& & {\bullet}\ar@{-}[rr]& &{y}}
\caption{A strip of  $2h_1-1$ triangles}
\label{striph}
\end{figure}

Take any pair $x,y $ of the  three boundary points of degree $4,$ and consider 
that part of the boundary between them which does not contain the third point 
of degree $4.$ Assume that $h_1$ is the length of this part of the boundary, 
which means that that part of the boundary has $h_1$ edges and $h_1+1$ vertex points.
The lengths of the other two parts of the boundary between the points of degree 
four are denoted by $h_2, h_3.$ Hence the overall length of the boundary of the 
given patch is $ b = h_1 + h_2 + h_3.$ 
By assumption, 
all inner points on this boundary 
(that is, the points on the boundary other than the three of degree $4$) 
as well as the points inside 
the patch at distance one from this boundary part must have degree $6.$ 
The segment $x'y'$ connecting the $h_1$ points of degree six 
which are of distance one from the boundary part of 
length $h_1$ separates the 
boundary from the rest of the triangle, and the area between the boundary 
and this segment consists of a strip of $2h_1-1$ triangles
(refer to Figure \ref{striph}).

By removing this strip from the patch, we obtain a smaller patch $P'$ which has 
three boundary segments of lengths $h_1-1, h_2-1, h_3-1$,   
and also contains three boundary points of degree $4$ on its boundary, 
namely, the third original boundary point, and the two end points $x',y'$.  
Hence the overall length of the boundary of $P'$ is (using part (i)) 
$$ 3h' = (h_1-1) + (h_2-1) + (h_3-1) = h_1+h_2+h_3-3 = b-3.$$ 
By the inductive assumption, we have that   
$$ h_1 - 1 = h_2-1 = h_3-1, $$ and hence 
$ h_1=h_2=h_3.$ We denote this common quantity by $ h= h_i.$ 

Again by the inductive assumption, we see that for the patch $P'$ 
with parameters
$(0,3,0,N')_{3h'}$, we get for $N'$ the equation
$$ N' = {(h-1)+2 \choose 2}-3.$$ 
Hence the total number of vertex points of degree $6$ in $P$ is 
$$ N = N' + (h_1+1) =  {(h-1)+2 \choose 2}-3 + h +1 = {h+2 \choose 2} -3. $$ 
This completes the proof of part (ii).

\vspace{0.5cm}
\subsubsection{Proof of part (iii):}

\begin{figure}[htb]
\hspace {3.5cm} 
\xymatrix @M=0ex@R=3ex@C=3ex{
&&&&{\bullet}\ar@{-}@{--}[dddrrr] \ar@{-}@{--}[dddlll]&&&&\\
&&&&&&&&\\
&&&&&&&&\\
& {x'} \ar@{-}[rr] \ar@{-}[dl] \ar@{-}[dr]& & {\bullet} \ar@{-}[rr] \ar@{-}[rr] \ar@{-}[dl] \ar@{-}[dr]
& & {\bullet} \ar@{-}[rr] \ar@{-}[rr] \ar@{-}[dl] \ar@{-}[dr]& &{y'}  \ar@{-}[dl] \ar@{-}[dr] &\\
{x} \ar@{-}[rr] & & {\bullet} \ar@{-}[rr]_(1.2){h_1}& & {\bullet} \ar@{-}[rr]& & {\bullet}\ar@{-}[rr]& &{y}}
\caption{A strip of triangles.}
\label{striph12}
\end{figure}
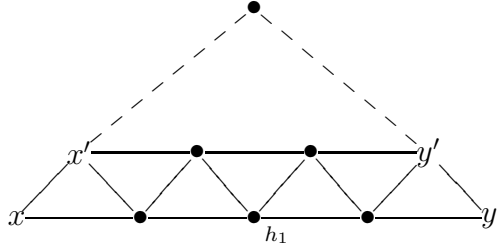

For the proof of (iii) we use induction on the integer $l>0$.   
We may assume throughout that $k>0$.   
That is, one of the points  
of degree four is not on the boundary of the patch.     
Let $x,y$ be the two boundary points of degree 
four, and assume that $ h_l, h_u $ are the lengths of the two boundary
parts between $x$ and $y,$ that is, the number of edges 
of the two parts of the boundary between $x$ and $y.$ 
The triangles which are 
next to the  boundary form a strip (refer to Figure \ref{striph12}) 
consisting   
of $ 2h_l-1$ or $2h_u-1$ triangles. 
After removing one of these strips from the patch, 
we get another patch of the kind $[0,0,k]$ or $[0,0,k-1]$,
depending on which strip is removed.   
We choose the notation in such a way that removing  
the boundary segment with $h_l$ edges 
gives rise to the remaining patch $P'$ of the form $[0,0,k]$.
In this case 
we get two new points $x', y' $ 
on the boundary of $P',$ which are of degree four. 
It is clear that removing the strip reduces the total boundary length by $3, $ 
hence we have $ b' = b-3, $ that is, $h' = h-1.$ 
From $k'=k$ it follows that $l'=l-1.$ 
Note that by the induction assumption of (iii), we get 
$ h_l' = 2h'-l' = 2h-2-l+1 = 2h-l-1.$ From Figure \ref{striph12} we easily infer that 
$ h_l = h_l'+1$, that is, 
we have $ h_l = 2h-l= h+k.$ Clearly, this implies 
that $ h_u = 3h- h_l = 
h+l$.   
Now we may use the inductive assumption for the formula on the number $N'$ to 
complete the proof, remembering that in the original patch, 
the points $x',\,y'$ have degree $6$ whereas the points $x,\,y$ have 
degree $4$. Then   
\begin{eqnarray*} 
N &=& (h_l+1) + N' \\
   &=& (h_l+1) + {h'+k'+2 \choose 2}-{ 2k'+1 \choose 2} + { k' \choose 2} - 3 \\
   &=& (h+k+1) + {(h-1)+k+2 \choose 2}-{ 2k+1 \choose 2} + { k \choose 2} - 3 \\    
    & = & {h+k+2 \choose 2}-{ 2k+1 \choose 2} + { k \choose 2} - 3.  
\end{eqnarray*}
This completes the proof of part (iii). 

\subsubsection{Proof of part (iv):} 
Again we use induction on $l>0$.   
Let $P$ be a patch with parameters $(0,3,0,N)_{3h}$ with $\beta_4=1$
and of the form $[0,\,l,\,k]$ where $0<l\le k<h$.  
If we remove all triangles which are incident with the boundary from $P$,
then we obtain a belt 
(refer to Figure \ref{belth}) consisting of $2b-3$ triangles. 
Removing this belt from $P$, we obtain another 
patch $P'$ of the form $[0,l-1,k-1]$
which is easily seen to have the boundary length $b'=b-3$.   


By the induction assumption we have that the number
of points of degree $6$ in $P'$ is  
$$ N' = {h'+k'+l'+2  \choose 2 } - 3 { k'+1 \choose 2} - 
3{ l'+ 1 \choose 2} - 3 . $$
Note that $ h'=h-1, l'=l-1, k'=k-1$.
Taking into consideration the fact that the boundary
of $P$ has one point of degree $4$ and
all other points of degree $6$,   
we obtain for the total number of points of degree $6$ in $P$ 
the expression 
\begin{eqnarray*}
N &=& (N'+1) + (b-1) \\
   &=& N' + 3h \\
   &=& 3h+ {h'+k'+l' +2 \choose 2 } - 3 { k'+1 \choose 2} - 3{ l'+1 \choose 2} - 3 \\
   &=&  3h+ {h+k+l -1 \choose 2 } - 3 { k \choose 2} - 3{ l \choose 2} - 3  \\
  &=& {h+k+l +2 \choose 2 } - 3 { k+1 \choose 2} - 3{ l+1 \choose 2} - 3 
\end{eqnarray*}
This completes the proof of (iv).

\section{Construction of certain patches \\
of types $(0,2,2),(0,1,4),(0,0,6)$ }

As a consequence of the above constructions, we can also construct 
certain patches of  types $(0,2,2), (0,1,4)$ and $(0,0,6).$ 
As defined earlier on,       
we let $ \beta_5 $ be the number of points of degree 
$5$ which are on the boundary. In this section we only consider the following 
kinds of patches:
\begin{description}
\item (i) $(0,2,2)$  with $\beta_5 = 2;$ 
\item (ii) $(0,1,4)$  with $\beta_5 = 4;$ 
\item (iii) $(0,0,6)$ with $\beta_5 = 6.$ \end{description}

\subsection{Type (1):~~$(0,2,2,N)_{b'}$ with $\beta_5 = 2 $ and $\beta_4 = 2.$ }
We can construct this type 
from  patches with parameters  $(0,3,0,\bar N)_{3h},$ 
with $ \beta_4 = 3,$ of type $P_h[0,0,0]$.   
Here we take a patch of type $ P_h[0,0,0],$ remove 
a subpatch of type $ P_{c-1}[0,0,0]$ ($2\le c\le h$)
which contains one of the 
boundary points of degree four, and take the topological closure 
to obtain a new patch $P'$ of type $(0,2,2)$.  
The edge-distance between the two 
boundary points of degree five in $P'$ is $c-1,$ 
and the edge-distance between a boundary point of degree 
four and  one of degree five is $ h-c+1$.  
The length of the boundary of $P'$ is 
$b'= (c-1)+2(h-c+1)+h= 3h-c+1$. 
By Theorem \ref{T1}(ii), the number of vertex points of degree
$6$ in $P_h[0,0,0]$
and $P_{c-1}[0,0,0]$ are 
$$\bar N= {h+2\choose 2}-3$$
and
$$N_{c-1}= {(c-1)+2\choose 2}-3={c+1\choose 2}-3,$$ 
respectively. Therefore the number of vertex points of degree $6$ in
$P'$ is 
$$N = \bar N-N_{c-1}+(c-2)-2
 ={h+2\choose 2}-{c\choose 2}-4.$$ 

\subsection{Type (2):~~$(0,2,2,N)_{b'}$ with $\beta_5 = 2$ and $\beta_4=1.$}
We can construct this type from the patch with parameters  
$(0,3,0,\bar N)_{3h},$
with $ \beta_4 = 2,$ of type $[0,0,k]$. 
Here we take a patch of type $ P_{h-k}[0,0,k],$ remove 
a subpatch of type $ P_{c-1}[0,0,0]$ ($c\le h-k+1$)
which contains one of the 
boundary points of degree four, and take the topological closure 
to obtain a new patch $P'$ of type $(0,2,2)$.   
Here $c-1$ is the number of edges on the boundary between the 
boundary point of degree four on
$ P_{h-k}[0,0,k]$ and (either of) the boundary 
points of degree five of the  
new patch $P'$.   
The length of the boundary of $P'$ is 
$ b' = 3h-c+1.$ The edge-distance between the two 
boundary points of degree five is $c-1$.  
By Theorem \ref{T1}, the number of vertex points of degree $6$
in $P_{h-k}[0,0,k]$ and 
$P_{c-1}[0,0,0]$ are 
$$\bar N={h+k+2\choose 2}-{2k+1\choose 2}+{k\choose 2}-3$$
and
$$N_{c-1}={(c-1)+2\choose 2}-3
  ={c+1\choose 2}-3,$$
respectively. Therefore the number of vertex points
of degree $6$ in $P'$ is
\begin{eqnarray*}
N &=& \bar N-N_{c-1}+(c-2)-2\\
 &=&{h+k+2\choose 2}-{2k+1\choose 2}+{k\choose 2}-{c\choose 2}-4.
\end{eqnarray*}
\medskip 

\subsection{Type (3):~~$(0,2,2,N)_{b'}$ with $\beta_5 = 2 $ and $\beta_4=0.$}
We can construct this type from the patch with parameters 
$(0,3,0,\bar N)_{3h},$ 
with $ \beta_4 = 1,$ of type $[0,l,k]$. 
Here we take a patch of type $ P[0,l,k]_{3h},$ remove 
a subpatch of type $P_{c-1}[0,0,0]$ ($c\le h+1$)
which contains 
one of the 
boundary points of degree four, and take the topological closure 
to obtain a new patch $P'$ of type $(0,2,2)$.  
Here $c-1$ is the number of edges on the boundary between the 
boundary point of degree four on
$ P[0,l,k]_{3h} $ and (either of) the boundary 
points of degree five of the  
new patch $P'$.    
The length of the boundary of $P'$ is 
$ b' = 3h-c+1.$ 
The edge-distance between the two 
boundary points of degree five is $c-1.$ 
By Theorem \ref{T1}, the number of vertex points of degree $6$
in $P[0,l,k]_{3h}$ and
$P_{c-1}[0,0,0]$ are 
$$\bar N={h+k+l+2\choose 2}-3{ k+1\choose 2}
  -3{l+1\choose 2}-3$$
and
$$N_{c-1}={c+1\choose 2}-3,$$
respectively. Therefore the number of vertex points of degree $6$
in $P'$ is  
\begin{eqnarray*}
N &=& \bar N-N_{c-1}+(c-2)-2\\
 &=& {h+k+l+2\choose 2} -3{k+1\choose 2}
     -3{l+1\choose 2} -{c\choose 2}-4.
\end{eqnarray*}
\medskip

\subsection{Type (4):~~$(0,1,4,N)_{b'}$ with $\beta_5 = 4 $ and $\beta_4 = 1.$}
We can construct this type 
from  patches with parameters  $(0,3,0,\bar N)_{3h},$ 
with $ \beta_4 = 3,$ of type $P_h[0,0,0]$. 
Here we take a patch of type $ P_h[0,0,0],$ remove 
two subpatches of types $P_{c_1-1}[0,0,0], P_{c_2-1}[0,0,0]$,  
each of which contains a boundary point of $P_h[0,0,0]$ of degree four, 
and take the topological closure to obtain a new patch $P'$ of type
$(0,1,4)$. Clearly we must have 
$(c_1-1)+(c_2-1)<h$.    
By Theorem \ref{T1}(ii), the number of vertex points of degree $6$
in $P_h[0,0,0],\,P_{c_1-1}[0,0,0]$ and
$P_{c_2-1}[0,0,0]$ are 
$$\bar N={h+2\choose 2}-3,$$ 
$$N_{c_1-1}={(c_1-1)+2\choose 2}-3={c_1+1\choose 2}-3$$
and
$$N_{c_2-1}={(c_2-1)+2\choose 2}-3={c_2+1\choose 2}-3,$$
respectively. 
Therefore the number of vertex points of degree $6$
in $P'$ is  
\begin{eqnarray*}
N &=& \bar N-N_{c_1-1}-N_{c_2-1}+(c_1-2)-2+(c_2-2)-2\\
 &=& {h+2\choose 2} -{c_1\choose 2}
     -{c_2\choose 2} -5.
\end{eqnarray*}
\medskip

Note that in the case $ c_1+c_2-2 = h $, there is another kind 
of patch of type $(0,2,2).$ This will be discussed below as 
type (7).

\subsection{Type (5):~~$(0,1,4,N)_{b'}$ with $ \beta_5 = 4 $ and $\beta_4 = 0.$ }
We can construct this type 
from  patches with parameters  $(0,3,0,\bar N)_{3h},$ 
with $ \beta_4 = 2,$ of type $[0,0,k].$ 
Here we take a patch of type $ P_{h-k}[0,0,k],$ remove 
two subpatches of types $ P_{c_1-1}[0,0,0]$ and $P_{c_2-1}[0,0,0]$,   
each of which contains one of the boundary points of $P_{h-k}[0,0,k]$
of degree four,     
and take the topological closure to obtain a new patch
$P'$ of type $(0,1,4)$.  Clearly we must have  
$(c_1-1)+(c_2-1)<h$.    
By Theorem \ref{T1}, the number of vertex points of degree $6$
in $P_{h-k}[0,0,k],\,P_{c_1-1}[0,0,0]$ and
$P_{c_2-1}[0,0,0]$ are 
$$\bar N={h+k+2\choose 2}-{2k+1\choose 2}+{k\choose 2}-3,$$ 
$$N_{c_1-1}={(c_1-1)+2\choose 2}-3={c_1+1\choose 2}-3$$

\noindent
and
$$N_{c_2-1}={(c_2-1)+2\choose 2}-3={c_2+1\choose 2}-3,$$
respectively. 
Therefore the number of vertex points of degree $6$
in $P'$ is  
\begin{eqnarray*}
N &=& \bar N-N_{c_1-1}-N_{c_2-1}+(c_1-2)-2+(c_2-2)-2\\
 &=& {h+k+2\choose 2} -{2k+1\choose 2}+{k\choose  2}-{c_1\choose 2}
     -{c_2\choose 2} -5.
\end{eqnarray*}
\medskip


\subsection{Type (6):~~$(0,0,6,N)_{b'}$ with $ \beta_5 = 6 $ and $\beta_4 = 0.$}
We can construct this type 
from  patches with parameters  $(0,3,0,\bar N)_{3h},$ 
with $ \beta_4 = 3,$ of type $P_h[0,0,0]$. 
Here we take a patch of type $P_h[0,0,0],$ remove 
three subpatches of type $ P_{c_1-1}[0,0,0], P_{c_2-1}[0,0,0], 
P_{c_3-1}[0,0,0]$,  
each of which contains one of the boundary points $P_h[0,0,0,]$
of degree four, 
and take the topological closure to obtain a new patch 
$P'$ of type $(0,6,6)$.   
Clearly we must have $c_1+c_2-2<h, c_1+c_3-2 < h, c_2+c_3-2<h.$ 
By Theorem \ref{T1}, the number of vertex points of degree $6$
in $P_h[0,0,k],\,P_{c_1-1}[0,0,0]$, 
$P_{c_2-1}[0,0,0]$ and  
$P_{c_3-1}[0,0,0]$ are 
$$\bar N={h+2\choose 2}-3,$$ 
$$N_{c_1-1}={(c_1-1)+2\choose 2}-3={c_1+1\choose 2}-3,$$
$$N_{c_2-1}={(c_2-1)+2\choose 2}-3={c_2+1\choose 2}-3,$$
and
$$N_{c_3-1}={(c_3-1)+2\choose 2}-3={c_3+1\choose 2}-3,$$
respectively. 
Therefore the number of vertex points of degree $6$
in $P'$ is  
\begin{eqnarray*}
N &=& \bar N-N_{c_1-1}-N_{c_2-1}-N_{c_3-1}+(c_1-2)-2+(c_2-2)-2\\
  &&\qquad +(c_3-2)-2\\
 &=& {h+2\choose 2} -{c_1\choose 2}
     -{c_2\choose 2} -{c_3\choose 2}-6.
\end{eqnarray*}
\medskip


\subsection{Type (7):~~$(0,2,2,N)_{b'}$ with $\beta_5 = 2 $ and $\beta_4 = 2.$ }
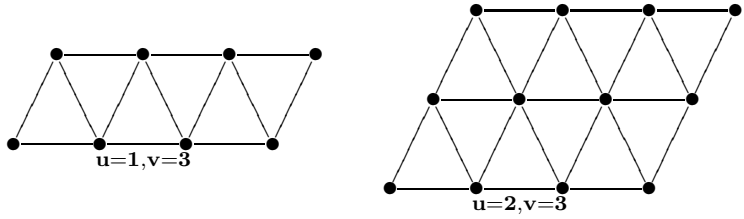
\begin{figure}[htb]
\[ \hspace {1.5cm} 
\vcenter 
{\xymatrix @M=0ex@R=5.5ex@C=2ex{
&{\bullet} \ar@{-}[rr] \ar@{-}[dr] \ar@{-}[dl] &&{\bullet} \ar@{-}[rr] \ar@{-}[dr] \ar@{-}[dl]&
	&{\bullet}\ar@{-}[rr] \ar@{-}[dr] \ar@{-}[dl]&&{\bullet} \ar@{-}[dl] \\
{\bullet}\ar@{-}[rr]&&{\bullet}\ar@{-}[rr]_-{\bf{u=1,v=3}}&&{\bullet}\ar@{-}[rr]&&{\bullet}&  }}
\quad 
\vcenter
{\xymatrix @M=0ex@R=5.5ex@C=2ex{
 && &{\bullet} \ar@{-}[rr] \ar@{-}[dr] \ar@{-}[dl] &&{\bullet} \ar@{-}[rr] \ar@{-}[dr] \ar@{-}[dl]&
	&{\bullet}\ar@{-}[rr] \ar@{-}[dr] \ar@{-}[dl]&&{\bullet} \ar@{-}[dl] \\
& &{\bullet} \ar@{-}[rr] \ar@{-}[dr] \ar@{-}[dl] &&{\bullet} \ar@{-}[rr] \ar@{-}[dr] \ar@{-}[dl]&
	&{\bullet}\ar@{-}[rr] \ar@{-}[dr] \ar@{-}[dl]&&{\bullet} \ar@{-}[dl] & \\
& {\bullet}\ar@{-}[rr]&&{\bullet}\ar@{-}[rr]_-{\bf{u=2,v=3}}&&{\bullet}\ar@{-}[rr]&&{\bullet}&  }}
\]
\caption{Patch with corner points $4,5,4,5$}
\label{rstrip}
\end{figure}

In contrast to type (1), the boundary points of degree $4$ and $5$ 
are arranged in cyclic order at the corners as $4,5,4,5$.    \vskip3mm

This type can be constructed by glueing $u$ {\it reduced strips } 
(refer to Figure \ref{rstrip}) consisting of  $2v$ triangles.

\section{Some applications towards the classification of patches}

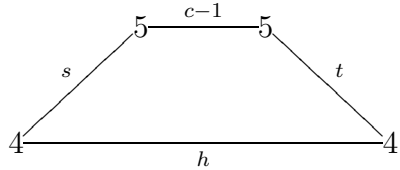
\begin{figure}[htb]
\hspace {4cm} 
\xymatrix @M=0ex@R=7ex@C=8ex{
& 5 \ar@{-}[dl]_-{s}  \ar@{-}[r]^-{c-1}  & 5 \ar@{-}[dr]^-{t} & \\
	4 \ar@{-}[r] &  \ar@{-}[r]_-{h} & \ar@{-}[r] & 4}
\caption{Patch with corner points $4,4,5,5$.}
\label{bd4455}
\end{figure}

We first consider patches of Type (1), that is, $(0,2,2,m)_{b'}$
with $\beta_5=2$ and $\beta_4=2.$
Let us additionally assume that the four boundary points 
of degree less than six are arranged in cyclic order 
such that the two points of degree $4$ precede the 
two points of degree $5$ (refer to Figure \ref{bd4455}). 
    
Let us choose positive integers 
$s,t,c,h$ such that the two points of degree $5$  
are the end points of a boundary string of $c-1$ edges, 
and such that the two points of degree $4$ are the 
end points of a string of $h$ edges. Let us assume that 
the remaining two pieces of the boundary between the 
two pairs of points of degrees $4$ and $5$ have lengths 
$s$ and $t$ (refer to Figure \ref{bd4455}). 

\begin{thm} \label{T2}
Under the assumptions on a patch $P$ with arrangement
of boundary points $4,4,5,5,$ and lengths as in Figure \ref{bd4455}
(as stated in the previous paragraph), we have the following: 
\begin{description}
\item{(i)} $ s=t; $
\item{(ii)} $ c-1 + s =h;$  
\item{(iii)} $ h \geq c;$ 
\item{(iv)} The number $N$ of vertex points of degree $6$ of the patch $P$ 
is given as 
\end{description}  
\begin{equation} \label{ER030}
N = { h+2 \choose 2 } - { c \choose 2 } - 4. \end{equation} 
\end{thm}

Proof: We take the union of the patch $P$ and another 
patch $Q$ of type $P_{c-1}[0,0,0]$ in such a way that 
one of the boundary segments between two points of degree 
$4$ of $Q$ is identified with the boundary segment of $P$ 
that lies between the two points of degree $5.$ 

The result of this identification is one patch $R = P \cup Q $ which has 
precisely three boundary points of degree $4$ and no other 
vertex points of degree different from $6.$ 
From the construction of $R$ we see that the three boundary 
segments between the boundary points of degree $4$ have lengths 
$ s+c-1, t+c-1$ and $h,$ respectively.  

Now we may apply part (ii) of Theorem \ref{T1} to the patch  $R.$ 
It follows that 
the lengths of the three boundary segments 
between the points of degree $4$ are all equal, 
which means that $ h = s+c-1 = t+c-1.$ 
This proves (i) and (ii). Part (iii) follows 
since $s,t >0.$ 
In the identification of $R= P \cup Q $ we observe that the two boundary points 
of degree $5$ of $P$ and two of the three boundary points of degree $4$ of $Q$ 
are identified, and become in $R$ boundary points of degree $6.$ By taking 
proper account of this fact we see 
that for the numbers of points of 
degree $6$ of the relevant patches,   
\begin{equation} \label{ERF}
N_R = N_P + N_Q - (c-2)+2. \end{equation}
By (\ref{E0}) we get 
$$ N_R = { h+2 \choose 2} - 3, $$
$$ N_Q = { (c-1)+2 \choose 2 } - 3.$$
Hence from (\ref{ERF}) 
\begin{eqnarray*}
 N_P &=& N_R - N_Q + c -4 \\
    &=& [{ h+2 \choose 2 } - 3] - [ {c+1 \choose 2} - 3] + c - 4 \\
    &=& { h+2 \choose 2} - { c \choose 2} - 4, 
\end{eqnarray*}
and the formula (\ref{ER030}) follows. This completes the proof of 
Theorem \ref{T2}.

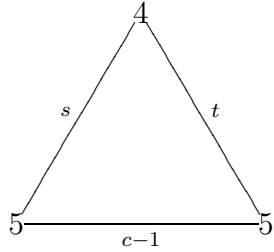
\begin{figure}[htb]
\hspace {5cm} 
\xymatrix @M=0ex@R=14ex@C=8ex{
& 4 \ar@{-}[dl]_-{s}  \ar@{-}[dr]^-{t} & \\
	5 \ar@{-}[r] _>{c-1} &  \ar@{-}[r] & 5}
\caption{Patch with boundary points $4,5,5$}
\label{bd455}
\end{figure}

\vspace{0.5cm}
In the case of Type (2), that is,  
$(0,2,2,m)_{b'}$ with $ \beta_5 = 2 $ and $\beta_4 = 1,$ 
let us choose positive integers 
$s,t,c$ such that the two points of degree $5$  
are the end points of a boundary string of $c-1$ edges 
and such that the boundary point of degree $4$ 
is one end point of two strings of $s,t$ edges with 
second endpoint of degree $5$ (refer to Figure \ref{bd455}). 

\begin{thm} \label{T3}
Under the assumptions on a patch $P$ of type $(0,2,2,m)_{b'}$
with the lengths of the boundary segments as in Figure \ref{bd455}
(as stated in the previous paragraph), we have the following: 
\begin{description}
\item{(i)} $ s+t+2(c-1),  2s-t + (c-1), 2t-s + (c-1)$ are multiples of $3;$ 
\item{(ii)} There exist integers $h,\,k,\,l$ such that    
\begin{eqnarray} \label{hndef}  
h &=& \frac{1}{3}(s+t+2(c-1)), \\
\label{kndef}
k &=& \frac{1}{3}(2s-t+(c-1)), \\
\label{lndef}
l &=& \frac{1}{3}(2t-s+(c-1)), 
\end{eqnarray}
and such that the patch $P$ has precisely 
\begin{equation} \label{Nndef}
N = { h+k+2 \choose 2} - { 2k+1 \choose 2} + {k \choose 2} - { c \choose 2} -4 
\end{equation} 
vertex points of degree $6.$  
\end{description}    \end{thm}

Proof: We take the union of the patch $P$ and another 
patch $Q$ of type $P_{c-1}[0,0,0]$ in such a way that 
one of the boundary segments between two points of degree 
$4$ of $Q$ is identified with the boundary segment of $P$ 
that lies between the two points of degree $5.$ 

The result of this identification is one patch $R = P \cup Q $ which has 
precisely two boundary points of degree $4,$ and it  contains precisely 
one point of degree $4$ in its interior, since by assumption $P \subset R$ 
is of type $(0,2,2,m).$ 
Hence $R$ is a patch of type $(0,3,0)$ with $ \beta_4=2.$ 
We see that the total length of the boundary of $R$ is 
$ s+t+2(c-1).$ 
By Theorem \ref{T1}(i) it follows that this boundary length
is a multiple of $3.$ From 
\begin{eqnarray*}
2t-s + (c-1) &=& 3(t + c-1) - (s+t+2(c-1)), \\
2s-t +(c-1)  &=& 3(s+c-1)   -(s+t+2(c-1)), \end{eqnarray*}
it follows that the other two  integers $2t-s+(c-1),2s-t+(c-1)$ are also multiples of $3.$ 

In the identification of $R= P \cup Q $ we observe that the two boundary points 
of degree $5$ of $P$ and two of the three boundary points of degree $4$ of $Q$ 
are identified, and become in $R$ boundary points of degree $6.$ By taking 
proper account of this fact we see 
that for the numbers of points of 
degree $6$ of the relevant patches,   
\begin{equation} \label{ERG}
N_R = N_P + N_Q - (c-2)+2. \end{equation}
By referring to Figure \ref{bd455} we see that the lengths of the two 
boundary segments between the two  boundary points of degree $4$ are 
given as $ s+c-1$ and $t+c-1.$ By using part (iii) of 
Theorem \ref{T1} we may take    
\begin{eqnarray*}
 h+l &=& t+c-1, \\
 2h-l &=& s+c-1. 
\end{eqnarray*}
By adding these two equations we see that 
$ 3h = s+t+2(c-1)$. By subtracting the second equation
from twice the first 
we have that 
$ 3l = 2t-s + c-1.$ 
From the fact that $k+l=h$ (by Theorem \ref{T1}(iii)), we have  
\begin{equation} \label{hrec}
k = \frac{1}{3} ( 2s-t + c-1).
\end{equation}
By (\ref{EET}) we get 
$$ N_R = { h+k+2 \choose 2} - {2k+1 \choose 2} + { k \choose 2} - 3, $$
and by (\ref{E0}) we have  
$$ N_Q = { (c-1)+2 \choose 2 } - 3.$$
Hence from (\ref{ERG}) 
\begin{eqnarray*}
 N_P &=& N_R - N_Q + c -4 \\
    &=& [{ h+k+2 \choose 2 } - {2k+1 \choose 2} + { k \choose 2} - 3] \\ 
   &     & - [{c+1 \choose 2} - 3] + c - 4 \\
    &=& { h+k+2 \choose 2} -{2k+1 \choose 2} + {k \choose 2}- { c \choose 2} - 4, 
\end{eqnarray*}
and the formula (\ref{Nndef}) follows. 
This completes the proof of Theorem \ref{T3}. \vskip3mm

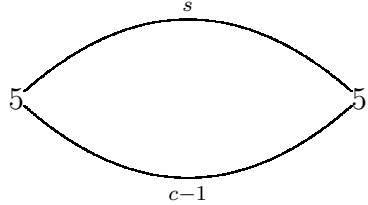
\begin{figure}[htb]
\hspace {4cm} 
\xymatrix @M=0ex@R=12ex@C=12ex{
5 \ar@/^2.5pc/@{-}[rr]^-{s}  & & 5 \ar@/^2.5pc/@{-}[ll]^-{c-1} }
\caption{Boundary of a patch with two boundary points of degree $5$}
\label{bd55}
\end{figure}

In the case of Type (3), that is,  
$(0,2,2,m)_{b'}$ with $\beta_5 = 2$ and $\beta_4 = 0,$ 
let us choose positive integers 
$s,c$ such that the two points of degree $5$  
are the end points of two boundary strings with $c-1$ 
and $s$ edges (refer to Figure \ref{bd55}). 

\begin{thm} \label{T4}
Under the assumptions on a patch $P$ of type $(0,2,2,m)_{b'}$
with the lengths of the boundary segments as in Figure \ref{bd55}
(as stated in the previous paragraph), we have the following: 
\begin{description}
\item{(i)} $ s+2(c-1), 2s +(c-1)$ are  multiples of $3;$ 
\item{(ii)} There exist integers $h,\,k,\,l$ such that 
\begin{equation} \label{klndef}  
h = \frac{1}{3}(s+2(c-1)) 
\end{equation}
and such that the patch $P$ has precisely 
\begin{equation} \label{Nkldef}
N = { h+k+l+2 \choose 2} - 3{ k+l \choose 2} - 3{l+1 \choose 2} - {c \choose 2} -4 
\end{equation} 
vertex points of degree $6.$  
\end{description}    \end{thm}

Proof: We take the union of the patch $P$ and another 
patch $Q$ of type $P_{c-1}[0,0,0]$ in such a way that 
one of the boundary segments between two points of degree 
$4$ of $Q$ is identified with the boundary segment of $P$ 
that lies between the two points of degree $5$ that contains  
$c-1$ edges.

The result of this identification is one patch $R = P \cup Q $ which has 
precisely one boundary point of degree $4,$ and it  contains precisely 
two points of degree $4$ in its interior, 
since by assumption $P \subset R$ 
is of type $(0,2,2,m)$.   
Hence $R$ is a patch of type $(0,3,0)$ with $ \beta_4=1.$ 
We see that the total length of the boundary of $R$ is 
$ s+2(c-1)$.   
By Theorem \ref{T1}(i) it follows that this boundary length is 
a multiple of $3$. That is,
$s+2(c-1)=3h$ for some integer $h$. From 
$$ 2s + (c-1) = 3(s+c-1) - (s+2(c-1)),$$   
it follows that $2s+(c-1)$ is also a multiple of $3$.   

In the identification of $R= P \cup Q $ we observe that the two boundary points 
of degree $5$ of $P$ and two of the three boundary points of degree $4$ of $Q$ 
are identified, and become in $R$ boundary points of degree $6.$ By taking 
proper account of this fact we see that 
for the numbers of points of 
degree $6$ of the relevant patches,   
\begin{equation} \label{ERH}
N_R = N_P + N_Q - (c-2)+2. \end{equation}
By part (iv) of Theorem \ref{T1} we see that 
there exist integers $k,l$ where $0<l\le k<h$ such that    
$$ N_R = { h+k+l+2 \choose 2} - 3{k+l \choose 2} - 3{ l+1 \choose 2} - 3. $$
By Theorem \ref{T1}(ii),      
$$ N_Q = { (c-1)+2 \choose 2 } - 3.$$
Hence from (\ref{ERH}) 
\begin{eqnarray*}
 N_P &=& N_R - N_Q + c -4 \\
    &=& [{ h+k+l+2 \choose 2 } - 3{k+l \choose 2} - 3{ l+1 \choose 2} - 3] - 
[{c+1 \choose 2} - 3] \\
 &&\qquad + c - 4 \\
    &=& { h+k+l+2 \choose 2} - 3{k+l \choose 2} - 3{l+1 \choose 2}- { c \choose 2} - 4, 
\end{eqnarray*}
and the formula (\ref{Nkldef}) follows. 
This completes the proof of Theorem \ref{T4}. \vskip5mm
  
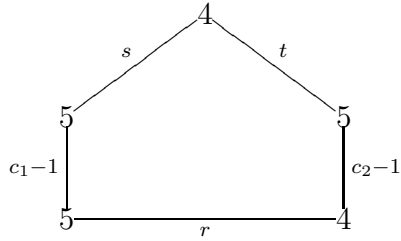
\begin{figure}[htb]
\hspace {4cm} 
\xymatrix @M=0ex@R=6ex@C=9ex{
& 4 \ar@{-}[dr]^-{t}  \ar@{-}[dl]_-{s} & \\
5 \ar@{-}[d]_-{c_1-1}  & & 5 \ar@{-}[d]^-{c_2-1} \\
5 \ar@{-}[r] &  \ar@{-}[r]_<{r} &  4}
\caption{Boundary of a patch with four boundary points of degree $5$, one of degree $4$.}
\label{bd45555}
\end{figure}

In the case of Type (4), that is,  
$(0,1,4,m)_{b'}$ with $\beta_5 = 4$ and $\beta_4 = 1,$ 
let us choose positive integers 
$r,s,t,c_1,c_2$ such that the four points of degree $5$  
are the end points of three boundary strings with $c_1-1,r,c_2-1$ 
edges in this cyclic order respectively (refer to Figure \ref{bd45555})
and such that the two segments of the boundary containing the boundary 
point of degree $4$ have $s$ and $t$ edges respectively 
(refer to Figure \ref{bd45555}). 

\begin{thm} \label{T5}
Under the assumptions on a patch $P$ of type $(0,1,4,m)_{b'}$
with the lengths of the boundary segments as in Figure \ref{bd45555}
(as stated in the previous paragraph), we have the following: 
\begin{description}
\item{(i)} The integer $ r+s+t+2(c_1-1)+2(c_2-1)$ 
is a  multiple of $3;$ 
\item{(ii)} The equality
\begin{equation} \label{rstdef}  
s+c_1-1 = t+c_2-1 = r + c_1+c_2-2
\end{equation} holds;
\item{(iii)} The patch $P$ has precisely 
\begin{equation} \label{Nrstdef}
N = { h+2 \choose 2} - { c_1 \choose 2} - {c_2 \choose 2 }  -5 
\end{equation} 
vertex points of degree $6$ where $h=s+c_1-1$.    
\end{description}    \end{thm}

Proof: We take the union of the patch $P$ and two 
more patches $Q_1, Q_2$ of types 
$P_{c_1-1}[0,0,0],\,P_{c_2-1}[0,0,0]$ 
respectively in such a way that one of the boundary segments 
between two points of degree 
$4$ of $Q_i$ is identified with the boundary segment of $P$ 
that lies between the two points of degree $5$ which contain  
$c_i-1$ edges ($i=1,2$).

The result of this identification is one patch 
$R = P \cup Q_1 \cup Q_2 $ which has 
precisely three boundary points of degree $4.$ 
Hence $R$ is a patch of type $(0,3,0)$ with $ \beta_4=3.$ 
We see that the total length of the boundary of $R$ is 
$ r+s+t+2(c_1-1)+2(c_2-1).$ 
By Theorem \ref{T1}(i) it follows that this boundary length
is a multiple of $3.$  

In the identification of $R= P \cup Q_1 \cup Q_2 $ 
we observe that the four boundary points 
of degree $5$ of $P$ and two of 
the three boundary points of degree $4$ of $Q_1, Q_2$ 
each are identified, and become in $R$ boundary points of degree $6$.
By taking 
proper account of this fact we see that 
for the numbers of points of 
degree $6$ of the relevant patches,   
\begin{equation} \label{ERL}
N_R = N_P + N_{Q_1} +N_{Q_2} - (c_1-2)-(c_2-2)+4. \end{equation}
By Theorem \ref{T1}(ii) we see that 
all three boundary parts of $R$ between the 
boundary points of degree $4$ are equal, that is,  
$$s+c_1-1 = t+c_2-1 = r + c_1-1+c_2-1=h,$$
which proves part (ii). 
 
For the number of vertex points of degree $6$ of the patch $R$ we see 
from Theorem \ref{T1}(ii) that 
$$ N_R = { h+2 \choose 2} - 3. $$
We also have by Theorem \ref{T1}(ii) that 
$$ N_{Q_1} = { (c_1-1)+2 \choose 2 } - 3, \quad
N_{Q_2} = { (c_2-1)+2 \choose 2 } - 3.$$
Hence from (\ref{ERL}) 
\begin{eqnarray*}
 N_P &=& N_R - N_{Q_1} - N_{Q_2} + c_1 + c_2 - 8 \\
    &=& [{ h+2 \choose 2 } - 3] \\
  &   &  - [{c_1+1 \choose 2} - 3] - [{c_2+1 \choose 2} - 3] + c_1 + c_2 - 8 \\
    &=& { h+2 \choose 2} - { c_1 \choose 2} - {c_2 \choose 2} - 5, 
\end{eqnarray*}
and the formula (\ref{Nrstdef}) follows. 
This completes the proof of Theorem \ref{T5}. \vskip5mm

\begin{figure}[htb]
\hspace {4cm} 
\xymatrix @M=0ex@R=9ex@C=15ex{
5 \ar@{-}[d]_-{c_1-1}  \ar@{-}[r]^-{s}  & 5 \ar@{-}[d]^-{c_2-1} \\
	 5 \ar@{-}[r]_-{t} & 5}
\caption{Boundary of a patch with four boundary points of degree $5$.}
\label{bd5555}
\end{figure}
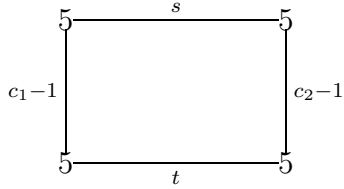

In the case of Type (5), that is,  
$(0,1,4,m)_{b'}$ with $ \beta_5 = 4 $ and $\beta_4 = 0,$ 
let us choose positive integers 
$s,t,c_1,c_2$ such that the four points of degree $5$  
are the end points of four boundary strings with $s,c_1-1,t,c_2-1$ 
edges in this cyclic order respectively (refer to Figure \ref{bd5555}). 

\begin{thm} \label{T6}
Under the assumptions on a patch $P$ of type $(0,1,4,m)_{b'}$
with the lengths of the boundary segments as in Figure \ref{bd5555}
(as stated in the previous paragraph), we have the following: 
\begin{description}
\item{(i)} The integers 
$$ s+t+2(c_1-1)+2(c_2-1), 2s-t +(c_1-1)+(c_2-1),2t-s + (c_1-1)+(c_2-1) $$ 
are  multiples of $3;$ 
\item{(ii)} There exist integers $h,\,k,\,l$ such that 
\begin{eqnarray} \label{hmdef}  
h &=& \frac{1}{3}(s+t+2(c_1-1)+2(c_2-1)), \\
\label{kmdef}
k &=& \frac{1}{3}(2s-t+(c_1-1)+(c_2-1)), \\
\label{lmdef}
l &=& \frac{1}{3}(2t-s+(c_1-1)+(c_2-1)), 
\end{eqnarray}
and such that the patch $P$ has precisely 
\begin{equation} \label{Nmdef}
N = { h+k+2 \choose 2} - { 2k+1 \choose 2} + {k \choose 2} - {c_1 \choose 2}
- { c_2 \choose 2 }  -5 
\end{equation} 
vertex points of degree $6.$  
\end{description}    \end{thm}

Proof: We take the union of the patch $P$ and two 
more patches $Q_1, Q_2$ of types 
$P_{c_1-1}[0,0,0], P_{c_2-1}[0,0,0]$   
respectively in such a way that one of 
the boundary segments between two points of degree 
$4$ of $Q_i$ is identified with the boundary segment of $P$ 
that lies between the two points of degree $5$ which contain  
$c_i-1$ edges ($i=1,2$).

The result of this identification is one patch $R = P \cup Q_1 \cup Q_2 $ which has 
precisely two boundary points of degree $4,$ and it  contains precisely 
one point of degree $4$ in its interior, since by assumption $P \subset R$ 
is of type $(0,1,4,m).$ 
Hence $R$ is a patch of type $(0,3,0)$ with $ \beta_4=2.$ 
We see that the total length of the boundary of $R$ is 
$ s+t+2(c_1-1)+2(c_2-1).$ 
By Theorem \ref{T1}(i) it follows that this boundary length
is a multiple of $3.$ From 
\begin{eqnarray*}
 &&2s -t + (c_1-1)+(c_2-1) \\
= &&3(s+(c_1-1)+(c_2-1)) - (s+t+2(c_1-1)+2(c_2-1)),\\
 &&2t -s + (c_1-1)+(c_2-1)\\
 = &&3(t+(c_1-1)+(c_2-1)) - (s+t+2(c_1-1)+2(c_2-1)),   
\end{eqnarray*}  
it follows that the other integers in (i) are also multiples of $3.$ 

In the identification of $R= P \cup Q_1 \cup Q_2 $ 
we observe that the four boundary points 
of degree $5$ of $P$ and 
two of the three boundary points of degree $4$ of $Q_1, Q_2$ 
each are identified, and become in $R$ boundary points of degree $6.$ By taking proper account of this fact we see that 
for the numbers of points of 
degree $6$ of the relevant patches,   
\begin{equation} \label{ERK}
N_R = N_P + N_{Q_1} +N_{Q_2} - (c_1-2)-(c_2-2)+4. \end{equation}
By Theorem \ref{T1}(iii), we see that the length of the boundary 
of $R$ is given as $ 3h=s+t+2(c_1-1)+2(c_2-1)$   
and that the lengths of the two boundary parts of $R$ between the 
boundary points of degree $4$ are 
$s+c_1+c_2-2=2h-l$ and 
$t+c_1+c_2-2 = h+l$. 
Thus $h$ and $l$ are as  
given in (\ref{hmdef}) and (\ref{lmdef}). 
Moreover, since $k+l=h$ (by Theorem \ref{T1}(iii)),
we have that $k$ is the integer in (\ref{kmdef}). 
For the number of vertex points of degree $6$ of the patch $R$ we see 
from Theorem \ref{T1}(iii) that 
$$ N_R = { h+k+2 \choose 2} - {2k+1 \choose 2} + { k \choose 2} - 3. $$
By Theorem \ref{T1}(ii),  
$$ N_{Q_1} = { (c_1-1)+2 \choose 2 } - 3,\quad
N_{Q_2} = { (c_2-1)+2 \choose 2 } - 3.$$
Hence from (\ref{ERK}) 
\begin{eqnarray*}
 N_P &=& N_R - N_{Q_1} - N_{Q_2} + c_1 + c_2 - 8 \\
    &=& [{ h+k+2 \choose 2 } - {2k+1 \choose 2} + { k \choose 2} - 3] \\
  &   &  - [{c_1+1 \choose 2} - 3] - [{c_2+1 \choose 2} - 3] + c_1 + c_2 - 8 \\
    &=& { h+k+2 \choose 2} - {2k+1 \choose 2} + {k \choose 2}- { c_1 \choose 2} 
     - {c_2 \choose 2} - 5, 
\end{eqnarray*}
and the formula (\ref{Nmdef}) follows. 
This completes the proof of Theorem \ref{T6}. \vskip5mm
  
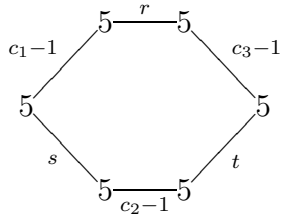
\begin{figure}[htb]
\hspace {5cm} 
\xymatrix @M=0ex{
& 5 \ar@{-}[dl]_-{c_1-1}  \ar@{-}[r]^-{r} & 5 \ar@{-}[dr]^-{c_3-1} & \\
	5 \ar@{-}[dr] _-{s} &  & &5 \ar@{-}[dl]^-{t} \\
& 5 \ar@{-}[r] _-{c_2-1}& 5 & }
\caption{Boundary of a patch with six boundary points of degree $5$.}
\label{bd555555}
\end{figure}

In the case of Type (6), that is, 
$(0,0,6,m)_{b'}$ with $ \beta_5 = 6 $ and $\beta_4 = 0,$ 
let us choose positive integers 
$r,s,t,c_1,c_2,c_3$ such that the six points of degree $5$  
are the end points of six boundary strings with $r,c_1-1,s,c_2-1,t,c_3-1$ 
edges in this cyclic order respectively (refer to Figure \ref{bd555555}). 

\begin{thm} \label{T7}
Under the assumptions on a patch $P$ of type $(0,0,6,m)_{b'}$
with the lengths of the boundary segments as in Figure \ref{bd555555}
(as stated in the previous paragraph), we have the following: 
\begin{description}
\item{(i)} The integer $ r+s+t+2(c_1-1)+2(c_2-1)+2(c_3-1) $ 
is a  multiple of $3;$ 
\item{(ii)} The equality
\begin{equation} \label{rc123def}  
r+ c_1+c_3-2  = s+ c_1 + c_2 -2 = t+c_2+c_3-2  
\end{equation} holds;
\item{(iii)} The patch $P$ has precisely 
\begin{equation} \label{Nnndef}
N = { h+2 \choose 2} - {c_1 \choose 2}- { c_2 \choose 2 } - {c_3 \choose 2}   - 6 
\end{equation} 
vertex points of degree $6$ where $h=r+c_1+c_3-2$.     
\end{description}    \end{thm}

Proof: We take the union of the patch $P$ and three 
more patches $Q_1, Q_2, Q_3$ of types 
$P_{c_1-1}[0,0,0], P_{c_2-1}[0,0,0], 
P_{c_3-1}[0,0,0]$ respectively in such a way that 
one of the boundary segments between two points of degree 
$4$ of $Q_i$ is identified with the boundary segment of $P$ 
that lies between the two points of degree $5$ which contain  
$c_i-1$ edges ($i=1,2,3$).

The result of this identification is one patch 
$R = P \cup Q_1 \cup Q_2 \cup Q_3 $ 
which has precisely three boundary points of degree $4.$ 
Hence $R$ is a patch of type $(0,3,0)$ with $ \beta_4=3.$ 
We see that the total length of the boundary of $R$ is 
$ r+s+t+2(c_1-1)+2(c_2-1)+2(c_3-1).$ 
By Theorem \ref{T1}(i) it follows that this boundary length
is a multiple of $3.$ 

In the identification of $R= P \cup Q_1 \cup Q_2 \cup Q_3 $ 
we observe that the six boundary points 
of degree $5$ of $P$ and two of the three boundary points 
of degree $4$ of $Q_1, Q_2, Q_3$ 
each are identified, and become in $R$ boundary points of degree $6.$ By taking 
proper account of this fact we see that 
for the numbers of points of 
degree $6$ of the relevant patches,    
\begin{equation} \label{ERN}
N_R = N_P + N_{Q_1} +N_{Q_2}+N_{Q_3} - (c_1-2)-(c_2-2)-(c_3-2)+6. 
\qquad\end{equation}
By Theorem \ref{T1}(ii) we see that the length of the boundary 
of $R$ is $ 3h=r+s+t+2(c_1-1)+2(c_2-1)+2(c_3-1)$  and   
that the lengths of the three boundary parts of $R$ between the 
boundary points of degree $4$ are equal, that is, 
$$ r+ c_1+c_3-2  = s+ c_1 + c_2 -2 = t+c_2+c_3-2 = h, $$
which shows part (ii). 
For the number of vertex points of degree $6$ of the patch $R$ we have  
by Theorem \ref{T1}(ii) that 
$$ N_R = { h+2 \choose 2} - 3. $$
We also have by Theorem \ref{T1}(ii) that   
$$ N_{Q_i} = { (c_i-1)+2 \choose 2 } - 3,\qquad i=1,2,3.$$
Hence from (\ref{ERN}) 
\begin{eqnarray*}
 N_P &=& N_R - N_{Q_1} - N_{Q_2} - N_{Q_3} + c_1 + c_2 + c_3 - 12 \\
    &=& [{ h+2 \choose 2 } - 3 ]-[{c_1+1 \choose 2} - 3] - [{c_2+1 \choose 2} - 3] - 
[{c_3+1 \choose 2} - 3]  \\
  &    & + c_1 + c_2 + c_3  - 12 \\
    &=& { h+2 \choose 2} - {c_1 \choose 2} - { c_2 \choose 2} 
     - {c_3 \choose 2} - 6, 
\end{eqnarray*}
and the formula (\ref{Nnndef}) follows. 
This completes the proof of Theorem \ref{T7}. \vskip5mm

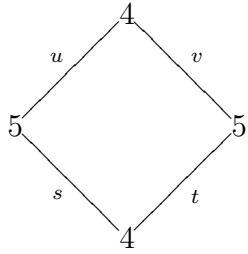
\begin{figure}[htb]
\hspace {5.5cm} 
\xymatrix @M=0ex@R=10ex@C=10ex @ur{
5 \ar@{-}[r]^-{u} & 4 \ar@{-}[d]^-{v}  \\
4 \ar@{-}[u]^-{s} &  5 \ar@{-}[l]^-{t} }
\caption{Boundary of a patch with boundary points $4,5,4,5$.}
\label{bd4545}
\end{figure}

Finally, in the case of Type (7), that is,  
$(0,2,2,m)_{b'}$ with $ \beta_5 = 2 $ and $\beta_4 = 2,$ 
where the boundary points of degree less than $6$ are 
arranged in cyclic order as $4,5,4,5,$ let us choose positive integers 
$u,v,s,t$ such that one of the boundary points of degree $4$  
has edge distances $u$ and $v$ from the two boundary points of degree 
$5$ and such that the second boundary point of degree 
$4$ has edge distances $s$ and $t$ from the two 
boundary points of degree $5$ (refer to Figure \ref{bd4545}). 

\begin{thm} \label{T8}
Under the assumptions on a patch $P$ of type $(0,2,2,m)_{b'}$
with the lengths of the boundary segments as in Figure \ref{bd4545}
(as stated in the previous paragraph), we have the following: 
\begin{description}
\item{(i)} The total length of the boundary 
$ b = u+v+s+t $ is an even number;  
\item{(ii)} The equality
\begin{equation} \label{rdef}  
u = t, v= s  
\end{equation}
holds;
\item{(iii)}
The patch $P$ has precisely  
\begin{equation} \label{Nuvdef}
N = (u+1)(v+1) - 4 
\end{equation} 
vertex points of degree $6.$  
\end{description}   \end{thm} 

\begin{figure}[htb]
\hspace {5.5cm} 
\xymatrix @M=0ex@R=10ex@C=10ex @ur{
5 \ar@{-}[r] \ar@{-}[dr]& 4 \ar@{-}[d]  \\
4 \ar@{-}[u]  &  5 \ar@{-}[l] }
\caption{Patch with $u=1, v=1$.}
\label{sm4545}
\end{figure}
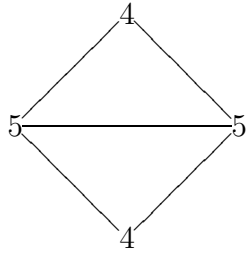

Proof: Clearly, $u,v>0$.    
The proof proceeds by induction on the integer $u+v$.   
If $u+v=2$ (that is,  
$u=v=1$), it is easy to see that there is precisely 
one patch $P$ (refer to Figure \ref{sm4545}) for which we 
clearly have $b=4$ is even and $ s=t=1,$ so that  $ u=t, s=v$ holds. 
Also $N=0=(2)(2)-4$ is immediate in this case. 

Assume now that $u+v>2$, and that the statements 
(i), (ii) and (iii) are already true for 
all patches of side lengths $ u',v'$ 
with $w=u'+v' < u+v  .$ Then for any patch $P$ under consideration with 
side lengths $u$ and $v,$ we have that at least one of $u,v $ is larger 
than one and, without loss of generality, 
we may assume that $u>1$.   

Now we consider the boundary segment of $P$ of length $v,$ 
which lies between a boundary point of degree $4$ 
and a boundary point of degree 
$5.$ There are precisely $2v$ triangles which are incident with this 
given boundary segment. These triangles form a reduced strip (refer to 
Figure \ref{rstrip}). If we remove this reduced strip from $P,$ and then 
take the topological closure, since $u>1$, we obtain a patch $R$ which 
again has two boundary points of degree $4,$ and two boundary points of degree 
$5.$ These points are cyclically arranged on the boundary as $4,5,4,5.$ 
For the new boundary point $q$ of degree $4$ of $R$, we see that the boundary 
segments incident with it have edge lengths $u-1$ and $v.$ 
We also see that the two other boundary segments of $R$ have edge lengths 
$t-1$ and $s.$ 
As the sum of $(u-1)$ and $v$ is $ (u-1)+v < u+v,$ we may apply the induction 
hypothesis to the patch $R.$ This implies that its total boundary  length  
$ b'=(u-1)+v+(t-1)+s $ is an even integer, and hence $b=b'+2$ is also even. 
This proves part (i). 
In a similar way for the patch $R$ 
we obtain $v=s$ and $ u-1 = t-1,$ which implies 
that $u=t,$ thus proving part (ii). 
Finally, the number of vertex points of degree $6$ of $R$ is 
$N' = (u-1+1)(v+1)-4=u(v+1)-4$, 
and the number $N$ of vertex points of degree $6$ 
of the patch $P$ is precisely $(v-1)+2=v+1$ larger than $N'$.   
Hence    
$$ N= N'+ v+1 = u(v+1)-4 + v+1 = (u +1)( v+1)-4. $$  
This proves part (iii).

\section{Construction of Generic Patches}

For this construction we consider the boundary points of degree four 
($\beta_4$) and of 
degree five ($\beta_5$) of a known patch $P.$
We can construct the new patch $P'$ from the known patch 
$P$ by putting a belt around 
the patch $P$ and pushing out the boundary points of degree four 
or degree five or of both degrees out of its original position. 
 
\subsection{Construction of Generic Patches with $\beta_4\ge 0$} 

Let us consider a patch $P$ with $\beta_4 = c,$ where $c>0.$ 
Then we can consider 
$k$ boundary points of degree $4$ ($\beta_4=k$),
$1 \leq k \leq c,$ for this construction, and thus,  
construct $c$ distinct patches $P'$
from $P.$
For the case of $\beta_4=k,$ let us denote the points of degree four on the boundary by 
$a_1, a_2,\dots,a_k.$ First we introduce one new point opposite 
(on the outside) to each boundary edge of the patch $P,$ and connect 
these new points with the vertex points of the corresponding 
edges of $P.$ Then we introduce two new points, say $x_i,y_i$ opposite 
(on the outside)
to the two edges
which are connected with the point $a_i$ for each $i.$ 
Connect the new points $x_i,y_i$
with each other and then connect them with the point $a_i$
for each $i.$ Now introduce $k$ more new points opposite 
(on the outside) to the edges of $x_iy_i,$ 
for each $i.$ Finally, connect all the new points by a belt. 
This gives the patch $P'.$ 
The patch $P'$ is called a generic patch of the patch $P$ with $\beta_4 >0.$
Note that if $\beta_4=0$ then we have a belt around the patch $P.$\\

\begin{lem}  \label{lemgen4}
Let $P$ be a patch of the type $(a_3,a_4,a_5,a_6)_b$ with $\beta_4 \geq 0.$ If we do the generic construction on $P,$ then we obtain a generic patch $P'$ of the type
 $({a_3}',{a_4}',{a_5}',{a_6}')_{b'}$ with $\beta_4 \geq 0.$ Here
\begin{description}
\item (i) ${a_i}' = a_i \mbox{~~~for~~~} i=3,4,5; $ 
\item (ii) ${a_6}' = a_6+b+3\beta_4;$ 
\item (iii) $b'=b+3\beta_4.$
\end{description}
\end{lem}

\begin{figure}[htb]
\hspace {4cm} 
\xymatrix @M=0ex@R=3.2ex@C=3ex{
& & &  {\bullet}\ar@{-}[dr] \ar@{-}[dl]& & &  \\
& & {\bullet}\ar@{-}[dr]\ar@{-}[dl]="a" \ar@{-}[rr] & & {\bullet}\ar@{-}[dr]="c" \ar@{-}[dl] & &  \\
& {\bullet} \ar@{-}[d]\ar@{-}[rr]& & {5}\ar@{-}[drr]\ar@{-}[dll]\ar@{-}[rr] & & {\bullet} \ar@{-}[d]&  \\
& {2} \ar@{-}[rrrr]\ar@{-}[dd]="b" \ar@{-}[drr]\ar@{-}[ddrr]& & & & {3}\ar@{-}[dd]="d"
	\ar@{-}[dll]\ar@{-}[ddll] & \\
& & & {1} \ar@{-}[d]& & & & \\
& {\bullet} \ar@{-}[rr]& & {4} & & {\bullet}\ar@{-}[ll] & &\\
&&&&&&&\ar@/_1.5pc/ @{-}"a";"b" \ar@/^1.5pc/ @{-}"c";"d" \ar@/_1.5pc/ @{-}"b";"d"}
\caption{Generic patch with $\beta_4=1$ of the patch $(1,1,1,2)_4$.}
\label{eggen1}
\end{figure}
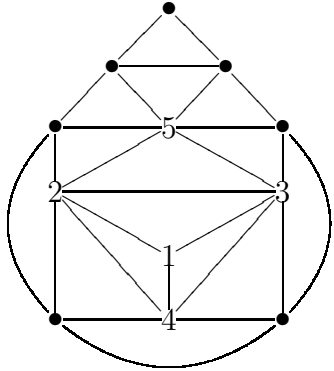

\noindent
{\bf Example:} Consider the patch of type $(1,1,1,2)_4$ given by $ \{123,124,134,235 \}$.
It has a boundary point of degree four $(\beta_4=1)$. Then by 
the generic construction described above   
we can get a patch $(1,1,1,9)_7.$  
This patch also contains a boundary point of degree four 
on its boundary ($ \beta_4=1$) (see Figure \ref{eggen1}). 

\subsection{Construction of Generic Patches with $\beta_5\ge 0$}
Let us consider a patch $P$ with $\beta_5 = c,$ where $c>0.$ 
Then we can consider 
$k$ boundary points of degree $5$ ($ \beta_5=k$),
$1 \leq k \leq c,$ for this construction, and thus, 
construct $c$ distinct patches $P'$
from the patch $P.$
For the case of $\beta_5=k,$ we introduce one new point opposite 
(on the outside) to each boundary edge of the patch $P$ 
and connect these new points with the vertex points of the corresponding 
edges of $P.$ Then introduce $k$ new points opposite 
(on the outside) to the $k$ points of degree
five which are on the boundary of $P$, 
and connect these new points to the corresponding points of degree five. 
Finally, connect all the new points by a belt. 
This gives the patch $P'.$ 
The patch $P'$ is called a generic patch of the patch $P$ with $\beta_5>0.$
Note that if $\beta_5=0$ then we have a belt around the patch $P.$

\begin{lem}  \label{lemgen5}
Let $P$ be a patch of the type $(a_3,a_4,a_5,a_6)_b$ with $\beta_5 \geq 0.$ 
If we do the generic construction on $P,$ then we obtain a generic patch $P'$ of the type
$({a_3}',{a_4}',{a_5}',{a_6}')_{b'}$ with $ \beta_5 \geq 0.$ Here 
\begin{description}
\item (i) ${a_i}' = a_i \mbox{~~~for~~~} i=3,4,5;$ 
\item (ii) ${a_6}' = a_6+b+\beta_5;$ 
\item (iii) $b'=b+\beta_5.$
\end{description}
\end{lem}

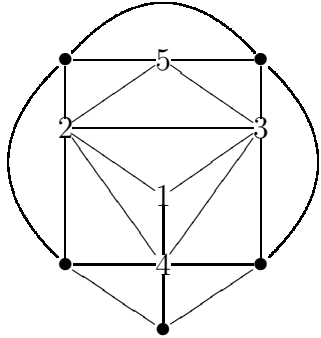
\begin{figure}[htb]
\hspace {4.2cm} 
\xymatrix @M=0ex@R=3.5ex@C=3ex{
& & &&& & \\
& {\bullet}\ar@{-}[d] &&{5} \ar@{-}[drr]\ar@{-}[dll]\ar@{-}[rr]="c"\ar@{-}[ll]="a"&&{\bullet}\ar@{-}[d] & \\
&{2}\ar@{-}[dd]="b"\ar@{-}[rrrr]\ar@{-}[drr]\ar@{-}[ddrr]& &&&{3}\ar@{-}[dd]="d" 
	\ar@{-}[dll] \ar@{-}[ddll]& \\
& &&{1}\ar@{-}[d] & && \\
&{\bullet}\ar@{-}[rr]\ar@{-}[drr]& &{4}\ar@{-}[rr] \ar@{-}[d]&&{\bullet}\ar@{-}[dll] & \\
& &&{\bullet} & &&\ar@/_1.8pc/ @{-}"a";"b" \ar@/^1.8pc/ @{-}"c";"d" \ar@/^1.8pc/ @{-}"a";"c"}
\caption{Generic patch with $\beta_5=1$ of the patch $(1,1,1,2)_4$.}
\label{eggen2}
\end{figure}

\noindent
{\bf Example:} Consider the patch of type $(1,1,1,2)_4$ 
given by $ \{123,124,134,235 \}$.
It has a boundary point of degree five $(\beta_5=1)$. 
Then by the generic construction described for $\beta_5>0$
we can get a patch $(1,1,1,7)_5.$   
This patch also contains a boundary point of degree five on its boundary 
($\beta_5=1$) (see Figure \ref{eggen2}).

\subsection{Construction of Generic Patches with $\beta_4 >0$ and $\beta_5 >0$}

Let us consider a patch $P$ with $\beta_4 >0$ and $\beta_5 >0.$  
Then we can construct a generic patch
$P'$ by considering two boundary points, one with degree four, 
say point $a$ and the other with degree five, say point $b.$  
Consider the case  
$\beta_4=1 $ and $\beta_5=1.$ 
First we introduce one new point opposite 
(on the outside) to each boundary edge of the patch $P,$ 
and connect these new points with the vertex points of the corresponding 
edges of $P.$ Then introduce two new points opposite (on the outside) 
to the two edges
which are connected with the point $a.$ Connect these two new points 
to each other (call this edge $xy$) and to the point $a.$
Then introduce another new point opposite (on the outside)
to the edge $xy$ and connect this new point with the vertex points of the edge $xy.$ Now introduce another new point opposite (on the outside) 
to the point $b$ and connect it with $b.$ 
Finally, connect all the new points by a belt. 
This gives the patch $P'.$ 
The patch $P'$ is called a generic patch of the patch $P$ with $\beta_4=1$ and $\beta_5=1.$
Note that we can generalize this construction for other values of $\beta_4$ and $\beta_5.$

\begin{lem}  \label{lemgen45}
Let $P$ be a patch of the type $(a_3,a_4,a_5,a_6)_b$ with $\beta_4 > 0$ and 
$\beta_5 > 0.$ If we do the generic construction on $P,$  
then we obtain a generic patch $P'$ of the type 
$({a_3}',{a_4}',{a_5}',{a_6}')_{b'}$ with $\beta_4 >0$ and $\beta_5>0.$ Here 

\begin{description}
\item (i) ${a_i}' = a_i \mbox{~~~for~~~} i=3,4,5;$ 
\item (ii) ${a_6}' = a_6+b+3\beta_4 + \beta_5;$
\item (iii) $b' = b+3\beta_4+\beta_5.$
\end{description}
\end{lem}

\begin{figure}[htb]
\hspace {3.9cm} 
\xymatrix @M=0ex@R=3ex@C=3ex{
& & &  {\bullet}\ar@{-}[dr] \ar@{-}[dl]& & &  \\
& & {\bullet}\ar@{-}[dr]\ar@{-}[dl]="a" \ar@{-}[rr] & & {\bullet}\ar@{-}[dr]="c" \ar@{-}[dl] & &  \\
& {\bullet} \ar@{-}[d]\ar@{-}[rr]& & {5}\ar@{-}[drr]\ar@{-}[dll]\ar@{-}[rr] & & {\bullet} \ar@{-}[d]&  \\
& {2} \ar@{-}[rrrr]\ar@{-}[dd]="b" \ar@{-}[drr]\ar@{-}[ddrr]& & & & {3}\ar@{-}[dd]="d"
	\ar@{-}[dll]\ar@{-}[ddll] & \\
& & & {1} \ar@{-}[d]& & & & \\
& {\bullet} \ar@{-}[rr] \ar@{-}[drr]& & {4} \ar@{-}[d]& & {\bullet}\ar@{-}[ll] \ar@{-}[dll]& &\\
& & & {\bullet} & & & &
\ar@/_1.8pc/ @{-}"a";"b" \ar@/^1.8pc/ @{-}"c";"d" }
\caption{Generic patch with $\beta_4=1$ and $\beta_5=1$ of the patch $(1,1,1,2)_4$.}
\label{eggen3}
\end{figure}
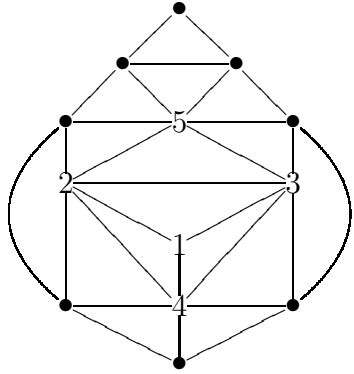

\noindent 
{\bf Example:} Consider the patch $(1,1,1,2)_4$ given by $ \{123,124,134,235 \}$.
It has one point of degree four $(\beta_4=1),$ and one of degree five
$(\beta_5=1)$ on its boundary. Then by the generic construction described
above for $\beta_4=1$ and $\beta_5=1$,    
we can get a patch with parameters   
$(1,1,1,10)_8$. This patch also contains one point of degree four 
and one of degree five on its boundary (see Figure \ref{eggen3}).

\section{Patches of Type $(1,1,1)$ } \label{con111}

In this section we give the formulas for some patches of type $(1,1,1),$
many of which we obtain by using the generic construction method. 

\subsubsection{1.~Type $A$}
$$ \left (1,1,1,2 ~[ ( \frac{k+1}{2})^2-2~ ] +km \right )_{(2k-1)+2m} 
\mbox{for ~ } k \geq 3 ,\mbox{~k~odd~and~~} m\geq 0. $$
This type of patch contains one point of degree $4$ on its boundary   
(that is, $\beta_4=1$) 
and one point of degree $5$ which is almost on the boundary.  

\begin{figure}[htb]
\hspace {4.5cm} 
\xymatrix @M=0ex{
& & x  \ar@{-}[d]   \ar@{-}[rr] &   & {\bullet} \ar@{-}[dd]  \\
{\bullet} \ar@{-}[urr] ="a" \ar@{-}[drr]="b" ^-{B}& & z \ar@{-}[r]  
\ar@{-}[drr]_(0.4){A} \ar@{-}[urr]& 
 {\bullet} \ar@{-}[dr] \ar@{-}[ur]  & \\
& & y  \ar@{-}[u] \ar@{-}[rr] & & {\bullet} \ar@/_1.5pc/ @{-}"a";"b"}
\caption{A patch of the type $(1,1,1,4)_5$.}
\label{Apatch_1}
\end{figure}
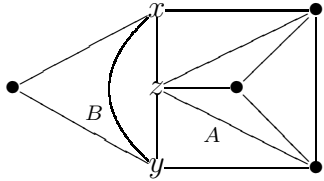

\begin{figure}[htb]
\hspace {4.5cm} 
\xymatrix @M=0ex{
& & {\bullet}\ar@{-}[d] \ar@{-}[r] & x  \ar@{-}[dl] \ar@{-}[d]   \ar@{-}[rr] &   & {\bullet} \ar@{-}[dd]  \\
{\bullet} \ar@{-}[urr] ="a" \ar@{-}[drr]="b" ^-{B}& & {\bullet}\ar@{-}[dr] \ar@{-}[d] \ar@{-}[r] 
& z \ar@{-}[r]  \ar@{-}[drr]_(0.4){A} \ar@{-}[urr]& 
 {\bullet} \ar@{-}[dr] \ar@{-}[ur]  & \\
& & {\bullet} \ar@{-}[r] & y  \ar@{-}[u] \ar@{-}[rr] & & {\bullet} \ar@/_1.5pc/ @{-}"a";"b"}
\caption{Patch of the type $A$ with $k=3$ and $m=1$.}
\label{Apatchb}
\end{figure}
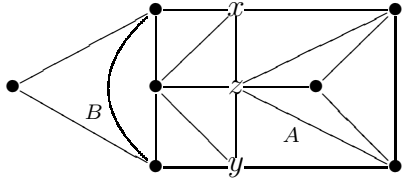

To construct a patch of this type, we first consider the patch $P$
of the type $(1,1,1,4)_5$ (refer to Figure \ref{Apatch_1}).
Let $A$ denote the part of the patch on the right-hand side of $xzy$
and $B$ the part of the patch on the left-hand side of $xzy$. 
Now insert a rectangular
strip of triangles of width $m$ $(m \geq 1)$ which has $3(m+1)$
points and $4m$ triangles in such a way that on side $A$ all the vertex
points of this strip have degree $6$ and on side $B,$ all the points
have the same original degree. 
This construction gives the patch
$(1,1,1,4+3m)_{5+2m}$ for $m\geq1.$
An illustration of the patch 
$(1,1,1,7)_7$ obtained by inserting the rectangular strip with $m=1$
to the patch $(1,1,1,4)_5$ is given in Figure \ref{Apatchb}.
In a similar way we can construct patches of type $A$ for other
odd values of $k\geq3.$ 
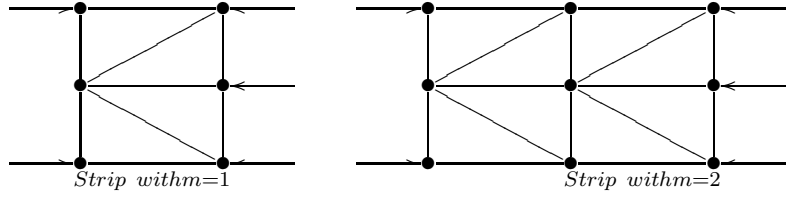
\begin{figure}[htb]
\[ \hspace {1cm} 
\vcenter 
{\xymatrix @M=0ex{
\ar@{-}@{_{>}}[r]  & {\bullet} \ar@{-}[rr]   \ar@{-}[d]& &  {\bullet} \ar@{-}[d] \ar@{-}[dll] & \ar@{-}@{^{>}}[l] \\
& {\bullet} \ar@{-}[d] \ar@{-}[rr] \ar@{-}[drr]& & {\bullet} \ar@{-}[d] & \ar[l] \\
\ar@{-}@{^{>}}[r] & {\bullet} \ar@{-}[rr]_-{Strip ~ with m=1} &  &{\bullet} & \ar@{-}@{_{>}}[l]  }}
\quad \quad
\vcenter
{\xymatrix @M=0ex{
\ar@{-}@{_{>}}[r]  & {\bullet} \ar@{-}[rr]   \ar@{-}[d]& &  {\bullet} \ar@{-}[rr]    \ar@{-}[dll] 
	\ar@{-}[d]& &{\bullet} \ar@{-}[d] \ar@{-}[dll] & \ar@{-}@{^{>}}[l] \\
& {\bullet} \ar@{-}[d] \ar@{-}[rr] \ar@{-}[drr]& & {\bullet} \ar@{-}[rr]   \ar@{-}[drr]  \ar@{-}[d]& 
	& {\bullet} \ar@{-}[d] & \ar[l] \\
\ar@{-}@{^{>}}[r] & {\bullet} \ar@{-}[rr]  &  & {\bullet} \ar@{-}[rr] _-{Strip ~ with m=2}  & 
	&{\bullet} & \ar@{-}@{_{>}}[l]  }}
\]
\caption{Rectangular strips of triangles.}
\label{MuM245}
\end{figure}

Patches of Types $B$ to $G$ as follows can be constructed in a 
similar manner.

\subsubsection{2.~Type $B$}
$$ \left (1,1,1,{k+1 \choose 2} -10+(2k-7)m \right )_{k+2m} 
\mbox{for ~ } k \geq 5 \mbox{,~} m \geq 0. $$
This type of patch contains one point of degree $5$ on its boundary,
that is, $\beta_5=1$ for all cases.

\subsubsection{3.~Type $C$}
$$ \left (1,1,1,{k+1 \choose 2} -20+(2k-11)m \right )_{k+2m} 
\mbox{for ~ } k \geq 7 \mbox{,~} m \geq 0. $$
This type of patch contains one point of degree $5$ on its boundary,
that is, $\beta_5=1$ for all cases.

\subsubsection{4.~Type $D$}
$$ \left (1,1,1,{k+1 \choose 2} -4+(2k-3)m \right )_{k+2m} 
\mbox{for ~ } k \geq 3 \mbox{,~} m \geq 0. $$
This type of patch contains one point of degree $5$ on its boundary,
that is, $\beta_5=1$ for all cases.

\subsubsection{5.~Type $E$}
$$ \left (1,1,1,{k+1 \choose 2} -40+(2k-15)m \right )_{k+2m} 
\mbox{for ~ } k \geq 10 \mbox{,~} m \geq 0. $$
This type of patch contains one point of degree $5$ on its boundary,
that is, $\beta_5=1$ for all cases.

\subsubsection{6.~Type $F$}
$$ \left (1,1,1,{k+1 \choose 2} -124+(2k-27)m \right )_{k+2m} 
\mbox{for ~ } k \geq 18 \mbox{,~} m \geq 0. $$
This type of patch contains one point of degree $5$ on its boundary,
that is, $\beta_5=1$ for all cases.

\subsubsection{7.~Type $G$}
$$ \left (1,1,1,{k+1 \choose 2} -76+(2k-21)m \right )_{k+2m} 
\mbox{for ~ } k \geq 14 \mbox{,~} m \geq 0. $$
This type of patch contains one point of degree $5$ on its boundary,
that is, $\beta_5=1$ for all cases. \\

\subsubsection{8.~Type $H$}
$$ \left (1,1,1,{k+1 \choose 2} -4 \right )_k \mbox{for ~ } k \geq 3. $$
This type of patch can be constructed by using the generic construction 
method with $\beta_5=1$ for $(k-3)$ times from the 
patch $(1,1,1,2)_3$. These patches contain 
one point of degree $5$ on their boundary, that is, $\beta_5=1$ for all cases.

\subsubsection{9.~Type $I$}
$$ \left (1,1,1,{k+1 \choose 2} -8 \right )_k \mbox{for ~ } k \geq 4. $$
This type of patch can be constructed by using the generic construction 
method with $\beta_5=1$ for $(k-4)$ times from the patch $(1,1,1,2)_4$.
These patches contain 
one point of degree $5$ on their boundary, that is, $\beta_5=1$ for all cases. 

\subsubsection{10.~Type $J$}
$$ \left (1,1,1,{k+1 \choose 2} -10 \right )_k \mbox{for ~ } k \geq 5. $$
This type of patch can be constructed by using the generic construction 
method with $\beta_5=1$ for $(k-5)$ times from the patch $(1,1,1,5)_5$.   
These patches contain 
one point of degree $5$ on their boundary, that is, $\beta_5=1$ for all cases. 

\subsubsection{11.~Type $K$}
$$ \left (1,1,1,{k+1 \choose 2} -16 \right )_k \mbox{for ~ } k \geq 6. $$
This type of patch can be constructed by using the generic construction 
method with $\beta_5=1$ for $(k-6)$ times from the patch $(1,1,1,5)_6$.
These patches contain 
one point of degree $5$ on their boundary, that is, $\beta_5=1$ for all cases. 

\subsubsection{12.~Type $L$}
$$ \left (1,1,1,{k+1 \choose 2} -20 \right )_k \mbox{for ~ } k \geq 7. $$
This type of patch can be constructed by using the generic construction 
method with $\beta_5=1$ for $(k-7)$ times from the patch $(1,1,1,8)_7$.
These patches contain 
one point of degree $5$ on their boundary, that is, $\beta_5=1$ for all cases. 

\subsubsection{13.~Type $M$}
$$ \left (1,1,1,{k+1 \choose 2} -26 \right )_k \mbox{for ~ } k \geq 8. $$
This type of patch can be constructed by using the generic construction 
method with $\beta_5=1$ for $(k-8)$ times from the patch $(1,1,1,10)_8$.
These patches contain 
one point of degree $5$ on their boundary, that is, $\beta_5=1$ for all cases. 

\subsubsection{14.~Type $N$}
$$ \left (1,1,1,{k+1 \choose 2} -34 \right )_k \mbox{for ~ } k \geq 9. $$
This type of patch can be constructed by using the generic construction 
method with $\beta_5=1$ for $(k-9)$ times from the patch $(1,1,1,11)_9$.
These patches contain 
one point of degree $5$ on their boundary, that is, $\beta_5=1$ for all cases. 

\subsubsection{15.~Type $O$}
$$ \left (1,1,1,{k+1 \choose 2} -40 \right )_k \mbox{for ~ } k \geq 10. $$
This type of patch can be constructed by using the generic construction 
method with $\beta_5=1$ for $(k-10)$ times from the patch $(1,1,1,15)_{10}$.
These patches contain 
one point of degree $5$ on their boundary, that is, $\beta_5=1$ for all cases. 

\subsubsection{16.~Type $P$}
$$ \left (1,1,1,{k+1 \choose 2} -44 \right )_k \mbox{for ~ } k \geq 10. $$
This type of patch can be constructed by using the generic construction 
method with $\beta_5=1$ for $(k-10)$ times from the patch $(1,1,1,11)_{10}$.
These patches contain 
one point of degree $5$ on their boundary, that is, $\beta_5=1$ for all cases.

\section{Some Known Patches } \label{allpat}

\subsection{Triangular Boundary} \label{tripat}
\begin{table}[h]
\begin{center}
\begin{tabular}{|c|l|l|c|} \hline 
$(a_3,a_4,a_5,a_6)_3$ & patch $P$ & \small{boundary} & \small{remark} \\  \hline \hline 
$(0,3,0,0)_3 $   & $~\{ 123 \} $  & $ 123 $ & $f_3=1, \beta_4=3$ \\ \hline
$(1,0,3,0)_3 $   & $~ \{ 124,134,234 \}$  
  & $123$ & $f_3=3, \beta_5=3 $ \\  \hline
$(1,1,1,2)_3$ & $ ~\{125,145,245,134,234 \} $ & $123$& $f_3=5, \beta_5=1$ \\
\hline
\end{tabular}
\end{center}
\caption{Some patches with boundary length $3$}
\label{tribd}
\end{table}

\newpage

\subsection{Rectangular Boundary} \label{rectpat}

\begin{table}[h]
\begin{center}
\begin{tabular}{|c|l|l|c|} \hline 
$(a_3,a_4,a_5,a_6)_4$ & patch $P$ & \small{boundary} & \small{remark} \\  \hline \hline 
$(0,1,4,0)_4 $   & $~\{ 125,145,235,345 \} $  & $ 1234 $&$f_3=4$ \\ 
		  &		& 		& $\beta_5=4$\\ \hline
$(0,2,2,0)_4 $   & $~ \{ 123,124 \}$  & $1324$&$f_3=2$ \\  
		  &		& 		& $\beta_4=2$\\ 
		  &		& 		& $\beta_5=2$\\ \hline
$(0,2,2,2)_4 $   & $~ \{ 136,145,156,236,245,256 \}$  & $1423$&$f_3=6$ \\ 
		  &		& 		& $\beta_5=2$\\ \hline
$(0,2,2,3)_4 $   & $~ \{ 127,145,157,236,267,345,356,567 \}$  & $1234$&$f_3=8$ \\ 
		  &		& 		& $\beta_5=1$\\ \hline
$(1,0,3,2)_4 $   & $~ \{ 126,135,156,245,256,345 \}$  & $1243$ & $f_3=6$\\  
		  &		& 		& $\beta_5=2$\\ \hline
$(1,0,3,4)_4 $   & $~ \{ 123,124,134,238,247,278 $   & $5678$ & $f_3=10$\\  
		     & $~~~346,356,358,467 \}$& & $\beta_5=1$\\ \hline
$(1,1,1,2)_4$    & $ ~\{123,124,134,235 \} $ & $ 2435 $ & $f_3=4$ \\
 		  &		& 		& $\beta_4=1$\\ 
		  &		& 		& $\beta_5=1$\\ \hline 		  
\end{tabular} 
\end{center}
\caption{Some patches with boundary length $4$}
\label{rectbd}
\end{table}

\newpage

\subsection{Pentagonal Boundary} \label{penpat}
\begin{center}
\begin{tabular}{|c|l|l|c|} \hline 
$(a_3,a_4,a_5,a_6)_5$ & patch $P$ & \small{boundary} & \small{remark} \\  \hline \hline 

$(0,0,6,0)_5 $ & $~ \{ 126,156,236,346,456 \}$  & $12345$ & $f_3=5$  \\ 
		  &		& 		& $\beta_5=5$\\ \hline
$(0,1,4,2)_5 $ & $~ \{ 137,146,167,237,256,267,456 \}$  & $14523$ & $f_3=7$  \\ 
		  &		& 		& $\beta_5=3$\\ \hline
$(0,1,4,3)_5 $ & $~ \{ 128,146,168,257,278,346,357,367,678 \}$& $12534$  & $f_3=9$  \\ 
		  &		& 		& $\beta_5=2$\\ \hline
$(0,1,4,4)_5 $ & $~ \{ 127,156,167,239,279,348, $ &$12345$& $f_3=11$   \\ 
   & $ ~~~389,456,468,678,789 \}$&   & $\beta_5=1$ \\ \hline 
$(0,1,4,5)_5 $ & $~ \{ 127,156,157,238,278,349,389, $ &$12346$& $f_3=13$   \\
  &  $~~~ 456,459,57A,59A,78A,89A \}$  &  & $\beta_5=1$ \\ \hline 
$(0,2,2,1)_5 $ & $ ~\{ 123,125,134 \} $& $15234$  & $f_3=3$  \\ 
		  &		& 		& $\beta_4=2$\\ 
		  &		& 		& $\beta_5=2$\\ \hline
$(0,2,2,2)_5 $ & $~\{ 125,126,136,246,346 \} $ & $ 13425$ & $f_3=5$  \\ 
		  &		& 		& $\beta_4=1$\\ 
		  &		& 		& $\beta_5=2$\\ \hline
$(0,2,2,4)_5 $ & $~\{ 128,146,148,237,278,345,347,456,478 \}$ & $ 12356$ & $f_3=9$  \\ 
		  &		& 		& $\beta_5=2$\\ \hline
$(0,2,2,5)_5 $   & $~ \{ 123,124,136,147,156,157,238, $ &$56897$ & $f_3=11$ \\
           & $~~~ 249,289,368,479  \}$&  & $\beta_5=1$\\ \hline
$(0,2,2,7)_5 $   & $~ \{ 136,138,156,145,149,189,236, $ &$789AB$ & $f_3=15$ \\
           & $~~~ 23B,256,245,24A,2AB,378,37B,49A  \}$&  & $\beta_5=1$\\ \hline
$(1,0,3,3)_5 $ & $~ \{ 124,127,136,137,145,156,237 \}$& $24563$ & $f_3=7$  \\ 
		  &		& 		& $\beta_5=3$\\ \hline
$(1,0,3,5)_5 $ & $~ \{ 128,156,158,247,248,345, $ &$12736$  & $f_3=11$   \\
   & $~~~ 347,356,459,489,589 \}$&  & $\beta_5=2$ \\ \hline 
\end{tabular}  
\end{center}
\begin{flushright}
{\it (continued on next page)}
\end{flushright}

\newpage 
{\it (continued from previous page)}
\begin{table}[h]
\begin{center}
\begin{tabular}{|c|l|l|l|} \hline
$(a_3,a_4,a_5,a_6)_5$ & patch $P$ & \small{boundary} & \small{remark} \\  \hline \hline 

$(1,0,3,6)_5 $ & $~ \{ 125,159,179,238,258,346,368,  $  &$12347$  & $f_3=13$   \\
    & $~~~ 469,479,569,56A,58A,68A  \}$&  & $\beta_5=1$ \\ \hline 
$(1,0,3,7)_5 $ & $~ \{ 12A,158,15A,237,27A,346,367,456,  $  & $12348$ &$f_3=15$   \\
   & $~~~  458,569,59A,67B,69B,79A,79B \}$&  & $\beta_5=1$ \\ \hline 
$(1,1,1,4)_5 $   & $~ \{ 127,145,157,235,257,345,346 \}$& $12364$ & $f_3=7$  \\ 
  	    &		& 		& $\beta_4=1$\\ \hline
$(1,1,1,5)_5 $   & $~ \{ 125,157,167,238,258,345,358,457,467 \}$& $12346$ & $f_3=9$   \\ 
		  &		& 		& $\beta_5=1$\\ \hline
$(1,1,1,7)_5 $   & $~ \{ 126,156,158,239,269,347,379, $ &$12348$ & $f_3=13$ \\
           & $~~~ 457,458,56A,57A,679,67A  \}$&  & $\beta_5=1$\\ \hline
$(1,1,1,11)_5 $   & $~ \{ 125,128,134,137,145,178,245, $ &$ABCDE$ & $f_3=21$ \\
           & $~~~ 234,239,289,376,369,67B,69E,6AB,$ &  & $\beta_5=1$ \\
           & $~~~ 6AE,78C,7BC,8CD,89D,9DE\}$&  & \\ \hline
\end{tabular}  
\end{center}
\caption{Some patches with boundary length $5$}
\label{penbd}
\end{table}

\clearpage
\subsection{Hexagonal Boundary}\label{hexpat}
\begin{center}
\begin{tabular}{|c|l|l|l|} \hline 
$(a_3,a_4,a_5,a_6)_6$ & patch$P$  & \small{boundary} & \small{remark} \\  \hline \hline
$(0,0,6,1)_6 $   & $~ \{ 123,127,134,145,156,167 \}$  & $234567  $ & $f_3=6$ \\  
		  &		& 		& $\beta_5=6$\\ \hline
$(0,0,6,2)_6 $   & $~ \{ 134,135,148,234,236,247,356,478  \}$  & $156278 $ & $f_3=8$\\  
		  &		& 		& $\beta_5=4$\\ \hline
$(0,0,6,3)_6 $   & $~ \{ 146,147,169,245,247,258,356,358,  $ &$172839$&  $f_3=10$ \\
    & $~~~ 369,456  \}$  &  & $\beta_5=3$\\ \hline 
$(0,0,6,4)_6 $   & $~ \{ 129,12A,136,139,16A,245,249,25A, $  &$374586$& $f_3=12$ \\
      &   $~~~ 379,479,58A,68A \}$  &  & $\beta_5=2$\\ \hline 
$(0,0,6,6)_6 $   & $~ \{12A,167,16A,23B,2AB,34C,3BC,458, $ &$123457$&$f_3=16$  \\
    & $~~~ 48C,567,568,689,69A,89C,9AB,9BC  \}$  &  & $\beta_5=1$ \\ \hline
$(0,1,4,2)_6 $ & $~\{ 126,127,136,256,346,456 \}$ & $ 134527$ &$f_3=6$ \\
 		  &		& 		& $\beta_4=1$\\ 
		  &		& 		& $\beta_5=3$\\ \hline
$(0,1,4,3)_6 $ & $~\{ 136,138,178,234,238,278,345,356 \}$ & $172456$ & $f_3=8$ \\ 
		  &		& 		& $\beta_5=4$\\ \hline
$(0,1,4,5)_6 $ & $~\{ 128,157,158,23A,28A,349,39A,456, $   &$123467$& $f_3=12$  \\ 
    & $~~~ 459,567,589,89A  \}$ &  & $\beta_5=2$\\ \hline
$(0,1,4,6)_6 $ & $~\{129,167,169,258,259,34A,358,35A, $ &$128347$&$f_3=14$   \\
    & $~~~ 467,46A,59B,5AB,69B,6AB  \}$ &  &$\beta_5=2$ \\ \hline
$(0,2,2,2)_6 $ & $~\{123,124,235,146  \} $ & $135246$ & $f_3=4$ \\ 
		  &		& 		& $\beta_4=2$\\ 
		  &		& 		& $\beta_5=2$\\ \hline
$(0,2,2,4)_6 $ & $~\{126,146,147,258,268,346,358,368 \}$ & $125347$ & $f_3=8$ \\  
		  &		& 		& $\beta_4=1$\\ 
		  &		& 		& $\beta_5=1$\\ \hline
$(0,2,2,5)_6 $ & $~\{ 125,159,179,258,268,345,358,368, $ &$126347$& $f_3=10$   \\
    & $~~~ 459,479 \}$ &  &  $\beta_5=2$\\ \hline
$(0,3,0,3)_6$    & $ ~\{123,125,134,236 \} $& $152634$& $f_3=4$  \\ 
		  &		& 		& $\beta_4=3$\\ \hline
$(0,3,0,4)_6$    & $ ~\{ 125,127,147,237,346,347 \} $& $152364$& $f_3=6$ \\ 
		  &		& 		& $\beta_4=2$\\ \hline
\end{tabular}
\end{center}
\begin{flushright}
{\it (continued on next page)}
\end{flushright}
\newpage
{\it (continued from previous page)}
\begin{table}[h]
\begin{center}
\begin{tabular}{|c|l|l|l|} \hline
$(a_3,a_4,a_5,a_6)_6$ & patch$P$  & \small{boundary} & \small{remark} \\  \hline \hline

$(0,3,0,6)_6 $ & $~\{ 126,156,157,238,268,349,389,469, $ &$123457$&$f_3=10$     \\
    & $~~~ 456,689 \}$ &  & $\beta_4=1$ \\ \hline
$(1,0,3,5)_6 $   & $~ \{ 129,145,159,235,259,358,368,458, $ &$123674$& $f_3=10$   \\
    & $ ~~~ 478,678 \} $&  & $\beta_5=2$  \\ \hline 
$(1,0,3,6)_6 $   & $~ \{ 145,147,158,256,258,269,346,347, $ &$182937$&$f_3=12$  \\
   & $~~~ 369,45A,46A,56A \}$ & &  $\beta_5=3$  \\ \hline  
$(1,0,3,7)_6 $   & $~ \{ 126,~156,159,26A,28A,347,37A,38A, $ &$128349$& $f_3=14$   \\
    & $~~~ 457,459,56B,57B,67A,67B \}$& & $\beta_5=2$ \\ \hline 
$(1,0,3,8)_6 $   & $~ \{ 126,156,159,237,267,378,38A,78B $ &$123A49$&$f_3=16$  \\
   & $ ~~~ 458,459,48A,56B,58B,67C,6BC,7BC \}$ & & $\beta_5=2$ \\ \hline 
$(1,0,3,9)_6 $   & $~ \{ 126,15B,16B,237,267,34C,37C,49A,  $ &$1234A5$&$f_3=18$   \\ 
     &   $~~~ 49C, 59A,59B,678,68B,7CD,8CD, $  &   &$\beta_5=1$   \\
    & $~~~78D,89B,89C \}$& &   \\ \hline 
$(1,1,1,5)_6 $   & $~ \{ 128,145,158,256,258,345,347,356 \}$& $126374$ & $f_3=8$ \\ 
		  &		& 		& $\beta_4=1$\\ 
		  &		& 		& $\beta_5=1$\\ \hline
$(1,1,1,13)_6 $   & $~ \{ 123,124,134,245,258,28G,23G,346, $ &$7BEDAF$&$f_3=24$  \\
   & $~~~ 369,39G,45C,46C,5BC,57B,578, $  &   & $\beta_5=1$  \\
   & $~~~6CD,69A,6AD,78F,89G,89F,9AF,$ &  &  \\
   & $~~~BCE,CDE \}$ & &    \\ \hline  

\end{tabular}  
\caption{Some patches with boundary length $6$}
\end{center}
\label{Hexabd}
\end{table}

\clearpage
\subsection{Heptagonal Boundary}\label{sevepat}
\begin{table}[h]
\begin{center}
\begin{tabular}{|c|l|l|l|} \hline 
$(a_3,a_4,a_5,a_6)_6$ & patch$P$  & \small{boundary} & \small{remark} \\  \hline \hline 
$(0,0,6,3)_7$   & $~ \{ 128,178,238,349,459,569,$& 1234567 &  $f_3=9$ \\
                & $~~~  679,389,789 \}$  &  & $\beta_5=5$   \\ \hline
$(0,0,6,4)_7$  &  $~\{  128,178,23A,34A,459,569,$& 1234567 & $f_3=11$ \\
               &  $~~~  679,28A,49A,789,89A \}$  &  & $\beta_5=4$    \\ \hline
$(0,0,6,5)_7$  &  $~ \{ 128,178,239,289,34A,39A,45A,$& 1234567 & $f_3=13$ \\
               &  $~~~  56A,67B,6AB,78B,89B,9AB \}$ &  &  $\beta_5=3$  \\ \hline
$(0,0,6,7)_7$  &  $~ \{ 128,178,239,289,34A,39A, $& 1234567 & $f_3= 17$  \\
               &  $~~~  45B,4AB,56B,67C,6BC,78C,  $&  &  $\beta_5=2$  \\ 
               &  $~~~  89D,8CD,9AD,ABD,BCD \} $  &  &   \\  \hline
$(0,0,6,8)_7$  &  $~ \{ 129,178,23A,34B,45B,56C, $& 1234567 & $f_3= 19$  \\
               &  $~~~  67D,189,29A,3AB,5BC,6CD, $&  &  $\beta_5=1$  \\ 
               &  $~~~  78D,89E,8DE,9AE,ABE,     $&    &  \\
               &  $~~~  BCE,CDE \}$          &  &  \\  \hline
$(0,1,4,4)_7$ &  $~ \{ 127,239,348,458,568, 679, $& 1234567 &  $f_3=9$ \\
              &  $~~~  389,689,279 \}$  &  & $\beta_4=1$   \\ 
 		  &		& 		& $\beta_5=2$\\ \hline
$(0,1,4,7)_7$ &  $~ \{  129,178,23A,34C,45C,56C, $& 1234567 & $f_3=15$ \\
              &  $~~~   67B,29A,3AC,6BC,78B,89B  $&  & $\beta_5=2$   \\
              &  $~~~   9AB,189,ABC \}$ &  &    \\ \hline
$(0,2,2,5)_7$ &  $~ \{ 128,239,345,569,679, $& 1234567 &  $f_3=9$ \\
              &  $~~~  178,289,789,359 \}$  &  & $\beta_4=1$ \\
		  &					  &  & $\beta_5=2$ \\ \hline
$(1,1,1,7)_7 $   & $~ \{ 129,145,159,235,259,358,368,458, $ &$632147A$&$f_3=11$  \\
  		 & $~~~ 478,678,67A \}$ & &  $\beta_4=1$  \\ \hline  
$(1,1,1,8)_7 $   & $~ \{ 123,124,134,245,346,457,467,579, $ &$236AB95$&$f_3=13$  \\
   		 & $~~~ 67A,789,78A,89B,8AB \}$ & &  $\beta_5=1$   \\ \hline

\end{tabular}  
\caption{Some patches with boundary length $7$}
\end{center}
\label{sevebd}
\end{table}

\clearpage 
\subsection{Some Other Used Patches}\label{otherpat}
\begin{table}[h]
\begin{tabular}{|c|l|l|l|} \hline 
$(a_3,a_4,a_5,a_6)_6$ & patch$P$  & \small{boundary} & \small{remark} \\  
\hline \hline
$(0,0,6,5)_8 $   & $ \{ 123,12A,134,14B,1AB,235,367, $ &$258964BA$& $f_3=12$ \\
   & $~~ 357,346,578,679,789 \}$ & &  $\beta_5=4$  \\ \hline
$(0,0,6,12)_8 $   & $ \{ 123,125,134,146,15I,16I,259,46A, $ &$234AGHF9$  &$f_3=26$  \\
   & $~~ 59C,57C,57I,68D,68I,6AD,78B, $  &   & $\beta_5=2$  \\
   & $~~78I,7BC,8BD,9CF,ADG,BCE,$ &  &  \\
   & $~~BDE,CEF,DEG,EFH,EGH \}$ & &    \\ \hline

$(1,0,3,15)_8 $   & $ \{ 126,12E,156,159,19D,1DE,26A,$ &$CDEFGHIJ$&$f_3=28$  \\
   & $~~ 2A8,28F,2EF,347,37A,3A8,38H,$ & &$\beta_5=2$ \\
   & $~~ 3HI,34I,457,459,49J,4JI,56B,$ & & \\
   & $~~ 57B,67A,67B,8GH,8GF,9CD,$ & & \\   
   & $~~ 9CJ\}$ & &    \\ \hline  

$(1,1,1,10)_8 $   & $ \{ 123,124,134,23A,25A,257,247,$ &$5BDC6897$& $f_3=16$ \\
   & $~~ 36A,368,348,479,489,5AB,6AC,$ & & $\beta_4=1$\\
   & $~~ ABC,BCD \}$ & &  $\beta_5=1$\\ \hline  
$(1,1,1,11)_{9} $ & $ \{ 123,124,236,246,356,467,569,589, $ &$1358CEA74$ &$f_3=17$ \\
	& $~~ 679,79A,89B,9AB,8BC,ABE,$ & &$\beta_5=1$  \\
	& $~~ BCD,BDE,CDE \}$ &  &  \\ \hline
$(1,1,1,11)_{10} $ & $ \{ 123,124,136,234,345,358,368,$ &$1457DECB96$& $f_3=16$\\
	& $~~ 578,689,78A,7AD,89A,9AB,$ & &$\beta_4=1$ \\
	& $~~ ABC,ACD,CDE \}$ &  &$\beta_5=1$  \\ \hline
$(1,1,1,15)_{10} $ & $ \{ 123,129,19B,234,24A,28A,289, $ &$1367FGHIEB$& $f_3=24$\\
	& $~~ 346,456,45D,4AD,567,57G,$ & &$\beta_4=1$ \\
	& $~~ 5GH,5DH,7FG,89A,9AC,9BC,$ & &$\beta_5=1$ \\
	& $~~ ACD,BCE,CDI,CEI,DHI \}$ &  & \\ \hline
$(1,1,1,14)_{12} $ & $ \{ 123,124,136,234,345,358,368,587,$ &$1457BGHFED96$& $f_3=20$\\
	& $~~ 689,78A,7AB,89A,9AD,ADC,$ & & $\beta_4=1$ \\
	& $~~ ABC,CDE,CEF,CFG,BCG,FGH\}$ &  & $\beta_5=1$ \\ \hline 
\end{tabular}  
\caption{Some other known used patches}
\label{otherbd}
\end{table}

\newpage

\chapter{Elliptic Triangulations of Spheres}\label{ETS}

Let $T$ be a triangulation of a sphere.   
For a point  $x \in T$, 
the number of triangles $ \sigma_1,\sigma_2,\dots,\sigma_d$ 
which contain $x$ is denoted by  $ d=d(x)$ and is called 
the {\it degree} of $x \in T.$ We can easily see that $ d \geq 3.$ 
For a given triangulation $T$ of a sphere,    
let us denote the face numbers of $T$ by
$ f_i=f_i(T)$ ($i=1,2,3$). That is, $f_1$ is the number of points,
$f_2$ the number of edges and $f_3$ the number of triangles in $T.$ 
Then we have the following relation between the face numbers:
\begin{equation} \label{EQ0}
2f_2(T) ~=~3f_3(T). \end{equation}
The integer  $ \chi(T) = f_1(T) - f_2(T) + f_3(T) $ is called   
{\it Euler's characteristic}.   
Now for $d\geq 3$, define the nonnegative integer $ \alpha_d $ as the 
number of points in $T$ which have degree equal to $d,$ that is,  
\begin{equation} \label{EQ1} 
\alpha_d(T) ~= ~ \alpha_d~=~ | \{ x \in T ~:~ d(x) = d~~\}  | ~.  
\end{equation} 
The tuple of numbers  
$$ (\alpha_3, \alpha_4, \alpha_5,\dots, \alpha_m ) $$ 
is called the {\it parameters} 
or the {\it system of parameters} of the triangulation.   
From the definition of $\alpha_d$ we obtain the equations 
\begin{eqnarray} \label{EE1}
\alpha_3 + ~\alpha_4 + ~~\alpha_5 + ~~\alpha_6 + ~~\alpha_7 + ~\alpha_8 +  ... + ~\alpha_m &=& f_1(T) ;\\ \label{E2}
3\alpha_3 + 4\alpha_4 + 5\alpha_5 + 6 \alpha_6 + 7\alpha_7+ 8\alpha_8 + ...+m \alpha_m &=& 2f_2(T); \\ \label{E3}
3\alpha_3 + 4\alpha_4 + 5\alpha_5 + 6 \alpha_6 + 7\alpha_7+ 8\alpha_8 + ...+m \alpha_m &=& 3f_3(T). \end{eqnarray}
By taking the combination $6(\ref{EE1})-3(\ref{E2}) + 2(\ref{E3}) $ 
we easily get 
\begin{equation} \label{E4} 
3\alpha_3 + 2\alpha_4 + \alpha_5  -\alpha_7 -2\alpha_8 - \dots -
(m-6)\alpha_m = 6\chi(T). \end{equation}

In the case of spherical triangulations   
we get $\chi(T)=2$ and equation  (\ref{E4})
becomes
\begin{equation} \label{Ee5} 
3\alpha_3 + 2\alpha_4 + \alpha_5 -\alpha_7 -2\alpha_8 - \dots -
(m-6)\alpha_m = 12. \end{equation}

A triangulation $T$ is called {\it elliptic} if it does not contain any point with degree
greater than $6$, that is, $d(x) \leq 6$ for every $x \in T.$
If $T$ is an elliptic triangulation, equation (\ref{Ee5}) becomes
\begin{equation} \label{Ee6} 
3\alpha_3 + 2\alpha_4 + \alpha_5 = 12.
\end{equation}
There are exactly $19$ solutions 
$(\alpha_3, \alpha_4, \alpha_5)$ for (\ref{Ee6}).
These are 
\begin{eqnarray*}
&&(0,0,12);(0,1,10);(0,2,8);(0,3,6);(0,4,4);(0,5,2);(0,6,0); \\
&&(1,0,9);(1,1,7);(1,2,5);(1,3,3);(1,4,1);(2,0,6);\\ 
&&(2,1,4);(2,2,2);(2,3,0);(3,0,3);(3,1,1);(4,0,0).
\end{eqnarray*} 
For any system of nonnegative integers  $(a_3,a_4,a_5)$ 
which satisfies the equation $ 3a_3 + 2a_4 + a_5 = 12,$
it is known by an old result of Eberhard (1910) ~\cite{Eber} 
that there exist nonnegative integer values  $N=a_6 \geq 0$ and 
a  triangulation $T$ with the property:
$$ (\alpha_3(T),\alpha_4(T), \alpha_5(T), \alpha_6(T)) = (a_3,a_4,a_5,a_6).$$ 
Let S be the following set of integral 4-tuples: 
\begin{eqnarray*} \label{E1}
S & = &  \{ (3,1,1,2m)| m \geq 0 \} \cup \{ (4,0,0,2m+1) | m \geq 0 \}\\ 
  &   &   \cup ~ \{ (0,0,12,1),(0,1,10,0),(0,1,10,1),(0,4,4,1),(0,5,2,1),\\  
  &   &   ~~(0,6,0,1),(1,0,9,0),(1,0,9,1), (1,0,9,2),(1,0,9,4), (1,1,7,0),\\
  &   &   ~~(1,1,7,1),(1,2,5,0),(1,4,1,0),(1,4,1,1),(2,0,6,1),(2,1,4,0), \\
  &   &   ~~(2,3,0,1),(2,3,0,3),(2,3,0,7),(2,3,0,15),(2,3,0,31),(3,0,3,0),\\
  &   &   ~~(3,0,3,2),(3,0,3,4),(3,0,3,12),(3,1,1,1),\\
  &   &   ~~(3,1,1,17),(4,0,0,2). ~ \} 
\end{eqnarray*} 
We shall show that if 
$ (a_3,a_4,a_5,a_6)$ is not contained in the set $S,$ then there 
exists an (elliptic) triangulation $T$ of the sphere with the 
system of parameters   
$$(\alpha_3(T),\alpha_4(T), \alpha_5(T), \alpha_6(T)) =(a_3,a_4,a_5,a_6).$$ 

We remark that the nonexistence of triangulations $(3,1,1,2m)$ and 
$(4,0,0$, $2m+1)$ for $ m \geq 0$ have been shown by Gr\"unbaum and Motzkin
in  \cite{GM} (see also \cite{GR}).

We shall use the following constructive methods to find the elliptic spherical 
triangulations:
\begin{enumerate}
\item Mutant, productive and self-reproductive configurations. 
\item Fullering constructions.
\item Glueing of Patches.
\end{enumerate}

\section{Mutant Configurations} \label{MU}
A configuration is called { \it mutant } 
if adding (deleting) one or more points to (of) it changes the 
parameters of its triangulation in a certain way. We shall
describe two types of mutant configurations where one point is added. 

\subsection{ Type  $M_1:$}  
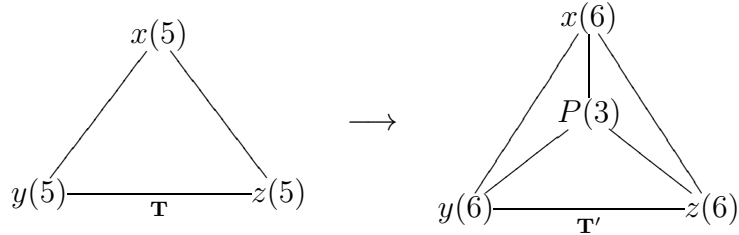
\begin{figure}[htb]
\[ \hspace {1cm} 
\vcenter 
{\xymatrix @M=0ex{
& {x(5)}  \ar@{-}[ddl]  \ar@{-}[ddr] & \\
& & \\
{y(5)} \ar@{-}[rr] _-{\bf{T}}& & {z(5)}}}
\quad \longrightarrow \quad
\vcenter
{\xymatrix @M=0ex{
& {x(6)} \ar@{-}[d] \ar@{-}[ddl]  \ar@{-}[ddr] & \\
&{P(3)}\ar@{-}[dl]  \ar@{-}[dr]   & \\
{y(6)} \ar@{-}[rr] _-{\bf{T'}}& & {z(6)}}}
\]
\caption{Mutant configuration of type $M_1$.}
\label{MuM2}
\end{figure}

Let $X = \{ x,y,z \}$ with 
a triangle $ xyz$ and assume that $d(x)= d(y)= d(z)=5 .$ Then the insertion of a
new point $P$ and the replacement of $xyz$ by 
$ pxy,pxz,pyz$ leads to a new triangulation $T'$ with 
$ d'(x)=d'(y)=d'(z)=6, d'(P) = 3$ (refer to Figure \ref{MuM2}). Hence 
$$ (\alpha_3(T'),\alpha_4(T'), \alpha_5(T'), \alpha_6(T')) = ( \alpha_3(T)+1, \alpha_4(T), \alpha_5(T)-3 ,\alpha_6(T)+3 ).$$
We say that $X$ is a mutant configuration of type $M_1.$

\subsection{ Type  $M_2:$}  
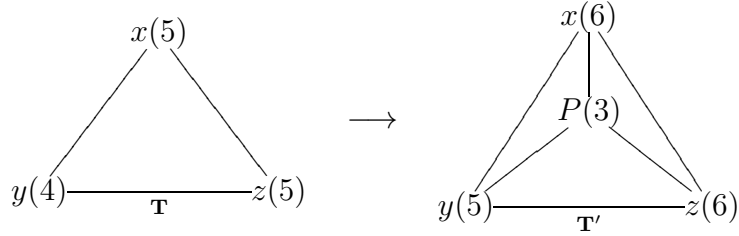
\begin{figure}[htb]
\[ \hspace {1cm}
\vcenter 
{\xymatrix @M=0ex{
& {x(5)}  \ar@{-}[ddl]  \ar@{-}[ddr] & \\
& & \\
{y(4)} \ar@{-}[rr] _-{\bf {T}}& & {z(5)}}}
\quad \longrightarrow \quad
\vcenter
{\xymatrix @M=0ex{
& {x(6)} \ar@{-}[d] \ar@{-}[ddl]  \ar@{-}[ddr] & \\
&{P(3)}\ar@{-}[dl]  \ar@{-}[dr]   & \\
{y(5)} \ar@{-}[rr] _-{\bf{T'}}& & {z(6)}}}
\]
\caption{Mutant configuration of type $M_2$.}
\label{MuM3}
\end{figure}

Let $X = \{ x,y,z \}$ with 
a triangle $ xyz$ and assume that $d(x)= d(z)= 5 ,d(y)=4 .$ Then the insertion of a
new point $P$ and the replacement of $xyz$ by 
$ pxy,pxz,pyz$ leads to a new triangulation $T'$ with 
$ d'(x)=d'(z)=6, d'(P) = 3,d'(y)=5$ (refer to Figure \ref{MuM3}). Hence 
$$ (\alpha_3(T'),\alpha_4(T'), \alpha_5(T'), \alpha_6(T')) =
 ( \alpha_3(T)+1, \alpha_4(T)-1, \alpha_5(T)-1 ,\alpha_6(T)+2 ). $$
We say that $X$ is a mutant configuration of type $M_2.$

\section{Productive and Self-Reproductive Configurations} 
Let $T$ be any triangulation with point set $V.$ 
A subset of points $X \subset V$  is called a {\it productive configuration} 
if it satisfies the following property: \\
{\it At least one point (say $z$) 
can be added to $X$ and certain changes involving triangles of 
$X \cup \{ z \}$ can be made to the triangulation $T$ such that the 
number of points of degree six is increased whereas
the number of points of degrees 3,4,5 remain unchanged.}\\

A subset of points $X \subset V$  is called a  
{ \it self-reproductive configuration} if it satisfies the following 
property: \\ 
{\it At least one point (say $z$)
can be added to $X$ and certain 
changes involving triangles of $X \cup \{ z \}$ 
can be made to the triangulation $T$ such that 
\begin{description}
\item (i) the number of points of degree six is increased 
whereas the number of points of degrees less than 
six remain unchanged; 
\item (ii) the resulting new triangulation 
contains the original configuration as a subconfiguration.
\end{description}
} 

\subsection{Type  $P_1:$}  
\begin{figure}[htb]
\[ \hspace {2cm} 
\vcenter 
{ \xymatrix @M=0ex @ur{
{x(6)} \ar@{-}[r] \ar@{-}[dr]& {p(4)} \ar@{-}[d]  \\
{q(5)} \ar@{-}[u] &  {y(6)} \ar@{-}[l] ^-{\bf{T}}}}
\qquad \longrightarrow \qquad
\vcenter 
{\xymatrix @M=0ex @ur{
{x(6)} \ar@{-}[r] \ar@{-}[dr]& {p(5)} \ar@{-}[d]  \\
{q(6)} \ar@{-}[u] &  {y(6)} \ar@{-}[l] ^-{\bf{T'}}}}
\]
\caption{Mutant configuration of type $P_1$.}
\label{Pp1}
\end{figure}
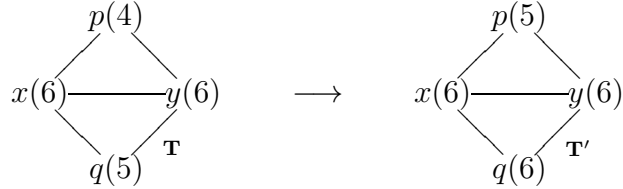

Let $X = \{ p,x,q,y \}$ with 
two triangles $ pxy,qxy$ and assume that $d(p)=4,$  $d(q)=5,
d(x)=d(y) =6.$ Then the insertion of a new point 
$z$ and the replacement of $pxy,qxy$ by 
$ pxz,pyz,qxz,qyz$ leads to a new triangulation $T'$ with 
$ d'(x)=d'(y)=d'(q)=6,\,d'(p)=5,\,d'(z) = 4$ 
(refer to Figure \ref{Pp1}). Hence 
$$ (\alpha_3(T'),\alpha_4(T'), \alpha_5(T'), \alpha_6(T')) 
= ( \alpha_3(T), \alpha_4(T), \alpha_5(T) ,\alpha_6(T)+1 ).$$
We say that $X$ is a productive configuration of type $P_1.$

\subsection{Type $P_2:$}  
\begin{figure}[htb]
\[ \hspace {0.5cm} 
\vcenter
{\xymatrix @M=0ex{
& v \ar@{-}[dl] \ar@{-}[dr]  \ar@{-} \ar@{-}[rr] &  & w \ar@{-}[dr] \ar@{-}[dl]& \\
u \ar@{-}[dr] \ar@{-}[rr] &  &P \ar@{-}[rr] \ar@{-}[dl] \ar@{-}[dr]& & x \ar@{-}[dl] \\
& z \ar@{-}[rr] _-{\bf{T}}& & y &  }}
\quad \longrightarrow \quad
\vcenter
{\xymatrix @M=0ex{
& v \ar@{-}[dl] ="b"\ar@{-}[dr]  \ar@{-} \ar@{-}[rr]="a" &  & w \ar@{-}[dr] \ar@{-}[dl]& \\
u \ar@{-}[dr] \ar@{-}[rr] &  &P \ar@{-}[rr] \ar@{-}[dl] \ar@{-}[dr]& & x \ar@{-}[dl] \\
& z \ar@{-}[rr] _-{\bf{T'}}& & y & \ar@{-}"a";"b"^(.55){q} }}
\]
\caption{Productive configuration of type $P_2$.}
\label{Pt1}
\end{figure}
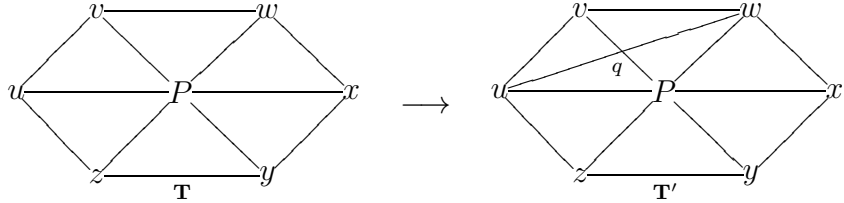

Let $X = \{u,v,w,p,x,y,z \}$ with 
six triangles $ upv,vpw,wpx,xpy,ypz,upz$ and assume that 
$d(p)=d(v)=6, d(u)=d(y)=5$ and 
$d(w)=d(x) =d(z)=4$. Then the insertion of a new point    
$q$ and the replacement of $upv,vpw$ by 
$ upq,uqv,pqw,qvw$ leads to a new triangulation $T'$ with 
$ d'(u)=6, d'(w)=5,d'(q)=4$ (refer to Figure \ref{Pt1}). Hence 
$$ (\alpha_3(T'),\alpha_4(T'), \alpha_5(T'), \alpha_6(T')) 
= ( \alpha_3(T), \alpha_4(T), \alpha_5(T) ,\alpha_6(T)+1 ).$$
We say that $X$ is a productive configuration of type $P_2$.   

Clearly, every self-reproductive configuration must contain 
at least one productive configuration. 
We give a list of self-reproductive configurations with points of degree six or less each.  These are of fundamental importance in the  constructions which follow.

\subsection{Type $A:$} 
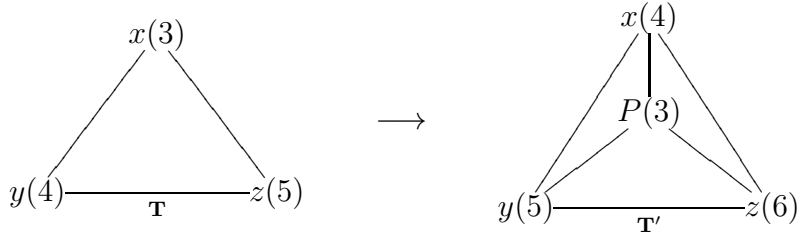
\begin{figure}[htb]
\[ \hspace {1cm}
\vcenter 
{\xymatrix @M=0ex{
& {x(3)}  \ar@{-}[ddl]  \ar@{-}[ddr] & \\
& & \\
{y(4)} \ar@{-}[rr] _-{\bf {T}}& & {z(5)}}}
\qquad \longrightarrow \qquad
\vcenter
{\xymatrix @M=0ex{
& {x(4)} \ar@{-}[d] \ar@{-}[ddl]  \ar@{-}[ddr] & \\
&{P(3)}\ar@{-}[dl]  \ar@{-}[dr]   & \\
{y(5)} \ar@{-}[rr] _-{\bf{T'}}& & {z(6)}}}
\]
\caption{Self-reproductive configuration of type $A$.}
\label{TyA}
\end{figure}

Let $xyz$ be any triple of points in the triangulation $T$ 
and assume that $X= \{ x,y,z \}$ forms a  triangle in $T$ with the
degrees given as $ d(x)=3,d(y) = 4,$ $ d(z) = 5.$  
If we delete the triangle $xyz$ from $T,$ and add a point $p$ and three new triangles $pxy,pxz,pyz $ then we get a new triangulation $T'$ with the degrees 
$ d'(p)=3, d'(x) = 4, d'(y) = 5, d'(z)= 6. $ The degrees 
of all the other points of $T$ remain the same. Hence $T'= T \cup \{p\} $ 
has the parameters  
$$ (\alpha_3(T'),\alpha_4(T'), \alpha_5(T'), \alpha_6(T')) = (\alpha_3(T),\alpha_4(T),\alpha_5(T),\alpha_6(T)+1).$$ 
We thus see that $X$ is a productive configuration  
as well as self-reproductive configuration, since in $T'$
the configuration $pxy$ is isomorphic to $X$ 
(refer to Figure \ref{TyA}). Such a configuration $X$ will be
called a self-reproductive configuration of type $A.$\\

\subsection{Type $B_1$:} 
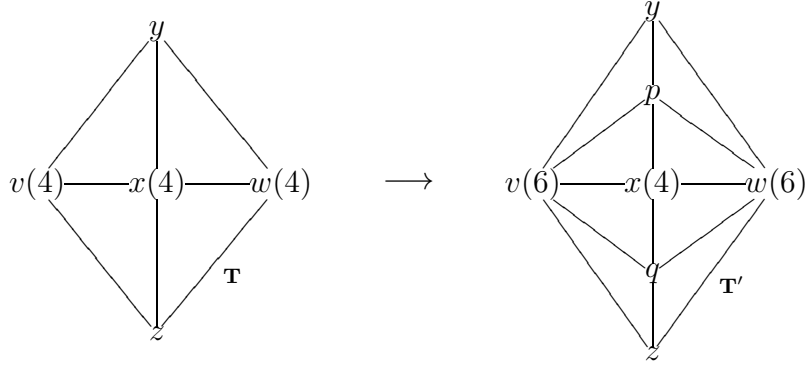
\begin{figure}[htb]
\[ \hspace {1cm} 
\vcenter 
{\xymatrix @M=0ex{
&y \ar@{-}[dd] \ar@{-}[ddl]  \ar@{-}[ddr] & \\
&  & \\
{v(4)} \ar@{-}[r] \ar@{-}[ddr]  & {x(4)} \ar@{-}[dd]  \ar@{-}[r]   
& {w(4)} \ar@{-}[ddl]^-{\bf{T}} \\
&  & \\
& z  & }}
\qquad \longrightarrow \qquad
\vcenter 
{\xymatrix @M=0ex{
&y \ar@{-}[d] \ar@{-}[ddl]  \ar@{-}[ddr] & \\
&p \ar@{-}[d] \ar@{-}[dl]  \ar@{-}[dr] & \\
{v(6)} \ar@{-}[r] \ar@{-}[dr] \ar@{-}[ddr]  & {x(4)} \ar@{-}[d]  \ar@{-}[r]   
& {w(6)} \ar@{-}[dl] \ar@{-}[ddl]^-{\bf{T'}} \\
&q \ar@{-}[d]  & \\
& z  & }}
\]
\caption{Self-reproductive configuration of type $B_1$.}
\label{TyB1}
\end{figure}

Let $vwxyz$ be any five points of the triangulation 
$T$ and assume that $X= \{ v,w,x,y,z \}$ contains exactly four 
triangles $ vxy, vxz,wxy,wxz $ of $T$ such that the degrees 
of the points are given as $d(v)=d(w)=d(x)=4$, $d(y)=d(z)=6$. If we  
delete these given four triangles from $T$ and add  
two new points $p,q$ and eight new triangles $vpx,vpy,vqx,vqz,wpx,wpy,
wqx,wqz$ then we get a new triangulation $T'$ with the degrees  
$ d'(p)=d'(q)=4, d'(v) =d'(w)=6$. The degrees of the 
other points of $T$ remain the same. Hence the triangulation  
$T'= T \cup \{ p,q \} $ will have the parameters 
$$ (\alpha_3(T'),\alpha_4(T'), \alpha_5(T'), \alpha_6(T')) = (\alpha_3(T),\alpha_4(T),\alpha_5(T),\alpha_6(T)+2).$$ 
Thus $X$ is a  productive as well as self-reproductive 
configuration, since in $T'$ the configuration $pqvwx$ is 
isomorphic to $X$ (refer to Figure \ref{TyB1}). Such a configuration $X$ will be called a self-reproductive configuration of type $B_1.$\\

\subsection{Type $B_2$:} 
\begin{figure}[htb]
\[ \hspace {1cm} 
\vcenter 
{\xymatrix @M=0ex{
&{x(5)}  \ar@{-}[ddl]  \ar@{-}[ddr] & \\
&  & \\
{v(4)} \ar@{-}[rr] \ar@{-}[ddr]  & & {w(4)} \ar@{-}[ddl]^-{\bf{T}} \\
&  & \\
& {y(5)}  & }}
\qquad \longrightarrow \qquad
\vcenter 
{\xymatrix @M=0ex{
&{x(6)} \ar@{-}[d] \ar@{-}[ddl]  \ar@{-}[ddr] & \\
&q \ar@{-}[d] \ar@{-}[dl]  \ar@{-}[dr] & \\
{v(5)} \ar@{-}[r]  \ar@{-}[ddr]  & {p} \ar@{-}[dd]  \ar@{-}[r]   
& {w(5)}  \ar@{-}[ddl]^-{\bf{T'}} \\
&  & \\
& {y(6)}  & }}
\]
\caption{Self-reproductive configuration of type $B_2$.}
\label{TyB2}
\end{figure}
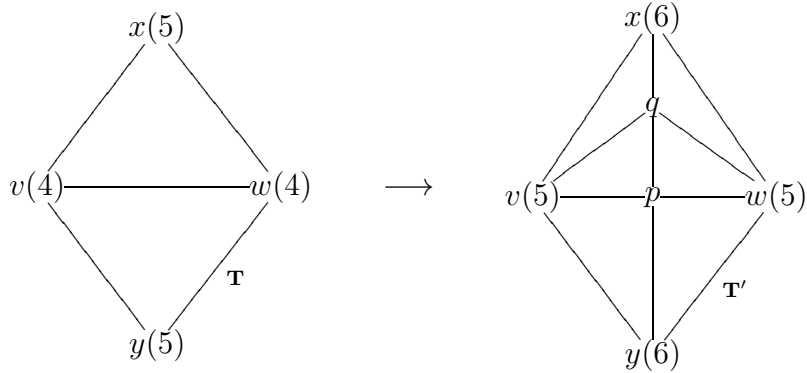

Let $vwxy$ be any four points of the triangulation 
$T$ and assume that $X= \{ v,w,x,y \}$ contains exactly two 
triangles $ vwx, vwy $ of $T$ such that the degrees 
of the points are given as $d(v)=d(w)=4, d(x) =d(y)=5.$ If we  
delete these two triangles from $T$ and add  
two new points $p,q$ and six new triangles $vqx,vpq,vpy,wqx,wpq,wpy $ 
then we get a new triangulation $T'$ with the degrees  
$ d'(p)=d'(q)=4, d'(v) =d'(w)=5, d'(x)=d'(y)=6.$ 
Hence the triangulation 
$T'= T \cup \{ p,q \} $ will have the parameters 
$$ (\alpha_3(T'),\alpha_4(T'), \alpha_5(T'), \alpha_6(T')) = (\alpha_3(T),\alpha_4(T),\alpha_5(T),\alpha_6(T)+2).$$ 
Thus $X$ is a  productive as well as self-reproductive 
configuration, since in $T'$ the configuration $pwqv$ is isomorphic to $X$ 
(refer to Figure \ref{TyB2}). Such a configuration $X$ will be
called a self-reproductive configuration of type $B_2.$\\

\subsection{Type $C:$} Let $vwxyz$ be any five points of the triangulation  
$T$ and assume that $X= \{ v,w,x,y,z \}$ contains three triangles  
$ vwx, wxy,xyz $ of $T$ such that the degrees of the points are given as 
$d(w)=4$ and $ d(v)=d(x)=d(y)=d(z)=5.$
\begin{figure}[htb]
\[ \hspace {0.2cm} 
\vcenter 
{\xymatrix @M=0ex@R=3ex@C=5ex{
&{x(5)}  \ar@{-}[ddl]  \ar@{-}[r] \ar@{-}[ddr] \ar@{-}[dddd]&  {z(5)} \ar@{-}[dd]\\
&   &\\
{v(5)}  \ar@{-}[ddr]  &   & {y(5)} \ar@{-}[ddl]^-{\bf{T}} &\\
&  & \\
& {w(4)}  & }}
\longrightarrow \quad
\vcenter 
{\xymatrix @M=0ex@R=3ex@C=2.5ex{
& &{x(5)}  \ar@{-}[ddll]  \ar@{-}[rr] \ar@{-}[dr] \ar@{-}[dddd]& & {z(6)} \ar@{-}[dl]\ar@{-}[dd]\\
&  & &  p \ar@{-}[dr] \ar@{-}[dddl]&\\
{v(5)}  \ar@{-}[ddrr] &  &  & & {y(5)} \ar@{-}[ddll]^-{\bf{T'}} \\
&   &&  & \\
&  &{w(5)}  &  &}}
\]
\caption{Self-reproductive configuration of type $C$.}
\label{striph2}
\end{figure}
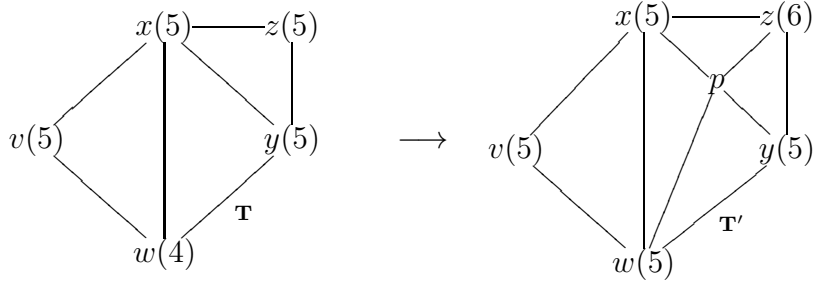

If we delete the two triangles $wxy,xyz $ from $T,$ and  
add a point $p$ with four new triangles $wpx,wpy,pxz,pyz$ 
then we get a new triangulation  $T'$ with the degrees of points as 
$ d'(p)=4, d'(w)=5 , d'(z)=6$ and the degrees of the other points 
remain unchanged. Hence the new triangulation 
$T'= T \cup \{p\} $ has the parameters 
$$ (\alpha_3(T'),\alpha_4(T'), \alpha_5(T'), \alpha_6(T')) = (\alpha_3(T),\alpha_4(T),\alpha_5(T),\alpha_6(T)+1).$$ 
Thus $X$  is a productive as well as self-reproductive configuration, since 
in $T'$  the configuration $pywvx$ is isomorphic to the configuration $X$
(refer to Figure \ref{TyC}). Such a configuration $X$ will be
called a self-reproductive configuration of type $C.$
Note that each such configuration of type $C$ contains a pentagonal triangle and a triangle of degrees $4,5,5.$

\subsection{Type $D:$} Let $puvwxyz$ be any seven points in the triangulation 
$T$ and assume that $X= \{ p,u,v,w,x,y,z \}$ has six  triangles  
$ puv,pvw,pwx,pxy,pyz,pzu$ with degrees of the points given as 
$d(p)=6,d(w)=d(x)=d(y)=d(z)=5,$ $d(u)=4, d(v)=6.$ If we  delete 
the two triangles $puv ,pvw $ from $T,$ and  add one point 
$q$ and four new triangles $pqu,pqw,quv,qvw,$ then we get 
a new triangulation  $T'$ with the degrees of points as
$ d'(q)=4,\,d'(u)=5,\,d'(w)=6$ and the degrees of the other points   
remain unchanged. Hence $T'= T \cup \{p\} $ has the  parameters 
$$ (\alpha_3(T'),\alpha_4(T'), \alpha_5(T'), \alpha_6(T')) = (\alpha_3(T),\alpha_4(T),\alpha_5(T),\alpha_6(T)+1).$$ 
\begin{figure}[htb]
\[ \hspace {-0.4cm} 
\vcenter
{\xymatrix @M=0ex{
& {v(6)} \ar@{-}[ddr] \ar@{-}[dl]   \ar@{-}[r] &  {w(5)} \ar@{-}[dr] \ar@{-}[ddl]  & \\
{u(4)} \ar@{-}[dr] \ar@{-}[rrr]_(0.4){P(6)} &  &  & {x(5)} \ar@{-}[dl] \\
& {z(5)} \ar@{-}[r] _-{\bf{T}}&  {y(5)} &  }}
\quad  \longrightarrow \quad
\vcenter
{\xymatrix @M=0ex{
& {v(6)} \ar@{-}[ddr] \ar@{-}[dl]   \ar@{-}[r] &  {w(6)} \ar@{-}[dr] \ar@{-}[ddl] \ar@{-}[dll]_(.5){q} & \\
{u(5)} \ar@{-}[dr] \ar@{-}[rrr]_(0.4){P(6)} &  &  & {x(5)} \ar@{-}[dl] \\
& {z(5)} \ar@{-}[r] _-{\bf{T'}}&  {y(5)} &  }}
\]
\caption{Self-reproductive configuration of type $D$.}
\label{TyD}
\end{figure}
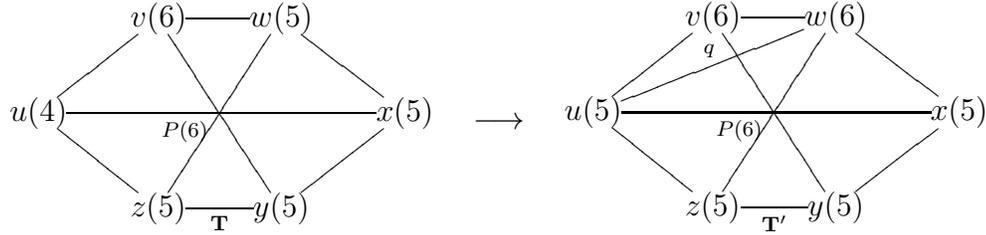

Thus $X$ is a productive as well as self-reproductive configuration, 
since in $T'$  
the configuration $pqwxyzu$ is isomorphic to the configuration $X$
(refer to Figure \ref{TyD}). Such a configuration $X$ will be
called a self-reproductive configuration of type $D.$\\

\subsection{Types $E_1$ and  $E_2:$} Let $puvwxyz$ be any seven points in a 
triangulation $T_1$ and assume that $X= \{ p,u,v,w,x,y,z \}$ contains 
six triangles $ puv,pvw,pwx,pxy,pyz,pzu$ with the degrees of the 
points given as 
$ d(u)=d(x)=4,\,d(v)=d(y)=5,\,d(p)=d(w)=d(z)=6$.    
Also let another triangulation $T_2$ be given with the set 
$X$ as above on seven points, where now the degrees are given as 
$  d(u)= d(w)=4, d(v)=d(z)=6,d(x)=d(y)=5,d(p) = 6.$ 
\begin{figure}[htb]
\[ \hspace {-0.4cm} 
\vcenter
{\xymatrix @M=0ex{
& {v(5)} \ar@{-}[ddr] \ar@{-}[dl]   \ar@{-}[r] &  {w(6)} \ar@{-}[dr] \ar@{-}[ddl]  & \\
{u(4)} \ar@{-}[dr] \ar@{-}[rrr]_(0.4){P(6)} &  &  & {x(4)} \ar@{-}[dl] \\
& {z(6)} \ar@{-}[r] _-{\bf{T_1}}&  {y(5)} &  }}
\quad  \longrightarrow \quad
\vcenter
{\xymatrix @M=0ex{
& {v(6)} \ar@{-}[ddr] \ar@{-}[dl] \ar@{-}[drr]^(.48){q}  \ar@{-}[r] &  {w(6)} \ar@{-}[dr] \ar@{-}[ddl]  & \\
{u(4)} \ar@{-}[dr] \ar@{-}[rrr]_(0.4){P(6)} &  &  & {x(5)} \ar@{-}[dl] \\
& {z(6)} \ar@{-}[r] _-{\bf{T'_1}}&  {y(5)} &  }}
\]
\end{figure}
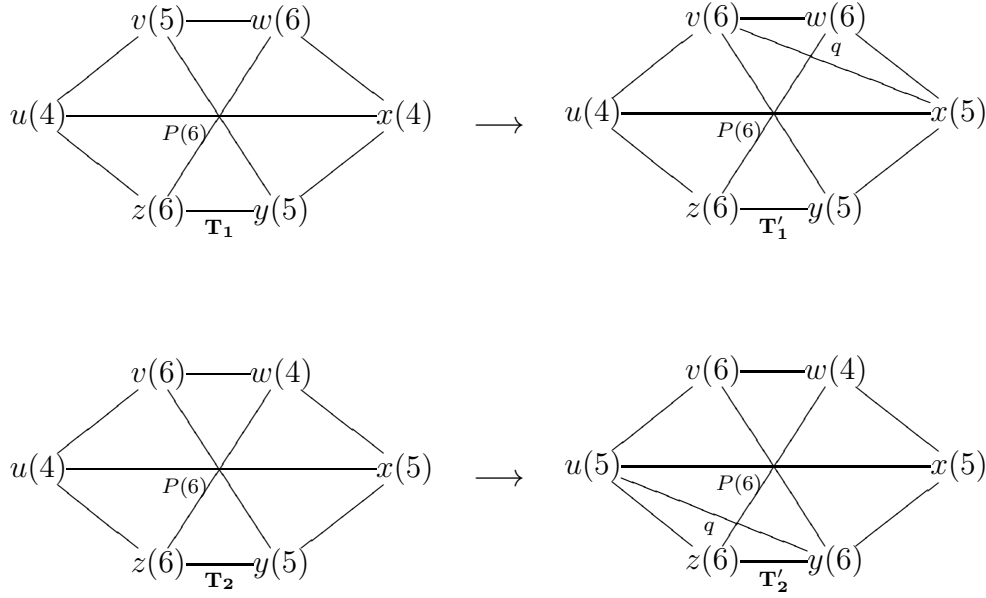
\begin{figure}[htb]
\[ \hspace {-0.4cm} 
\vcenter
{\xymatrix @M=0ex{
& {v(6)} \ar@{-}[ddr] \ar@{-}[dl]   \ar@{-}[r] &  {w(4)} \ar@{-}[dr] \ar@{-}[ddl]  & \\
{u(4)} \ar@{-}[dr] \ar@{-}[rrr]_(0.4){P(6)} &  &  & {x(5)} \ar@{-}[dl] \\
& {z(6)} \ar@{-}[r] _-{\bf{T_2}}&  {y(5)} &  }}
\quad  \longrightarrow \quad
\vcenter
{\xymatrix @M=0ex{
& {v(6)} \ar@{-}[ddr] \ar@{-}[dl]   \ar@{-}[r] &  {w(4)} \ar@{-}[dr] \ar@{-}[ddl]  & \\
{u(5)} \ar@{-}[dr] \ar@{-}[rrr]_(0.4){P(6)} \ar@{-}[drr]_(.5){q}&  &  & {x(5)} \ar@{-}[dl] \\
& {z(6)} \ar@{-}[r] _-{\bf{T'_2}}&  {y(6)} &  }}
\]
\caption{Self-reproductive configuration of types $E_1$ and $E_2$.}
\label{TyE1}
\end{figure}

If we delete the two triangles $pvw,pwx $ from $T_1,$ and  
add a point $q$ and four new triangles  $pqv,pqx,qvw,qxw$ 
then we get a new triangulation  $T_1'$ with the degrees 
$ d'(q)=4, d'(x)=5 , 
d'(v)=6$ and the degrees of the remaining points unchanged. 
Hence $T_1'= T_1 \cup \{q\} $ has the parameters  
$$ (\alpha_3(T_1'),\alpha_4(T_1'), \alpha_5(T_1'), \alpha_6(T_1')) = (\alpha_3(T_1),\alpha_4(T_1),\alpha_5(T_1),\alpha_6(T_1)+1).$$ 
Thus $X$ is a self-reproductive configuration because the 
configuration $pqvuzyx$ is isomorphic to the
configuration $X$ (in the triangulation $T_2$) (refer to Figure \ref{TyE1}). \\

If we delete the two triangles $pyz,pzu $ from $T_2,$ and  
add one point  $q$ and four new triangles  $pqy,pqu,qyz,qzu $ 
then we get another new triangulation $T_2'$ with the degrees 
$ d'(q)=4, d'(u)=5 , d'(y)=6$ and the degrees of the 
remaining points unchanged. Hence  $T_2'= T_2 \cup \{q\} $ 
has the parameters   
$$ (\alpha_3(T_2'),\alpha_4(T_2'), \alpha_5(T_2'), \alpha_6(T_2')) 
= (\alpha_3(T_2),\alpha_4(T_2),\alpha_5(T_2),\alpha_6(T_2)+1).$$ 
Thus $X$ is a self-reproductive configuration because 
the configuration $pquvwxy$ is isomorphic to the configuration $X$ 
in the  triangulation $T_1$ (refer to Figure \ref{TyE1}).\\

Hence from the given triangulation $T_1$   
we obtain a new triangulation $T_2'$ 
which has the parameters:
$$ (\alpha_3(T_2'),\alpha_4(T_2'), \alpha_5(T_2'), \alpha_6(T_2')) = (\alpha_3(T_1),\alpha_4(T_1),\alpha_5(T_1),\alpha_6(T_1)+2).$$ 
Since we can get the configuration $X$ (in the triangulation $T_1$) again
after two steps,    
$X$ is called a self-reproductive configuration of types
$E_1$ and $E_2$.

\subsection{Type $E_3:$} 
Let $puvwxyz$ be seven points in a 
triangulation $T_1$ and assume that  
$X= \{ p,u,v,w,x,y,z \}$ 
contains six triangles $ puv,pvw,pwx,pxy,pyz,pzu$ 
with the degrees of the points given as 
$d(u)=d(x)=4$, $d(v)=d(z)=5,\,d(w)=d(y)=d(p)=6$.     

Let a triangulation $T_2$ be given with the set 
$X$ as above on seven points, where now the degrees are given as 
$  d(u)= d(w)=4,\,d(v)=d(y)=d(p)= 6$, $ d(x)=d(z)=5$.   
Also let a triangulation $T_3$ be given with the set 
$X$ as above on seven points, where now the degrees are given as 
$d(v)=4$, $d(u)=d(x)= d(w)=d(z)=5$, $d(p)=d(y)=6$.    

If we delete the two triangles $pvw,pwx $ from $T_1,$ and  
add a point $q$ and four new triangles  $pqv,pqx,qvw,qxw$ 
then we get a new triangulation $T_1'$ with the degrees
$ d'(q)=4, d'(x)=5 , d'(v)=6$ and the degrees 
of the remaining points unchanged. 
Hence $T_1'= T_1 \cup \{q\} $ has the parameters  
$$ (\alpha_3(T_1'),\alpha_4(T_1'), \alpha_5(T_1'), \alpha_6(T_1')) 
= (\alpha_3(T_1),\alpha_4(T_1),\alpha_5(T_1),\alpha_6(T_1)+1).$$ 
Thus $X$ is a self-reproductive configuration because the 
set $pqvuzyx$ is a configuration which is isomorphic to the
configuration $X$ in the triangulation $T_2$ 
(refer to Figure \ref{TyE3}).
\begin{figure}[htb]
\[ \hspace {-0.4cm} 
\vcenter
{\xymatrix @M=0ex{
& {v(5)} \ar@{-}[ddr] \ar@{-}[dl]   \ar@{-}[r] &  {w(6)} \ar@{-}[dr] \ar@{-}[ddl]  & \\
{u(4)} \ar@{-}[dr] \ar@{-}[rrr]_(0.4){P(6)} &  &  & {x(4)} \ar@{-}[dl] \\
& {z(5)} \ar@{-}[r] _-{\bf{T_1}}&  {y(6)} &  }}
\quad  \longrightarrow \quad
\vcenter
{\xymatrix @M=0ex{
& {v(6)} \ar@{-}[ddr] \ar@{-}[dl] \ar@{-}[drr]_(.42){q}  \ar@{-}[r] &  {w(6)} \ar@{-}[dr] \ar@{-}[ddl]  & \\
{u(4)} \ar@{-}[dr] \ar@{-}[rrr]_(0.4){P(6)} &  &  & {x(5)} \ar@{-}[dl] \\
& {z(5)} \ar@{-}[r] _-{\bf{T'_1}}&  {y(6)} &  }} 
\]
\end{figure}
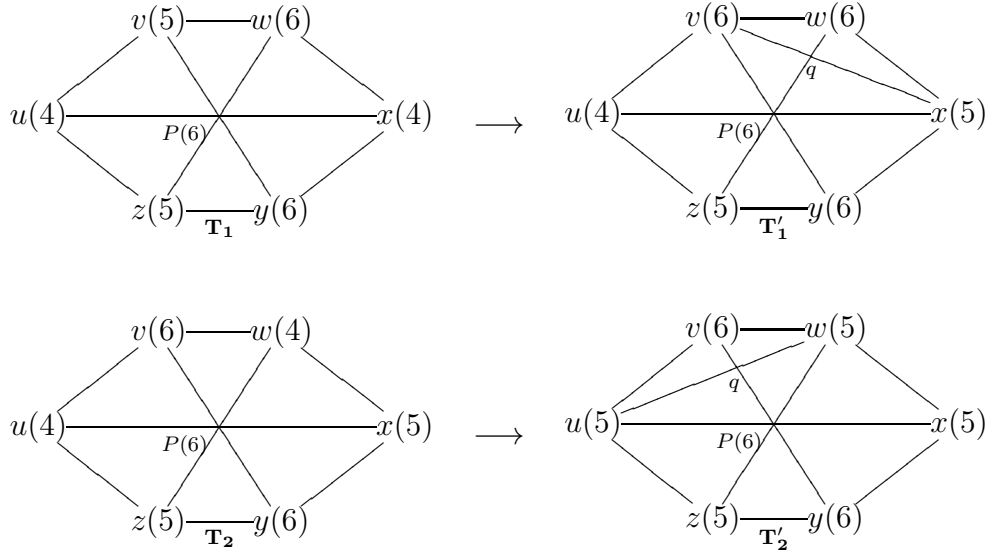
\begin{figure}[htb]
\[ \hspace {-0.4cm} 
\vcenter
{\xymatrix @M=0ex{
& {v(6)} \ar@{-}[ddr] \ar@{-}[dl]   \ar@{-}[r] &  {w(4)} \ar@{-}[dr] \ar@{-}[ddl]  & \\
{u(4)} \ar@{-}[dr] \ar@{-}[rrr]_(0.4){P(6)} &  &  & {x(5)} \ar@{-}[dl] \\
& {z(5)} \ar@{-}[r] _-{\bf{T_2}}&  {y(6)} &  }}
\quad  \longrightarrow \quad
\vcenter
{\xymatrix @M=0ex{
& {v(6)} \ar@{-}[ddr] \ar@{-}[dl]   \ar@{-}[r] &  {w(5)} \ar@{-}[dr] \ar@{-}[ddl] \ar@{-}[dll]^(.44){q} & \\
{u(5)} \ar@{-}[dr] \ar@{-}[rrr]_(0.4){P(6)} &  &  & {x(5)} \ar@{-}[dl] \\
& {z(5)} \ar@{-}[r] _-{\bf{T'_2}}&  {y(6)} &  }}
\]
\caption{Self-reproductive configuration of types $E_3$.}
\label{TyE3}
\end{figure}

If we delete the two triangles $pvw,puv $ from $T_2,$ and 
add one point  $q$ and four new triangles  $puq,quv,pwq,qvw $ 
then we get a new triangulation $T_2'$ with the degrees
$ d'(q)=4, d'(u)=5 , d'(w)=5$ and the degrees of the 
remaining points unchanged. Hence  $T_2'= T_2 \cup \{q\} $ 
has the parameters 
$$ (\alpha_3(T_2'),\alpha_4(T_2'), \alpha_5(T_2'), \alpha_6(T_2')) 
= (\alpha_3(T_2),\alpha_4(T_2)-1,\alpha_5(T_2)+2,\alpha_6(T_2)).$$ 
Thus $X$ is a self-reproductive configuration because 
the configuration $pqwxyzu$ is isomorphic to the configuration $X$ 
in the  triangulation $T_3$ (refer to Figure \ref{TyE3}).

\begin{figure}[htb]
\[ \hspace {-0.4cm} 
\vcenter
{\xymatrix @M=0ex{
& {v(4)} \ar@{-}[ddr] \ar@{-}[dl]   \ar@{-}[r] &  {w(5)} \ar@{-}[dr] \ar@{-}[ddl]  & \\
{u(5)} \ar@{-}[dr] \ar@{-}[rrr]_(0.4){P(6)} &  &  & {x(5)} \ar@{-}[dl] \\
& {z(5)} \ar@{-}[r] _-{\bf{T_3}}&  {y(6)} &  }}
\quad  \longrightarrow \quad
\vcenter
{\xymatrix @M=0ex{
& {v(4)} \ar@{-}[ddr] \ar@{-}[dl]   \ar@{-}[r] &  {w(5)} \ar@{-}[dr] \ar@{-}[ddl]  & \\
{u(5)} \ar@{-}[dr] \ar@{-}[rrr]_(0.4){P(6)} &  &  & {x(6)} \ar@{-}[dl] \ar@{-}[dll]_(.58){q}\\
& {z(6)} \ar@{-}[r] _-{\bf{T'_3}}&  {y(6)} &  }}
\]
\caption{Self-reproductive configuration of types $E_3$.}
\label{TyE32}
\end{figure}
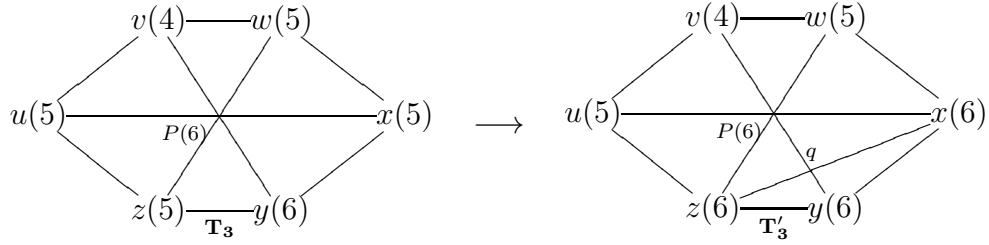

Again if we delete the two triangles $pzy,pyx $ from $T_3,$ and  
add one point  $q$ and four new triangles  $pqz,qyz,pqx,qyx $ 
then we get a new triangulation $T_3'$ with the degrees
$ d'(q)=4, d'(x)= d'(z)=6$ and the degrees of the 
remaining points unchanged. Hence  $T_2'= T_2 \cup \{q\} $ 
has the parameters 
$$(\alpha_3(T_3'),\alpha_4(T_3'), \alpha_5(T_3'), \alpha_6(T_3'))=
(\alpha_3(T_3),\alpha_4(T_3)+1,\alpha_5(T_3)-2,\alpha_6(T_3)+2).$$
Since the configuration $pqzuvwx$ is isomorphic to the configuration $X$ 
in the  triangulation $T_1$ (refer to Figure \ref{TyE32}), it follows that   
$X$ is a self-reproductive configuration. 

Hence from the given triangulation $T_1$    
we obtain a new triangulation $T_3'$    
which has the parameters:
$$ (\alpha_3(T_3'),\alpha_4(T_3'), \alpha_5(T_3'), \alpha_6(T_3')) = 
(\alpha_3(T_1),\alpha_4(T_1),\alpha_5(T_1),\alpha_6(T_1)+3).$$ 
Since we can get the configuration $X$ (in the triangulation $T_1$)
again after three steps, 
$X$ is called a self-reproductive configuration of type $E_3$.   

\subsection{Type $G:$} Let $puvwxyz$ be any seven points in a 
triangulation $T$ and assume that $X= \{ p,u,v,w,x,y,z \}$ 
contains six triangles 
$ puv,pvw,pwx,pxy,pyz,pzu$ with the degrees 
of the points given as 
$d(x)=d(v)=d(z)=d(p)=6$    
and $d(u)=d(w)=d(y)=4.$  
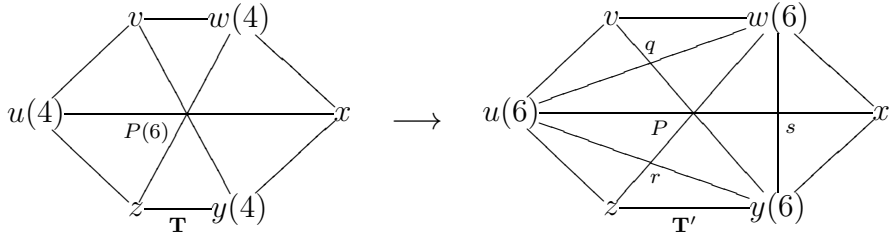
\begin{figure}[htb]
\[ \hspace {-0.1cm} 
\vcenter
{\xymatrix @M=0ex{
& {v} \ar@{-}[ddr] \ar@{-}[dl]   \ar@{-}[r] &  {w(4)} \ar@{-}[dr] \ar@{-}[ddl]  & \\
{u(4)} \ar@{-}[dr] \ar@{-}[rrr]_(0.36){P(6)} &  &  & {x} \ar@{-}[dl] \\
& {z} \ar@{-}[r] _-{\bf{T}}&  {y(4)} &  }}
\quad  \longrightarrow \quad
\vcenter
{\xymatrix @M=0ex{
& {v} \ar@{-}[ddrr] \ar@{-}[dl]   \ar@{-}[rr] & &  {w(6)} \ar@{-}[dr] \ar@{-}[ddll]  \ar@{-}[dlll]_(.46){q} 
\ar@{-}[dd]^(.58){s} & \\
{u(6)} \ar@{-}[dr] \ar@{-}[rrrr]_(0.4){P} \ar@{-}[drrr]_(.56){r} & &  &  & {x} \ar@{-}[dl] \\
& {z} \ar@{-}[rr] _-{\bf{T'}} & &  {y(6)} &  }}
\]
\caption{Self-reproductive configuration of types $G$.}
\label{TyG}
\end{figure}

If we delete the six triangles from $T$ and add  
three new points $q,r,s$ and twelve new triangles 
$upq,uqv,pqw,qvw,psw,psy,sxw,sxy,pur,pry,ruz,ryz $ 
then we get a triangulation $T'$ with 
the degrees $ d'(u)=d'(w)=d'(y)=6 , d'(q)=d'(r)=d'(s)=4$ and 
the degrees of the remaining points  
unchanged. Hence $T'= T \cup \{ q,r,s \} $ 
has the parameters 
$$ (\alpha_3(T'),\alpha_4(T'), \alpha_5(T'), \alpha_6(T')) =
(\alpha_3(T),\alpha_4(T),\alpha_5(T),\alpha_6(T)+3).$$  
Thus $X$ is a productive and self-reproductive configuration, since in $T'$  
the configuration $pqwsyru$ is isomorphic to the configuration $X$ 
(refer to Figure \ref{TyG}). Such a configuration $X$ will be called a 
self-reproductive configuration of type $G.$

\section{Fullering Constructions}

Let $T$ be a triangulation and let $a,b$ be two integers such that $3 \leq a < b $ and  
$ \alpha_d = 0 $ for $ d \not= a,b.$ The 
triangulation $T$  is called a  {\it fullering 
triangulation } if $ 6 \in \{ a,b \}.$  \\

\noindent
{\bf Example:} Note that 
$$T = \{ 156,157,167,256,258,268,357,358,378,467,468,478 \}$$ 
is a triangulation with  $ d(x) = 3$ for $ x=1,2,3,4$ and with  
$ d(x) = 6 $ for $ x=5,6,7,8.$ Hence  $T$ is a fullering triangulation. \\

\subsection{Face-fullering of Triangulations}
Let $T$ be a triangulation with face numbers $f_k(T)$ ($k=1,2,3$). That is,
$f_1(T)$ is the number of points, $f_2(T)$ the number of edges and  
$f_3(T)$ the number of triangles in $T.$ In this section we show that  
there exists a triangulation $ T_1$ with face numbers $f_k(T_1)$ ($k=1,2,3$)
and parameters as follows: 
\begin{equation} \label{NM1}
f_1(T_1) = f_1(T) + f_3(T); ~~ f_2(T_1) = 3f_2(T); ~~ f_3(T_1)= 3f_3(T);
\end{equation}
\begin{equation} \label{NM2}
\alpha_d(T_1) = \alpha_d(T) \mbox{ for } d \not= 6~;~~ 
\alpha_6(T_1) = \alpha_6(T) + f_3(T). 
\end{equation}
We also show that if T is a fullering (elliptic) triangulation, 
then $T_1$ is also a fullering (elliptic) triangulation. 

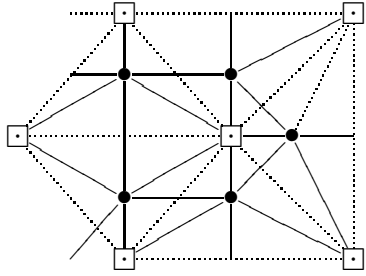
\begin{figure}[htb]
\hspace {4cm} 
\xymatrix @M=0ex@R=3.2ex@C=3ex{
& & {\boxdot} \ar@{.}[ddrr] \ar@{-}[d] \ar@{.}[ddll] \ar@{.}[l] \ar@{.}[rrrr] & & 
	 & & {\boxdot} \ar@{-}[dll] \ar@{.}[dddd] \ar@{.}[ddll]\ar@{.}[ddl]   \\
& & {\bullet} \ar@{-}[drr] \ar@{-}[dd] \ar@{-}[dll] \ar@{-}[l] \ar@{-}[rr] & & {\bullet}
	 \ar@{-}[d] \ar@{-}[u] \ar@{-}[dr]& &  \\
{\boxdot} \ar@{.}[ddrr] \ar@{-}[drr] \ar@{.}[rrrr] & & 
	 & & {\boxdot} \ar@{-}[dll] \ar@{-}[d] \ar@{-}[r] \ar@{.}[ddll] \ar@{.}[ddrr]  & 
	 {\bullet} \ar@{-}[r] \ar@{-}[ddr] \ar@{-}[dl]&\\
& & {\bullet} \ar@{-}[dl] \ar@{-}[d] \ar@{-}[rr] & & {\bullet} \ar@{-}[drr] \ar@{-}[d] \ar@{-}[dll] & &  \\
& & {\boxdot} \ar@{.}[rrrr] & & & & {\boxdot} }	
\caption{Face-fullering of a triangulation}
\label{sem-FFT}
\end{figure}

For the construction of $T_1$ we add one vertex point $ x_{\sigma}$ for each triangle $ \sigma \in T$ to the vertex points of $T:$ The set of vertices $V_1$ of  $T_1$ is accordingly $ V_1= V \cup \{ x_{\sigma}~ |~ \sigma \mbox{~is a triangle in~} T \}.$ 
The new triangles  for $T_1$ will be all 
$ \{ x , x_{\sigma}, x_{\tau} \} $ where 
$ | \sigma \cap \tau | = 2 $ and $x  \in  \sigma \cap \tau$
(see Figure \ref{sem-FFT}). 

We then obtain  $d_1(x) = d(x)$ for points $x \in T$ and 
$ d_1(x_{\sigma}) = 6$ for the new points $ x_{\sigma} \in T_1.$ 
Hence $ T_1$ possesses the property of being a fullering triangulation
if $T$ is a fullering triangulation and we also get (\ref{NM2}).   
The number of pairs  $ \{ \sigma, \tau \}$ which are neighbours (that is,  
the pairs $\{\sigma,\,\tau\}$ 
with the property $ | \sigma \cap \tau | = 2$)  
is the same as the integer $ f_2(T)$ (that is, the 
number of edges of $T$).  Hence we obtain from  (\ref{EQ0})
the equality $ f_3(T_1) = 2f_2(T) = 3 f_3(T)$,   
which verifies the third equation in (\ref{NM1}). 
Applying (\ref{EQ0}) once more we get  
$$ 2f_2(T_1) = 3f_3(T_1) = 9f_3(T) = 6f_2(T).$$
Cancelling the factor $2$ from both ends of this equation, we arrive 
at the second equation of (\ref{NM1}).  
We may also verify the first equation of (\ref{NM1}) by observing that
\begin{eqnarray*}
2 &=& f_1(T_1) - f_2(T_1) + f_3(T_1) = f_1(T_1) - 3f_2(T) + 3f_3(T)\\ 
 &=& f_1(T_1) - 3f_1(T) + 6. 
\end{eqnarray*}
This gives us 
$$ f_1(T_1)= 3f_1(T) -4 = f_1(T) + 2f_1(T) - 4 = f_1(T) + f_3(T).$$ 

This construction will be denoted by $ T_1 = FF(T)$ and 
will be called either the {\it face-fullering} or the {\it leapfrog} of 
$T.$  \\[2mm]

\subsection{Edge-fullering of Triangulations}\label{EFT}  
In this section we construct a   
triangulation $ T_2$ with face numbers $f_k(T_2)$ ($k=1,2,3$)
and parameters as follows: 

\begin{equation} \label{NM3}
f_1(T_2) = f_1(T) + f_2(T); ~~ f_2(T_2) = 4f_2(T); ~~ f_3(T_2)= 4f_3(T);
\end{equation}
\begin{equation} \label{NM4}
\alpha_d(T_2) = \alpha_d(T) \mbox{ for } d \not= 6~;~~ 
\alpha_6(T_2) = \alpha_6(T) + f_2(T). \end{equation}
As in the previous section it may be shown
that if  $T_2$ is a fullering (elliptic) triangulation, 
then $T_2$ is also a fullering (elliptic) triangulation. 

\begin{figure}[htb]
\hspace {4cm} 
\xymatrix @M=0ex@R=3.2ex@C=3ex{
& & & & & & & & \\
& & & {\boxdot} \ar@{.}[dr] \ar@{.}[u] \ar@{.}[dl] \ar@{.}[l] \ar@{.}[rr] & & {\bullet} \ar@{-}[dr] \ar@{-}[dl] 
	\ar@{-}[u] \ar@{-}[ur] \ar@{.}[rr] & & {\boxdot} \ar@{.}[dd] \ar@{.}[l] \ar@{.}[r] \ar@{.}[u] \ar@{.}[dl]&  \\
& & {\bullet} \ar@{.}[dl] \ar@{-}[l] \ar@{-}[ul] \ar@{-}[dr] \ar@{-}[rr] & & {\bullet} \ar@{.}[dr] \ar@{-}[dl] \ar@{-}[rr] 
	& &{\bullet} \ar@{-}[dd] \ar@{-}[dr] \ar@{.}[dl] & &  \\
& {\boxdot} \ar@{.}[dr] \ar@{.}[ul] \ar@{.}[dl] \ar@{.}[l] \ar@{.}[rr] & & {\bullet} \ar@{-}[dr] \ar@{-}[dl] 
	 \ar@{.}[rr] & & {\boxdot} \ar@{.}[dl] \ar@{.}[dr] & &{\bullet}\ar@{-}[dl] \ar@{.}[dd] \ar@{-}[dr] 
	 \ar@{-}[r]&  \\
& & {\bullet} \ar@{-}[dl] \ar@{-}[l] \ar@{-}[rr] \ar@{.}[dr] & & {\bullet} \ar@{-}[dr] \ar@{.}[dl] \ar@{-}[rr] 
	& &{\bullet} \ar@{-}[dl] \ar@{.}[dr] & &  \\
& & & {\boxdot} \ar@{.}[l] \ar@{.}[rr] \ar@{.}[dl]  & & {\bullet} \ar@{-}[dl] \ar@{-}[d] 
	\ar@{.}[rr] & & {\boxdot} \ar@{.}[d] \ar@{.}[dr] &  \\
& & & & & & & & }
\caption{Edge-fullering of a triangulation}
\label{sem-EFT}
\end{figure}
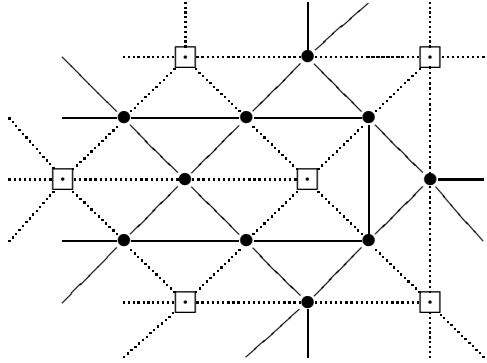

For the construction  of $T_2$ we add a point $x_e $ for each 
edge $e=xy$ of $T.$ Hence the point set  $V_2$ of $T_2$ is just  
$ V_2= V \cup VE,$ where $ VE = \{ x_e ~:~ e \in Edge(T) \}.$ The 
new triangles of  $T_2$ are  
all $ \{ x , x_e, x_f \} $ where $ e, f \subset \sigma $ and 
where $x \in  \sigma \cap \tau$ for some triangles $\sigma,~ \tau $ in $T,$
together with triangles 
$ \{ x_e, x_f, x_g \}$ where $\{ e,f,g \}$ are the edges of 
an old triangle (see Figure \ref{sem-EFT}). Note that 
for the old vertices $x \in V$ we get the same degree 
$ d_2(x)= d(x),$ and for the new vertices $x_e$ we get 
the degree $d_2(x_e)=6.$  
This construction will be denoted by $ T_2 = EF(T)$ and 
will be called the {\it edge-fullering} of $T.$ \\[2mm]

\subsection{Edge-fullering of Patches}
Let $P$ be a patch with face numbers $f_k(P)$ ($k=1,2,3$). That is,
$f_1(P)$ is the number of points, $f_2(P)$ the number of edges and  
$f_3(P)$ the number of triangles in $P$.   
By using the same construction method as for the edge-fullering of 
triangulations discussed in the previous section, 
it can be shown that 
there exists a patch $ P_2$ with face numbers $f_k(P_2)$ ($k=1,2,3$)
and parameters as follows: 

\begin{equation} \label{NM3}
f_1(P_2) = f_1(P) + f_2(P); \, f_2(P_2) = 2f_2(P)+3f_3(P);\,
 f_3(P_2)= 4f_3(P);\,\,  
\end{equation}
\begin{equation} \label{NM4}
\alpha_d(P_2) = \alpha_d(P) \mbox{ for } d \not= 6~;~~ 
\alpha_6(P_2) = \alpha_6(P) + f_2(P). \end{equation}
Moreover, it can be shown that 
the boundary length $b(P_2)$ of $P_2$ is twice that of $P,$
that is, $b(P_2)=2b(P).$

This construction will be denoted by $ P_2 = EF(P)$ 
and will be called the {\it edge-fullering} of $P.$ \\[2mm]

\section{Glueing of Patches} \label{glupat}
Let us consider two patches $P:(\alpha_3, \alpha_4, \alpha_5, \alpha_6)_b$
and $P':(\alpha_3', \alpha_4', \alpha_5', \alpha_6')_{b'}$. 
If both patches
have the same boundary length (that is, $b=b'$), 
then we can glue them
to get a new patch 
$(\alpha_3+ \alpha_3',\,\alpha_4+ \alpha_4',\,\alpha_5+ \alpha_5',\,
(\alpha_6+ \alpha_6')+bm),$
where $m \geq 0.$
The glueing can be done 
by using a circular strip of $2b$ triangles. 

We can also glue the patches by some special methods, which we describe below: 

\subsection{Method A}
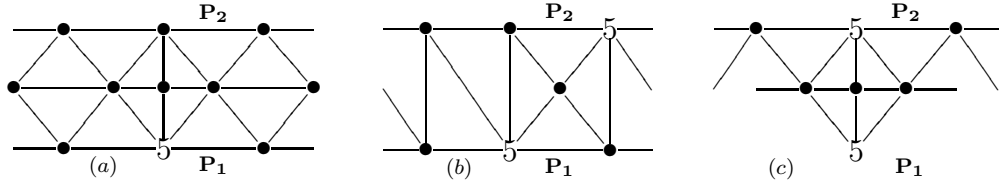
\begin{figure}[htb]
\[ \hspace{-0.6cm}
\vcenter 
{\xymatrix @M=0ex@R=3.2ex@C=2.5ex{
&  {\bullet} \ar@{-}[dr] \ar@{-}[dl] \ar@{-}[l] \ar@{-}[rr] & & {\bullet} \ar@{-}[dr] \ar@{-}[dl] \ar@{-}[d] 
	\ar@{-}[rr]^-{\bf{P_2}}	& & {\bullet} \ar@{-}[dr] \ar@{-}[dl] \ar@{-}[r] &  \\
{\bullet} \ar@{-}[dr] \ar@{-}[rr] & & {\bullet} \ar@{-}[dr] \ar@{-}[dl] \ar@{-}[r] 
	&{\bullet} \ar@{-}[d] \ar@{-}[r] & {\bullet} \ar@{-}[dr] \ar@{-}[dl] \ar@{-}[rr] &  &{\bullet}\ar@{-}[dl]\\
&  {\bullet} \ar@{-}[l] \ar@{-}[rr]_(.4){(a)} & & {5} \ar@{-}[rr]_-{\bf{P_1}} & & {\bullet} \ar@{-}[r] &} }
\quad \quad
\vcenter
{\xymatrix @M=0ex@R=3.2ex@C=2.5ex{
&  {\bullet} \ar@{-}[dd] \ar@{-}[ddrr] \ar@{-}[l] \ar@{-}[rr] & & {\bullet} \ar@{-}[dd] \ar@{-}[dr] 
	\ar@{-}[rr]^-{\bf{P_2}}	& & {5} \ar@{-}[dr] \ar@{-}[dl] \ar@{-}[r] \ar@{-}[dd]&  \\
 & & & & {\bullet} \ar@{-}[dr] \ar@{-}[dl]  &  & \\
&  {\bullet} \ar@{-}[l] \ar@{-}[ul] \ar@{-}[rr] _(.4){(b)}& & {5} \ar@{-}[rr]_-{\bf{P_1}} & & {\bullet} \ar@{-}[r] &}}
\quad \quad
\vcenter
{\xymatrix @M=0ex@R=3.2ex@C=2.5ex{
&  {\bullet} \ar@{-}[dr] \ar@{-}[dl] \ar@{-}[l] \ar@{-}[rr] & & {5} \ar@{-}[dr] \ar@{-}[dl] \ar@{-}[d] 
	\ar@{-}[rr]^-{\bf{P_2}}	& & {\bullet} \ar@{-}[dr] \ar@{-}[dl] \ar@{-}[r] &  \\
& & {\bullet} \ar@{-}[dr] \ar@{-}[l] \ar@{-}[r] 
	&{\bullet} \ar@{-}[d] \ar@{-}[r] & {\bullet} \ar@{-}[dl] \ar@{-}[r]&  &\\
&   & & {5} \ar@{}[rr]_-{\bf{P_1}} \ar@{}[lll]^-{(c)} & &  &} } 
\]
\caption{Glueing patches of type $(1,1,1,N)$ by different methods}
\label{con230-1}
\end{figure}

Consider the patches of types $(1,1,1,M)_b$ and $(1,1,1,N)_b$ 
which contain a point of 
degree five each on their boundary. If we glue these two patches by a belt, 
then we can get a triangulation of type $(2,2,2,M+N+b)$.    
This triangulation can be arranged in such a way that the two points of
degree $5$ are edge-related, that is,   
are connected by a sequence of edges  
(refer to Figure {\ref{con230-1}}(a)). 
Hence by inserting a new point of degree four 
between the two existing points of degree five, 
we obtain a triangulation of type $(2,3,0,M+N+b+2).$ 

\subsection{Method B}
Consider the patches of types $(1,1,1,M)_b$ and $(1,1,1,N)_b$ 
which contain a point of 
degree five each on their boundary. 
If we glue these two patches without a belt, 
we can get a triangulation of type $(2,2,2,M+N)$.   
Then by inserting a point of degree $4$ between the two points
of degree 5 on the boundary, 
we get a triangulation of type $(2,3,0,M+N+2)$     
(refer to Figure {\ref{con230-1}}(b)).  

\subsection{Method C}
Consider the patch $(1,1,1,M)_b$ with one point of degree five 
almost on the boundary 
and the patch $(1,1,1,N)_b$ which contains a point of degree five 
on its boundary.
If we glue these two patches, 
then we can get a triangulation of type $(2,2,2,M+N)$
(refer to Figure {\ref{con230-1}}(c)).
This triangulation can be arranged in such a way that 
the points of degree five 
are edge-related. 
Hence by inserting a point of degree four between them, 
we get a triangulation of type $(2,3,0,M+N+2).$

\subsection{Connected sum of triangulations}
Let $T_1$ be a triangulation with an inside triangle $a_1a_2a_3$ 
(that is, the triangle
$a_1a_2a_3$ is not on the boundary of $T_1$) with 
$d_i$ the degree of $a_i$ $(i=1,2,3)$.  
Let $T_2$ be another triangulation whose 
outside edges form a triangle $b_1b_2b_3$ 
with $d_i'$ the degree of $b_i$ $(i=1,2,3)$. 
Then we can glue these two triangulations
by identifying the points $a_i$  and $b_i$ 
so that the degrees of these identified
points become $\overline{d}_i = d_i + d'_i - 2$ ~$(i=1,2,3).$ 
This construction will be called the
connected sum of the triangulations $T_1$ and $T_2.$

\vspace{0.3cm}
\section{Triangulations of type $(0,0,12)$} \label{tri0012} 
In this section we shall see that the 
triangulation $(0,0,12,N)$ exists for all values of $N$ except for $N=1.$
To construct triangulations of this type, we will consider the following $6$
cases:~$0(6),1(6),2(6),3(6),4(6),5(6).$
\begin{table}[h]
\begin{center}
\begin{tabular}{|c|l|} \hline 
$N=a_6 $  &  Triangulation \\ \hline
0  & $~ \{ 123,~134,~145,~156,~126,~238,~267,~278,~349,~389,45A,49A, $  \\
   & $ ~~~56B,~5AB,~67B,~78C,~7BC,~89C,~9AC,~ABC ~ \}$    \\ \hline 
2  & $ ~\{ 123,~134,~145,~156,~126,~238,~287,~267,~349,~389,45B,49A, $  \\
   & $ ~~~4AB,~56D,~5BD,~67D,~89C,~78C,~9AC,~ACE,~7CE, $ \\
   & $ ~~BDE, 7DE,ABE, \}$    \\ \hline 
3  & $~ \{ 123,~126,~134,~145,156,23C,~2C7,~267,~3AB,~3BC,~34A,  $ \\ 
   & $ ~~ 458,~48A,~589,~569,~679,~89F,~8AE,~8EF,~79F,~ABE,~BCD,$     \\
   & $ ~~7CD,~EFD,~7DF,~BED \} $ \\ \hline 
4 & $ ~ \{ 123,~126,~16C,~13D,~1CD,~234,~245,~256,~34E,~3ED,~49E, $   \\
   & $  ~~489,~458,~567,~578,~67B,~6CB,~7AB,~78A,~89A,~ABG,~A9F,$   \\
   & $  ~~BCG,~CDG,~EDF,~9EF,~DFG,~AGF  \} $  \\  \hline 

\end{tabular}  
\end{center}
\caption{Some triangulations of type $(0,0,12)$}
\label{tab0012}
\end{table}

\subsection{$0(6)$ type}
There exists a triangulation $T=(0,0,12,0+6m)$ for all $m \geq 1.$
This can be constructed as follows:\\
Glueing of the patches $(0,0,6,3)_6$ and $(0,0,6,3)_6$ gives a triangulation 
$T_1=(0,0,12,6+6m)$ for all $m \geq 0.$
The construction of the patch $(0,0,6,3)_6$ 
is given in Table $2.4$.
The triangulation $T_1=(0,0,12,6+6m)$ for all $m \geq 0,$
is the same as the triangulation 
$T=(0,0,12,0+6m)$ for all $m \geq 1$.    
Hence, the triangulation $T$ gives all the values of the type $0(6)$ except the 
initial value $N=0.$ For $N=0,$ the triangulation
is given in Table \ref{tab0012} above.

\subsection{$1(6)$ type}
There exists a triangulation $T=(0,0,12,1+6m)$ for all $m \geq 1.$
This can be constructed as follows:\\
Glueing of the patches $(0,0,6,3)_6$ and $(0,0,6,4)_6$ gives a triangulation 
$T_1=(0,0,12,7+6m)$ for all $m \geq 0.$
The construction of the patches $(0,0,6,3)_6$ and $(0,0,6,4)_6$ 
are given in Table $2.4.$
The triangulation  $T_1=(0,0,12,7+6m)$ for all $m \geq 0$
is the same as the triangulation  
$T=(0,0,12,1+6m)$ for all values of $m \geq 1$.   
Hence, the triangulation $T$ gives all the values of type $1(6)$ except the 
initial value of $N=1.$ For $N=1,$ the triangulation does not exist as 
has been shown by Gr\"unbaum and Motzkin \cite{GM} (see also \cite{GR}). 

\subsection{$2(6)$ type}
There exists a triangulation $T=(0,0,12,2+6m)$ for all values of $m \geq 0.$
This can be constructed as follows:\\
Glueing $m$ belts together with the triangulation $(0,0,12,2)$
gives the triangulation $T:(0,0,12,2+6m)$ for all $m \geq 0.$
The construction of the triangulation $(0,0,12,2)$ is given in 
Table \ref{tab0012} above.   

\subsection{$3(6)$ type}
There exists a triangulation $T=(0,0,12,3+6m)$ for all $m \geq 0.$
This can be constructed as follows: \\
Glueing of the patches $(0,0,6,1)_6$ and $(0,0,6,2)_6$ gives the triangulation 
$T=(0,0,12,3+6m)$ for all $m \geq 0.$
The construction of the patches $(0,0,6,1)_6$ 
and $(0,0,6,2)_6$ are given in Table 2.4.   

\subsection{$4(6)$ type}
There exists a triangulation $T=(0,0,12,4+6m)$ for all of $m \geq 0.$
This can be constructed as follows:\\
Glueing of the patches $(0,0,6,1)_6$ and $(0,0,6,3)_6$ gives the triangulation 
$T=(0,0,12,4+6m)$ for all $m \geq 0.$
The construction of the patches $(0,0,6,1)_6$ and $(0,0,6,3)_6$ 
are given in Table $2.4$.  

\subsection{$5(6)$ type}
There exists a triangulation $T=(0,0,12,5+6m)$ for all $m \geq 0.$
This can be constructed as follows:\\
Glueing of the patches $(0,0,6,2)_6$ and $(0,0,6,3)_6$ gives the triangulation 
$T=(0,0,12,5+6m)$ for all $m \geq 0.$
The construction of the patches $(0,0,6,2)_6$ and $(0,0,6,3)_6$ 
are given in Table $2.4.$   

\section{Triangulations of type $(0,1,10)$}  
In this section we show that the triangulation $(0,1,10,N)$ 
exists for all values of $N$ except for $N=0,1$. 
This can be constructed as follows:\\
The triangulation $T_1:(0,1,10,3)$ as given in
Table \ref{tab0110} below 
contains a self-reproductive configuration  of 
type $D$, namely, $ X = \{ 1,3,4,6,7,D,E \}$.    
This self-reproductive configuration gives the
triangulation $T:(0,1,10,N)$ for all  $N \geq 3$ from the triangulation
$T_1$.  

\begin{table}[h]
\begin{center}
\begin{tabular}{|c|l|} \hline 
$N=a_6 $  &  Triangulation \\ \hline
2  & $ ~\{ 134,~137,~145,~167,~15D,~16D,~258,~25D,~269,~26D,28B,29B, $  \\
   & $ ~~~34A,~37C,~3AC,~458,~48A,~679,~79C,~8AB,~9BC,~ABC \} $ \\  \hline 
    
3  & $~ \{~134,~137,~14E,~1ED,16D,167,~258,~25D,~269,~26D,~28B,  $ \\ 
   & $ ~~ 29B,~34A,~37C,~3AC,~458,~45E,~48A,~679,~79C,~8AB,~9BC, $ \\
   & $ ~~~ABC,~5DE \} $   \\  \hline 

\end{tabular}  
\end{center}
\caption{Some triangulations of type $(0,1,10)$}
\label{tab0110}
\end{table}
For $N=2$, the triangulation $(0,1,10,2)$  
is given in Table \ref{tab0110}. For $N=0,1$, 
the triangulations cannot be constructed as has been shown 
by Eberhard \cite{Eber} 
and Br\"uckner \cite{bruck}.
This has been stated by Gr\"unbaum \cite{GR} and is easy to prove.

\section{Triangulations of type $(0,2,8)$}  
There exists a triangulation $T:(0,2,8,N)$ for all  $N \geq 0.$  
This can be constructed as follows:\\
The triangulation $T_1:(0,2,8,0)$ 
as given in Table \ref{tab028} below    
contains 
a self-reproductive configuration  of 
type $C$, namely, $X = \{ 2,3,4,5,9 \}$.    
This self-reproductive configuration gives the
triangulation $T:(0,2,8,m)$ for all  $m \geq 0$ from the triangulation
$T_1$.    

\begin{table}[h]
\begin{center}
\begin{tabular}{|c|l|} \hline 
$N=a_6 $  &  Triangulation \\ \hline
0  & $ ~\{ 126,~129,~149,~167,~147,~239,~235,~256,~349,~358,~348,~478, $  \\
   & $ ~~~56A,~58A,~67A,~78A \} $ \\  \hline 
    
1  & $~ \{~123,~12B,~13A,~178,17A,18B,~235,~259,~29B,~345,~34A,  $ \\ 
   & $ ~~ 456,~467,~47A,~569,~678,~689,~89B \} $   \\  \hline 
\end{tabular}  
\end{center}
\caption{Some triangulations of type $(0,2,8)$}
\label{tab028}
\end{table}

\section{Triangulations of type $(0,3,6)$} \label{Tri036} 
There exists a triangulation $T:(0,3,6,N)$ for all  $N \geq 0.$  
This can be constructed as follows:\\
The triangulation $T_1:(0,3,6,0)$ as given in
Table \ref{tab036} below    
contains a self-reproductive configuration  of 
type $C$, namely, $X = \{ 2,3,4,6,8 \}$.    
This self-reproductive configuration gives the
triangulation $T:(0,3,6,N)$ for all  $N \geq 0$ from the triangulation
$T_1$.    
Note that the triangulation $T:(0,3,6,N)$ can be constructed
so that it contains a pentagonal triangle.

\begin{table}[h]
\begin{center}
\begin{tabular}{|c|l|} \hline 
$N=a_6 $  &  Triangulation \\ \hline
0  & $ ~\{ 123,~126,~167,~234,~246,~345,~458,~468,~678,~139,~179,~359, $  \\
   & $ ~~~589,~789 \} $ \\  \hline 
    
1  & $~ \{~123,~126,~167,~234,~24A,~26A,~345,~458,~48A,~678,~68A,  $ \\ 
   & $ ~~ 139,~179,~359,~589,~789, \} $   \\  \hline 
   
\end{tabular}  
\end{center}
\caption{Some triangulations of type $(0,3,6)$}
\label{tab036}
\end{table}

\section{Triangulations of type $(0,4,4)$}  
In this section we show that the triangulation $(0,4,4,N)$ 
exists for all values of $N$ except for $N=1.$
We first show the existence of
a triangulation $T$ of the type $(0,4,4,N)$ for all  $N \geq 4.$  
This can be constructed as follows:\\

\begin{table}[h]
\begin{center}
\begin{tabular}{|c|l|} \hline 
$N=a_6 $  &  Triangulation \\ \hline
0  & $~ \{ 125,~156,~127,~168,~178,~235,~237,~345,~347,~456,478,468 ~ \} $  \\ \hline 
    
2  & $ ~\{ 125,~156,~127,~178,~168,~259,~29A,~27A,~37A,~39A,~359,~345, $  \\
   & $ ~~~347,~456,~478,~468 ~ \} $ \\ \hline 
  
3  & $~ \{~125,~156,~127,~178,~168,~259,~29A,~27A,~37A,~39A,~359,  $ \\ 
   & $ ~~ 35B,~34B,~347,~46B,~478,~468,~56B ~ \} $ \\ \hline 

4  & $~ \{~125,~15C,~16C, ~127,~178,~168,~259,~29A,~27A,~37A,~39A, $ \\ 
   & $ ~~~359,~35B,~34B,~347,~46B,~478,~468,~5BC, ~6BC ~ \} $ \\ \hline 

\end{tabular}  
\end{center}
\caption{Some triangulations of type $(0,4,4)$}
\label{tab044}
\end{table}

The triangulation $T_1:(0,4,4,4)$ 
as given in Table \ref{tab044}
contains 
a self-reproductive configuration  of 
types $E_1$ and $E_2$, namely, $ X = \{ 1,2,3,5,9,B,C \} .$ 
This self-reproductive configuration gives the
triangulation $T:(0,4,4,N)$ for all  $N \geq 4$ from the triangulation
$T_1$.    
 
For $N=0,2,3,$ the construction of a triangulation of type 
$(0,4,4,N)$ is also given in Table \ref{tab044}. For $N=1$,   
a triangulation of the type $(0,4,4,1)$ does not exist. This can be
proven as follows:

Let $X=\{abcdefg \}$ with the triangles $abc,acd,ade,aef,afg,abg$ and  
degree $d(a)=6$. 

\noindent
Case (I): If the degrees of the points $d,\,e$ are five, 
then there may be three new points $h,\,i,\,j$ and five triangles
$cdh,\,dhi,\,dei,\,eij,\,efj$. This would give the
number of points $f_1>9$, which contradicts the fact that the
triangulation $(0,4,4,1)$ contains only $9$ points.
\smallskip

\noindent
Case (II): If the degrees of the points $d,\,e$ are
$d(d)=5,\,d(e)=4$, then there may be two new points
$h$ and $j$ and four triangles $cdh,\,dhj,\,dej,\,efj$. Then the
points $h,\,j$ must have either degrees $4$ or $5$. In both cases
there may be one or more edges from these two points which are
not connected with the other points of this construction.
This would give the number of points in the
construction as $f_1>9$, which is again a contradiction.
 
Hence, the triangulation $(0,4,4,1)$ cannot be constructed.

\section{Triangulations of type $(0,5,2)$}  
We will show that the triangulation $(0,5,2,N)$ 
exists for all values of $N$ except for $N=1$.  
To construct triangulations of this type we will consider two cases: 
(1) $N$ even, ~(2)~$N$ odd. 

\begin{table}[h]
\begin{center}
\begin{tabular}{|c|l|} \hline 
$N=a_6 $  &  Triangulation \\ \hline
0  & $~ \{ 134,~145,~156,~167,~137,~234,~237,~245,~256,~267 ~ \} $  \\ \hline 
    
2  & $ ~\{ 137,~138,~156,~167,~158,~239,~237,~267,~259,~256,~348,~349, $  \\
   & $ ~~~459,~458 ~ \} $ \\ \hline 
  
3  & $~ \{~128,~12A,~145,~148,~156,~16A,~237,~27A,~239,~289,~349,  $ \\ 
   & $ ~~ 356,~367,~345,~489,~67A ~ \} $ \\ \hline 

4  & $~ \{~136,~138,~167, ~189,~19A,~17A,~246,~267,~248,~289,~29B,  $ \\ 
   & $ ~~~27B,~356,~358,~456,~458,~9AB,~7AB ~ \} $ \\ \hline 
   
\end{tabular}  
\end{center}
\caption{Some triangulations of type $(0,5,2)$}
\label{tab052}
\end{table}

\subsection{$N$ even }
There exists a triangulation $T:(0,5,2,2m)$ for all  $m \geq 0.$  
This can be constructed as follows:

The triangulation $T_1:(0,5,2,0)$
as given in Table \ref{tab052}   
contains a self-reproductive configuration  of 
type $B_1$, namely, $ X = \{ 1,2,5,6,7 \}$.    
This self-reproductive configuration gives the
triangulation $T:(0,5,2,2m)$ for all  $m \geq 0$ from the triangulation
$T_1$.  

\subsection*{$N$ Odd}
In this case we will consider three types: $1(6),3(6),5(6).$ 

\subsection{$1(6)$ type}
There exists a triangulation $T:(0,5,2,1+6m)$ for all  $m \geq 1.$  
This can be constructed as follows:\\
The triangulation $T_1:(0,5,2,4)$ 
as given in Table \ref{tab052} above     
contains a self-reproductive configuration  of 
type $E_3,$ namely, $ X = \{ 1,3,6,7,8,9,A \}$.   
This self-reproductive configuration gives a
triangulation $T_2:(0,5,2,4+3h)$ for $h \geq 0$ from the triangulation
$T_1.$ 
Note that the triangulation $T_1:(0,5,2,4)$ 
can also be constructed as in the above $N$ 
even case with $m=2$.    

If we let $h=1+2m$ for all $m \geq 0$ in the triangulation 
$T_2:(0,5,2,4+3h),$ then we can get a triangulation $(0,5,2,7+6m)$ 
for all  $m \geq 0,$ which is the same as the triangulation
$T:(0,5,2,1+6m)$ for all  $m \geq 1.$ \\
 
The triangulation $T$ gives all the values of type $1(6),$ except the initial
value $N=1.$ For $N=1,$ the triangulation cannot be constructed. This
can be proven as follows:\\ 

Let $X=\{abcdefg \}$ with the triangles $abc,acd,ade,aef,afg,abg$ and  
degree $d(a)=6$. Suppose without loss of generality that the point
$d$ has degree $5$. There must be a new point (say $h$) which is 
connected to $d$. Otherwise, since $d(d)=5$, we would have the triangle
$def$ or $dcb$. But then $e$ or $c$ would have degree $3$ 
which contradicts the
fact that the triangulation we want to construct 
cannot have any points of   
degree $3$.  
We may assume, without loss of generality, that the points $d,\,h,\,e$
form a triangle. Since $d$ has degree $5$, the points $d,\,c,\,h$
cannot form a triangle. Thus there must exist at least one more new
point (say $i$) so that we would have the triangles $dhi$, $dic$;
hence making the degree of $d$ to be $5$.      
But then there would be altogether at least $9$ points in this
triangulation which contradicts the fact that the 
triangulation 
$(0,\,5,\,2,\,1)$ has only got $8$ points.   

Hence, the triangulation $(0,5,2,1)$ cannot be constructed.

\subsection{$3(6)$ type}
There exists a triangulation $T:(0,5,2,3+6m)$ for all  $m \geq 0.$  
This can be constructed as follows:\\
The triangulation $T_1:(0,5,2,3)$ 
as given in Table \ref{tab052}
above contains 
a self-reproductive configuration  of 
type $G$, namely, $ X = \{ 1,2,4,5,6,8,A \}$.    
This self-reproductive configuration gives the
triangulation $T_2:(0,5,2,3+3h)$ for $h \geq 0$ from the triangulation
$T_1.$ 
If we let $h=2m$ for all $m \geq 0$  in the triangulation 
$T_2:(0,5,2,3+3h),$ then  we can get the triangulation $T:(0,5,2,3+6m)$ 
for all  $m \geq 0.$ 
The triangulation $T$ gives all the values of type $3(6).$

\subsection{$5(6)$ type}
There exists a triangulation $T:(0,5,2,5+6m)$ for all  $m \geq 0.$  
This can be constructed as follows:\\
The triangulation $T_1:(0,5,2,2)$ 
as given in Table \ref{tab052} above     
contains a self-reproductive configuration of
type $G$, namely, $ X = \{ 1,2,3,4,7,8,9 \} .$ 
This self-reproductive configuration gives the
triangulation $T_2:(0,5,2,2+3h)$ for $h \geq 0,$ from the triangulation
$T_1$. 
If we let $h=1+2m$ for all $m \geq 0$  in the triangulation 
$T_2:(0,5,2,2+3h),$ then  we can get the triangulation $T:(0,5,2,5+6m)$ 
for all  $m \geq 0.$ 
The triangulation $T$ gives all the values of type $5(6).$

\section{Triangulations of type $(0,6,0)$}  \label{Tri060}
In this section we show that the triangulation $(0,6,0,N)$ exists for all values of $N$ except for $N=1.$ 
To construct this type of triangulation we will consider two cases: 
(1) $N$ even, ~(2)~$N$ odd.

\begin{table}[h]
\begin{center}
\begin{tabular}{|c|l|} \hline 
$N=a_6 $  &  Triangulation \\ \hline
0  & $~ \{ 123,~126,~135,~156,~234,~246,~345,~456 ~ \} $  \\ \hline 
    
2  & $ ~\{ 127,~157,~237,~537,~238,~538,~248,~548,~126,~246,~156,~456 ~ \} $ \\ \hline
  
3  & $~ \{~123,~126,~135,~156,~234,~246,~345,~469,~479,~457,~578,  $ \\ 
   & $ ~~ 568,~689,~789 \} $ \\  \hline
4  & $~ \{~174,~179,~16A,~169,~26A,~27A,~237,~238,~248,~246,~456,  $ \\ 
   & $ ~~ 537,~538,~548,~569,~579 ~ \} $ \\ \hline 
  
\end{tabular} 
\end{center}
\caption{Some triangulations of type $(0,6,0)$}
\label{tab060}
\end{table}

\subsection{$N$ even}
There exists a triangulation $T:(0,6,0,2m)$ for all  $m \geq 0.$  
This can be constructed as follows:\\
The triangulation $T_1:(0,6,0,0)$ 
as given in Table \ref{tab060} above   
contains a self-reproductive configuration  of 
type $B_1$, namely, $ X = \{ 1,2,3,4,5 \}$.    
This self-reproductive configuration gives the
triangulation $T:(0,6,0,2m)$ for all  $m \geq 0$.

\subsection*{$N$ odd}
In this case we will consider three types: $1(6),3(6),5(6).$ 

\subsection{$1(6)$ type}
There exists a triangulation $T:(0,6,0,1+6m)$ for all  $m \geq 1.$  
This can be constructed as follows:\\
The triangulation $T_1:(0,6,0,4)$ 
as given in Table \ref{tab060} above       
contains a self-reproductive configuration of 
type $G$, namely, $ X = \{ 2,3,4,6,7,8,A \}$.    
This self-reproductive configuration gives the
triangulation $T_2:(0,6,0,4+3h)$ for all  $h \geq 0$ from the triangulation
$T_1$.  
If we let $h=1+2m$ for all $m \geq 0$ in the triangulation 
$T_2:(0,6,0,4+3h),$ then we can get a triangulation $(0,6,0,7+6m)$ 
for all  $m \geq 0,$ which is the same as the triangulation
$T:(0,6,0,1+6m)$ for all  $m \geq 1.$   

The triangulation $T$ gives all the values of type $1(6)$ except the initial
value of $N=1.$ For $N=1,$ the triangulation does not exist as 
has been shown by Gr\"unbaum and Motzkin \cite{GM} (see also \cite{GR}). 

\subsection{$3(6)$ type} \label{sbtri060}
There exists a triangulation $T:(0,6,0,3+6m)$ for all  $m \geq 0.$  
This can be constructed as follows:\\
Using the construction of the connected sum of the triangulations 
$(0,6,0,i)$ and $(0,6,0,i)$ or $(0,6,0,i+3)$ for all $i= 3h $ where $h \geq 0,$
we get a triangulation $T_1:(0,6,0,3h)$ for all $h \geq 0.$  

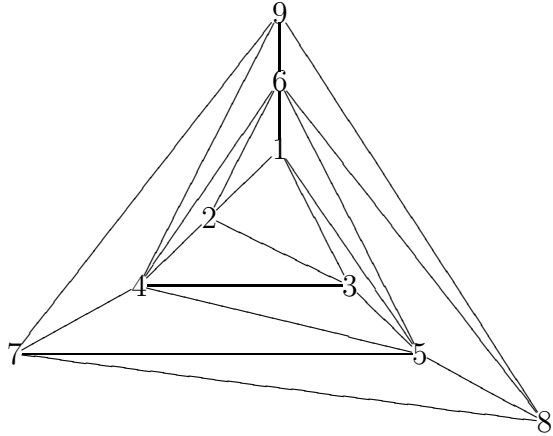
\begin{figure}[htb]
\hspace {3.5cm} 
\xymatrix @M=0ex@R=3.5ex@C=4ex{
& && & 9\ar@{-}[d] \ar@{-}[dddddllll] \ar@{-}[ddddll] \ar@{-}[ddddddrrrr]& & & & \\
& && & 6 \ar@{-}[d] \ar@{-}[ddl] \ar@{-}[dddll] \ar@{-}[dddddrrrr] \ar@{-}[ddddrr]& & & &\\ 
& & & & 1 \ar@{-}[dl] \ar@{-}[ddr] \ar@{-}[dddrr] & & & &\\
& & & 2\ar@{-}[dl]  \ar@{-}[drr]  & & & &\\
& & 4 \ar@{-}[dll]  \ar@{-}[drrrr]  \ar@{-}[rrr]&&      & 3 \ar@{-}[dr]& & &\\
7 \ar@{-}[drrrrrrrr] \ar@{-}[rrrrrr]&& & & & &5\ar@{-}[drr] & &\\ 
 && & & & & & &8}
\caption{Triangulation of type $(0,6,0,3)$.}
\label{Tr0603}
\end{figure}

For example, we can construct the triangulation $(0,6,0,3)$ by the connected 
sum of the triangulations $(0,6,0,0)$ and $(0,6,0,0).$ This is 
illustrated in Figure \ref{Tr0603}. 

If we let $h=1+2m$ for all $m \geq 0$ in the triangulation 
$T_1:(0,6,0,3h),$ then we can get the 
triangulation $T:(0,6,0,3+6m)$ for all  $m \geq 0.$ 
The triangulation $T$ gives all the values of type $3(6).$

\subsection{$5(6)$ type}
There exists a triangulation $T:(0,6,0,5+6m)$ for all  $m \geq 0.$  
This can be constructed as follows:\\
The triangulation $T_1:(0,6,0,2)$ 
as given in Table \ref{tab060} above  
contains a self-reproductive configuration  of
type $G$, namely, $ X = \{ 1,2,3,4,6,7,8 \}$.    
This self-reproductive configuration gives a
triangulation $T_2:(0,6,0,2+3h)$ for all  $h \geq 0,$ from the triangulation
$T_1.$ 
If we let $h=1+2m$ for all $m \geq 0$ in the triangulation 
$T_2:(0,6,0,2+3h),$ then we can get the triangulation $T:(0,6,0,5+6m)$ 
for all  $m \geq 0.$ 
The triangulation $T$ gives all the values of type $5(6).$

\section{Triangulations of type $(1,0,9)$}  
We will show that the triangulation $(1,0,9,N)$ exists for all values 
of $N$ except for $N=0,1,2$ and possibly $N=4$.   

\begin{table}[h]
\begin{center}
\begin{tabular}{|c|l|} \hline 
$N=a_6 $  &  Triangulation \\ \hline
3  & $~ \{ 128,~12D,~139,~13D,~18B,~19B,~235,~256,~268,~23D,~357, $ \\  
   & $~~~379,~456,~457,~46A,~47C,~4AC,~68A,~79C,~8AB,~9BC,$ \\
   & $~~~ABC ~\} $ \\ \hline 
6  & $~ \{~123,~12G,~16G,~134,~145,156,23C,~2C7,~267,~26G,~3AB,  $ \\ 
   & $ ~~~3BC,~34A,~458,~48A,~589,~569,~679,~89F,~8AE,~8EF,~79F,$  \\
   & $ ~~~ABE,~BCD,~7CD,~EFD,~7DF,~BED \} $ \\ \hline 
  
\end{tabular} 
\end{center}
\caption{Some triangulations of type $(1,0,9)$}
\label{tab109}
\end{table}

First we show that the 
triangulation $T:(1,0,9,N)$ exists for all $N \geq 6.$ This can be  
constructed as follows:\\
By glueing the patches $(1,0,3,i)_6$ for $i=5,6,7,8,9$
and $(0,0,6,j)_6$ for  $j=1,2,3,4,6 ,$  we can get a triangulation
$T_1:(1,0,9,k+6m)$ for  $k=6,7,8,9,10,11$ and for $m \geq 0$ (here
$k=i+j$). The triangulation $T_1$ is the same as the triangulation
$T:(1,0,9,N)$ for all $N \geq 6.$ The construction of the patches $(1,0,3,i)_6$ for $i=5,6,7,8,9$ 
and $(0,0,6,j)_6$ for $j=1,2,3,4,6$ are given in Table $2.4.$ 

The triangulation $T$ gives all the values of $ N,$ 
except the initial values of 
$N=0,1,2,3,4,5.$ For $N=0,1,2,$ the triangulations cannot be
constructed as has been shown by Eberhard \cite{Eber} 
and Br\"uckner \cite{bruck}.
This has been stated by Gr\"unbaum \cite{GR} and is easy to prove. 
For $N=3,$ the triangulation $(1,0,9,3)$ is given in Table \ref{tab109} above. 
It is not known yet whether the triangulation
$(1,0,9,4)$ exists.
For $N=5,$ the triangulation $(1,0,9,5)$ can be constructed as follows: \\

Glueing of the patches $(0,0,6,0)_5$ and $(1,0,3,5)_5$ gives a triangulation 
$T_1:(1,0,9,5+5m)$ for all $m \geq 0.$
The construction of the patches $(0,0,6,0)_5$ and $(1,0,3,5)_5$ are given in Table \ref{penbd}. 
If we let $m=0$ in the triangulation $T_1:(1,0,9,5+5m),$ 
then we can get the triangulation $(1,0,9,5).$

\section{Triangulations of type $(1,1,7)$}  \label{Tri117}

There exists a triangulation $T:(1,1,7,N)$ for all  $N \geq 2.$  
This can be constructed as follows:\\
The triangulation $T_1:(1,1,7,2)$ 
as given in Table \ref{tab117} below     
contains a self-reproductive configuration  of 
type $C$, namely, $ X = \{ 2,3,9,A,B \}$.    
This self-reproductive configuration gives the
triangulation $T:(1,1,7,N)$ for all  $N \geq 2$ from the triangulation
$T_1$.

\begin{table}[h]
\begin{center}
\begin{tabular}{|c|l|} \hline 
$N=a_6 $  &  Triangulation \\ \hline
2  & $~ \{ 124,~12A,~1AB,~147,~17B,~245,~235,~23A,~359,~3AB,~39B,~456, $ \\  
   & $~~~468,~478,~568,~589,~789,~79B ~\} $ \\ \hline 

3  & $~ \{ 124,~12A,~1AB,~147,~17B,~245,~235,~23A,~359,~3AC,~39C,~456, $ \\  
   & $~~~468,~478,~568,~589,~789,~79B, ~9BC,~ABC ~\} $ \\ \hline 
   
\end{tabular}  
\end{center}
\caption{Some triangulations of type $(1,1,7)$}
\label{tab117}
\end{table}

The triangulation $(1,1,7,N)$ for $N=0,1$ cannot be
constructed. This has been shown by Eberhard \cite{Eber} 
and Br\"uckner \cite{bruck}, and is stated 
in Gr\"unbaum \cite{GR}.

\section{Triangulations of type $(1,2,5)$} \label{Tri125}
In this section we show that there exists a triangulation $T:(1,2,5,N)$ for all  $N \geq 1.$  
This can be constructed as follows:\\
The triangulation $T_1:(1,2,5,1)$ 
as given in Table \ref{tab125} below  
contains a self-reproductive configuration  of
type $C$, namely, $ X = \{ 1,2,4,5,6 \}$.   
This self-reproductive configuration gives the
triangulation $T:(1,2,5,N)$ for all  $N \geq 1$ from the triangulation
$T_1.$ 

The triangulation $(1,2,5,0)$ cannot be
constructed as has been shown by Eberhard \cite{Eber} 
and Br\"uckner \cite{bruck}, and stated in 
\cite{GR}.

\begin{table}[h]
\begin{center}
\begin{tabular}{|c|l|} \hline 
$N=a_6 $  &  Triangulation \\ \hline

1  & $~ \{ 125,~156,~128,~189,~169,~235,~237,~278,~345,~348,~378,~456, $ \\  
   & $~~~489,~469 ~\} $ \\ \hline
2  & $~ \{ 12A,~16A,~128,~189,~169,~235,~237,~25A,~278,~345,~348,~378,$ \\
   & $~~~456,~489, 469, ~56A ~\} $ \\ \hline 
   
\end{tabular} 
\end{center}
\caption{Some triangulations of type $(1,2,5)$}
\label{tab125}
\end{table}

\section{Triangulations of type $(1,3,3)$} 
In this section we show that the triangulation 
$(1,3,3,N)$ exists for all values of $N\ge 0$.    

\begin{table}[h]
\begin{center}
\begin{tabular}{|c|l|} \hline 
$N=a_6 $  &  Triangulation \\ \hline
0  & $~ \{ 126,~127,~134,~137,~146,~235,~237,~256,~345,~456 ~\} $ \\ \hline 

1  & $~ \{ 125,~128,~137,~138,~156,~167,~234,~238,~245,~347,~456,~467 ~\} $ \\ \hline 
 
2  & $~ \{ 128,~189,~149,~134,~137,~127,~286,~256,~237,~235,~345,~456,$ \\
   & $~~~469,~689 ~\} $ \\ \hline
   
\end{tabular}  
\end{center}
\caption{Some triangulations of type $(1,3,3)$}
\label{tab133}
\end{table}

First we show that
there exists a triangulation $T:(1,3,3,N)$ for all $N \geq 3.$
This can be constructed as follows:\\
The construction method of the 
triangulations of type $(0,3,6)$ as given in Section
\ref{Tri036} shows that there exists a 
triangulation $T_1:(0,3,6,m)$ for all $m \geq 0$ which contains a 
pentagonal triangle. By using the mutant configuration of type $M_1$ on the triangulation $T_1:(0,3,6,m)$ 
for all $ m \geq 0,$ we can get a triangulation 
$(1,3,3,m+3)$ for all $m \geq 0,$
which is the same as the triangulation $T:(1,3,3,N)$ for all $N \geq 3.$ 

For $N=0,1,2,$ the triangulations $(1,3,3,N)$
are given in Table \ref{tab133} above.

\section{Triangulations of type $(1,4,1)$}  
We show that the triangulation $(1,4,1,N)$ exists for all values of $N$ except for $N=0,1.$
To construct triangulations of this type, we will consider two cases: 
(1) $N$ odd, ~(2)~$N$ even.

\begin{table}[h]
\begin{center}
\begin{tabular}{|c|l|} \hline 
$N=a_6 $  &  Triangulation \\ \hline
2  & $~ \{ 123,~124,~148,~138,~235,~245,~358,~456,~467,~478,~568,~678 ~\} $ \\ \hline 

3  & $~ \{ 125,~128,~134,~145,~189,~139,~256,~267,~238,~237,~389,~347, $ \\
   & $ ~~~467,~456 ~\} $ \\ \hline 
 
4  & $~ \{ 123,~134,~127,~145,~156,~167,~234,~249,~289,~278,~45A,~49A, $ \\
   & $ ~~~56A,~67A,~78A,~89A ~ \}  $  \\ \hline 

10 & $~ \{~129,~12C,~15C, ~15D,~189,~18D,~239,~23G,~24C,~24G,~346,  $ \\ 
   & $ ~~~34G,~36A,~39A,~456,~45C,~567,~57D,~67E,~78B,~7BE,~89F, $ \\
   & $ ~~~8BF,~9AF,~ABE,~ABF,~78D,~6AE ~ \} $ \\ \hline 
   
\end{tabular}  
\end{center}
\caption{Some triangulations of type $(1,4,1)$}
\label{tab141}
\end{table}

\subsection{$N$ odd}
There exists a triangulation $T:(1,4,1,1+2m)$ for all $m \geq 1.$
This can be constructed as follows:\\
The triangulation $T_1:(1,4,1,3)$ 
as given in Table \ref{tab141} above  
contains a self-reproductive configuration of   
type $B_1$, namely, $ X = \{ 2,4,5,6,7 \}$.   
This self-reproductive configuration gives a
triangulation $T_2:(1,4,1,3+2m)$ for all  $m \geq 0$ from the triangulation
$T_1$. The triangulation $T_2:(1,4,1,3+2m)$ for $m\ge 0$
is the same as the
triangulation $T:(1,4,1,1+2m)$ for $m \geq 1$.   

The triangulation $T$ gives all the odd values of $N,$ except the initial 
value of $N=1$. For $N=1,$ the triangulation cannot be
constructed (see Eberhard \cite{Eber} 
and Br\"uckner \cite{bruck}, or 
Gr\"unbaum \cite{GR}).

\subsection*{$N$ even}
We will divide the even case into $3$ types:~$0(6),2(6),4(6).$

\subsection{$0(6)$ type}
There exists a triangulation $T:(1,4,1,0+6m)$ for all $m \geq 1.$
This can be constructed as follows:\\
The triangulation $T_1:(1,4,1,3)$ 
as given in Table \ref{tab141} above    
contains a self-reproductive configuration  of
type $G$, namely, $ X = \{ 1,2,3,5,6,7,8 \}$.    
This self-reproductive configuration gives a
triangulation $T_2:(1,4,1,3+3h)$ for all  $h \geq 0$ from the triangulation
$T_1$. If we let $h=1+2m$ in the triangulation $T_2:(1,4,1,3+3h),$
then we can get a triangulation $(1,4,1,6+6m)$ for all $m \geq 0,$
which is the same as the triangulation $T:(1,4,1,0+6m)$ for all $m \geq 1.$

The triangulation $T$ gives all the values of type $0(6)$, except the 
initial value of $N=0.$ For $N=0,$ the triangulation cannot be constructed
(see Eberhard \cite{Eber} 
and Br\"uckner \cite{bruck}, or 
Gr\"unbaum \cite{GR}).

\subsection{$2(6)$ type}
There exists a triangulation $T:(1,4,1,2+6m)$ for all $m \geq 0.$
This can be constructed as follows:\\
The triangulation $T_1:(1,4,1,2)$ 
as given in Table \ref{tab141}
contains a triangle with degrees $4,4,4$.    
Therefore, we can construct a triangulation $T_2:(1,4,1,2+3h)$ for all 
$h \geq 0,$ by taking the connected sum with the copies of  the 
triangulation $T_3:(0,6,0,3)$ 
(as given in Table \ref{tab060}).   
If we let $h=2m$ for all $m \geq 0$ 
in the triangulation $T_2: (1,4,1,2+3h)$, then   
we can get the triangulation $T:(1,4,1,2+6m)$ for all $m \geq 0$.   
The triangulation $T$ gives all the values of type $2(6).$

\subsection{$4(6)$ type}
There exists a triangulation $T:(1,4,1,4+6m)$ for all $m \geq 2.$
This can be constructed as follows:\\
Glueing of the patches $(1,1,1,13)_6$ and $(0,3,0,3)_6$ gives a triangulation 
$T_1:(1,4,1,16+6m)$ for all $m \geq 0.$
The construction of the patch $(1,1,1,13)_6$ is given in  
Table $2.4.$ The patch $(0,3,0,3)_6$
can be constructed from Theorem \ref{T1}(ii) with the substitution of $h=2.$ 
The triangulation $T_1:(1,4,1,16+6m)$ for all $m \geq 0$
is the same as the triangulation  $T:(1,4,1,4+6m)$ for all $m \geq 2.$ \\
The triangulation $T$ gives all the values of type $4(6),$ except the 
initial values of $N=4,10.$ For $N=4,10,$ the triangulations
are given in Table \ref{tab141} above.

\section{Triangulations of type $(2,0,6)$}  
In this section we show that the triangulation $(2,0,6,N)$ exists for all values of $N$ except for $N=1.$

\begin{table}[h]
\begin{center}
\begin{tabular}{|c|l|} \hline 
$N=a_6 $  &  Triangulation \\ \hline
0  & $~ \{ 124,~127,~135,~137,~145,~236,~237,~246,~356,~458,~468,~568 ~\} $ \\ \hline 

2  & $~ \{ 123,~126,~134,~145,~159,~169,~267,~278,~23A,~28A,~348,~38A, $ \\
   & $ ~~~457,~478,~567,~569 ~\} $ \\ \hline 
 
3  & $~ \{ 146,~149,~156,~157,~178,~189,~236,~238,~256,~25A,~278,~27A, $ \\
   & $ ~~~346,~34B,~389,~39B, ~49B,~57A ~ \}  $  \\ \hline 
    
\end{tabular} 
\end{center}
\caption{Some triangulations of type $(2,0,6)$}
\label{tab206}
\end{table}

First we show that there exists a triangulation $T=(2,0,6,N)$ for all $N \geq 4.$
This can be constructed as follows:\\
The construction method of the triangulations of type $(1,1,7)$ as given 
in Section \ref{Tri117} shows that there exists a triangulation $T_1:(1,1,7,m)$ 
for all $m \geq 2$ which contains a 
self-reproductive configuration of type $C;$ hence it contains
a triangle of degrees $4,5,5.$ If we insert a new point $x$ 
of degree $d(x)=3$ into this triangle, then we can get a 
triangulation $(2,0,6,N)$ for all $N \geq 4.$ 

For $N=0,2,3,$ the construction of the triangulations 
$(2,0,6,N)$ are given in Table \ref{tab206} above. 
For $N=1,$ the triangulation cannot be constructed. 
We can prove this as follows: 

Let $X=\{abcdefg \}$ with the triangles $abc,acd,ade,aef,afg,abg$ 
and let the degrees of the points $a,d$ be $d(a)=6,d(d)=5$.   
We consider the following two cases: \\
Case (I): Suppose that the degrees of the points $e,\,c$ are three.
Then the degree of the remaining points $b,\,g,\,f$ must be five.
Suppose that $bcd,\,def$ form triangles. Since $d(b)=5$, we cannot
have the triangle $dbg$. Thus we introduce a new point $h$ so that
$dhb$ forms a triangle. Since $d$ has degree $5$, we cannot form any
more triangles involving $d$. In order for the construction to be
of type $(2,0,6,1)$, $h$ must have degree $5$ and we must have
another point $i$ so that $hib$ forms a triangle. But then 
$gbi$ will also form a triangle contradicting the fact that 
$b$ has degree $5$.  

Suppose without loss of generality that only $bcd$ forms a triangle.
Since $d$ has degree $5$, there must be a new point $h$ so that 
$dhb$ and $edh$ form triangles. Since $d(b)=5$, we must also have the 
triangle $hbg$. Since $d(g)=5$, there must be a new point $i$ so that 
$hig$ and $fgi$ form triangles. If we were to introduce new points
in this construction, we would not be able to get a triangulation
of the type $(2,0,6,1)$ which has only $9$ points. But then $f$
would only have degree $3$ which is also a contradiction.

Suppose now that $bcd$ and $def$ do not form triangles. Then $d$ must
be connected to two new points $h,\,i$ so that the triangles
$dhi,\,edh,\,cdi$ are formed. In order for the triangulation
to be of type $(2,0,6,1)$, the points $h$ and $i$ must have degree $5$.  
Thus more triangles involving $h$ and $i$ must be formed.
However, any other triangle formed involving $h$ or $i$ will result
in $c$ or $e$ having degree more than $3$; a contradiction.
\smallskip
 
\noindent
Case (II): Suppose that $d(e)=5$ and $d(c)=3$. Then there must be two new
points $h,\,i$ and the triangles $edh,\,hei,\,ief$.   
Since $d(c)=3$, we must then have the triangles $bcd,\,hdb$.
In order for the triangulation to be of the type $(2,0,6,1)$,
$b$ must have degree $5$ and we must have the triangle $hbg$. 
Since we do not want to introduce anymore new points in this
construction and $h$ has degree $5$, the triangle $hif$ must be
formed. Thus $d(i)=3$ and it follows that $f,\,g$ must have degree $5$.
But in order for $f,\,g$ to have degree $5$, we have to introduce at least
another new point; thus giving rise to a triangulation not of the
form $(2,0,6,1)$.  
\smallskip

Hence, the triangulation $(2,0,6,1)$ cannot be constructed.
\medskip

\section{Triangulations of type $(2,1,4)$}  
We show that the triangulation $(2,1,4,N)$ exists for all values of $N$ except for $N=0.$

\begin{table}[h]
\begin{center}
\begin{tabular}{|c|l|} \hline 
$N=a_6 $  &  Triangulation \\ \hline
1  & $~ \{ 123,~126,~136,~234,~248,~268,~345,~357,~367,~457,~478,~678,~\} $ \\ \hline 
    
2  & $~ \{ 126,~128,~134,~137,~178,~146,~235,~237,~256,~278,~345,~459, $ \\
   & $ ~~~469,~569 ~\} $ \\ \hline 
3  & $~ \{ 125,~128,~135,~139,~178,~179,~246,~24A,~25A,~268,~345,~346, $ \\
   & $ ~~~367,~379,~45A,~678 ~ \}  $  \\ \hline 
     
\end{tabular}  
\end{center}
\caption{Some triangulations of type $(2,1,4)$}
\label{tab214}
\end{table}

First we show that there exists a triangulation $T:(2,1,4,N)$ for all $N \geq 3.$
This can be constructed as follows:\\
The triangulation $T_1:(1,2,5,N)$ for all $N \geq 1$ 
contains a self-reproductive configuration of type $C$ and hence contains 
a triangle of degrees $4,5,5$. 
The construction of the triangulation
$T_1:(1,2,5,N)$ for all $N \geq 1$ is given in Section
\ref{Tri125}. If we apply the mutant configuration of type $M_2$ to  
the triangle of degrees $4,5,5$
in the triangulation $T_1,$ then we can get a 
triangulation $(2,1,4,N+2)$ for all $N \geq 1,$  which
is the same as the triangulation $T:(2,1,4,N)$ for all $N \geq 3.$ 

For $N=1,2$, the construction of the triangulations 
$(2,1,4,N)$ are given in Table \ref{tab214} above. 
For $N=0,$ the triangulation cannot be constructed
(see Eberhard \cite{Eber} 
and Br\"uckner \cite{bruck}, or 
Gr\"unbaum \cite{GR}).

\section{Triangulations of type $(2,2,2)$}  
In this section we show that there exists a triangulation $T:(2,2,2,N)$ for all $N \geq 0.$
This can be constructed as follows:\\
The triangulation $T_1:(2,2,2,0)$ 
as given in Table \ref{tab222} below     
contains a self-reproductive configuration of type  $A$, namely,  
$X = \{ 1,2,3 \}$. This self-reproductive configuration 
gives the triangulation $T:(2,2,2,N)$ for all $N \geq 0$.

\begin{table}[h]
\begin{center}
\begin{tabular}{|c|l|} \hline 
$N=a_6 $  &  Triangulation \\ \hline
0  & $~ \{ 123,~124,~146,~135,~156,~234,~345,~456 ~\} $ \\ \hline 

1  & $~ \{ 127,~137,~237,~124,~146,~135,~156,~234,~345,~456 ~\} $ \\ \hline 

2  & $~ \{ 127,~137,~278,~378,~238,~124,~146,~135,~156,~234,~345,~456 ~ \}  $  \\ \hline 
   
\end{tabular}  
\end{center}
\caption{Some triangulations of type $(2,2,2)$}
\label{tab222}
\end{table}

\section{Triangulations of type $(2,3,0)$}  
In this section we show that the triangulation $(2,3,0,N)$ 
exists for all values of $N$ except for  $N=1$ and
possibly $N=3,7,15,31$.

\begin{table}[h]
\begin{center}
\begin{tabular}{|c|l|} \hline 
$N=a_6 $  &  Triangulation \\ \hline
0  & $~ \{ 123,~124,~134,~235,~245,~345 ~\} $ \\ \hline 

13  & $~ \{ 123,~12G,~134,~14K,~15G,~15K,~23E,~267,~26G,~27E,~389, $ \\
    & $~~~38E,~349,~49A,~45A,~45K,~56A,~56G,~67D,~6AD,~78B, $ \\
    & $~~~78E,~7BD,~89F,~8BF,~9AC,~9CF,~ACD,~BCF,~BCH,$ \\
    & $~~~BDH,~CDH ~ \}  $  \\ \hline 
  
\end{tabular}  
\end{center}
\caption{Some triangulations of type $(2,3,0)$}
\label{tab230}
\end{table}

To construct triangulations of this type we will consider two cases: 
(1) $N$ even, ~(2)~$N$ odd.

\subsection{$N$ even} \label{ty562002}
There exists a triangulation $T:(2,3,0,2m)$ for all $m \geq 0.$
This can be constructed as follows:\\
The triangulation $T_1:(2,3,0,0)$ 
as given in Table \ref{tab230} above   
contains a self-reproductive configuration of type $B_1$,   
namely, $X = \{ 1,2,3,4,5 \}.$
This self-reproductive configuration 
gives the triangulation $T:(2,3,0,2m)$ for all $m \geq 0$.

\subsection*{$N$ odd}
In this case we will consider three types: $1(6),3(6),5(6).$ 

\subsection{$1(6)$ type} 

To construct this type, we will consider the triangulations of the 
form    
$(2,3,0,$ $X+Ym)$ with $Y=24,30,48,60$ and all possible values of $X.$ 
The least common multiple of $24,30,48,60$ is $240.$ Therefore the values
of $N$ in this $1(6)$ type will be 
$N=1,7,13,19,25,31,37,43,\dots,217,223,229,235.$ 
For the construction of these types we will consider the following cases:

\subsubsection{Case 1}
There exists a triangulation $T=(2,3,0,1+24m)$ for all 
$m \geq 1.$ This can be constructed as follows:\\
Using the edge-fullering construction on the triangulation $(2,3,0,2h),$ 
for all $h \geq 0,$ we can get a triangulation $T_1=(2,3,0,9+8h)$ for all $h \geq 0,$
that is, $EF(2,3,0,2h) = (2,3,0,9+8h).$ 
The construction of the triangulation $(2,3,0,2h)$ for all
$h \geq 0$ is given in Section \ref{ty562002}.
If we substitute $h=2+3m$ for $m \geq 0$ in $T_1=(2,3,0,9+8h)$,   
then we can get a triangulation $(2,3,0,25+24m)$ for $m \geq 0,$
which is the same as the triangulation $T=(2,3,0,1+24m)$ for  
$m \geq 1$.    

The triangulation $T$ gives all the values of the case $1(24)$ except the 
initial value of $N=1.$ For $N=1$ the triangulation cannot be constructed.
We can prove this as follows:\\ 

Let $X=\{abcdefg \}$ with the triangles $abc,acd,ade,aef,afg,abg$ 
so that the degree of the 
point $a$ is $d(a)=6$. This construction   
contains $f_1 \geq 7$ points which is a contradiction, 
since the triangulation 
$(2,3,0,1)$ contains only $6$ points.
Hence, the triangulation $(2,3,0,1)$ cannot be constructed.

\subsubsection{Case 2} 
There exists a triangulation $T=(2,3,0,7+30m)$ for all $m \geq 1.$
This can be constructed as follows:\\
Glueing of the patches $(1,1,1,4+3h)_{5+2h}$ and $(1,1,1,11+7h)_{5+2h}$ by 
the method $C$ (as explained in Section \ref{glupat}) gives a triangulation 
 $T_1=(2,3,0,17+10h)$ for all $h \geq 0.$  
The patch $(1,1,1,4+3h)_{5+2h}$ is of type $A$ with $k=3$,
and $(1,1,1,11+7h)_{5+2h}$ is of type $D$ with $k=5$ (types $A$ and $D$ are
as given in Section \ref{con111}).  

If we let $h=2+3m$ for $m \geq 0$ in the triangulation
$T_1=(2,3,0,17+10h),$ then we can get a triangulation 
$(2,3,0,37+30m)$ for $m \geq 0,$ 
which is the same as the triangulation $T=(2,3,0,7+30m)$ for $m\geq 1$.   
The triangulation $T$ gives all the values of the case $7(30),$ except the 
initial value of $N=7.$ We have not been able to construct a triangulation
of the type $(2,3,0,7).$

\subsubsection{Case 3}
There exists a triangulation $T=(2,3,0,19+24m)$ for all $m \geq 1.$
This can be constructed as follows:\\
Glueing of the patches $(1,1,1,11+3h)_{9+2h}$ and $(1,1,1,14+5h)_{9+2h}$ by 
the method $C$ (as explained in Section \ref{glupat}) gives a triangulation 
$T_1=(2,3,0,27+8h)$ for all $h \geq 0.$
The patch $(1,1,1,11+3h)_{9+2h}$ can be constructed 
from the patch $(1,1,1,5+3h')_{5+2h'}$ with the substitution of $h'=2+h$
where the patch $(1,1,1,5+3h')_{5+2h'}$ is of type C with $k=5$ and 
$(1,1,1,14+5h)_{9+2h}$ is of type A with $k=5$ (types $A$ and $C$ are
as given in Section \ref{con111}).    

If we let $h=2+3m$ for $m \geq 0$ in the triangulation
$T_1=(2,3,0,27+8h),$ then we can get a triangulation 
$(2,3,0,43+24m)$ for $m \geq 0$,   
which is the same as the triangulation $T=(2,3,0,19+24m)$ for $m \geq 1$.   
The triangulation $T$ gives all the values of the case $19(24),$ except the 
initial value of $N=19.$ 
For $N=19,$ the triangulation $(2,3,0,19)$ can be constructed 
by glueing of the patches $(1,1,1,5)_5$ and $(1,1,1,7)_5$ 
via method $A$ (as explained in Section \ref{glupat}).
The construction of the patches $(1,1,1,5)_5$ and $(1,1,1,7)_5$ 
are given in Table \ref{penbd}.

\subsubsection{Case 4}
There exists a triangulation $T=(2,3,0,13+48m)$ for all $m \geq 1.$
This can be constructed as follows:\\
Glueing of the patch $(1,1,1,18+7h)_{7+2h}$ with itself by 
using the method $A$ (as explained in Section \ref{glupat}) gives a triangulation 
$T_1=(2,3,0,45+16h)$ for all $h \geq 0.$
The patch  $(1,1,1,18+7h)_{7+2h}$ is of type C with $k=7$
where type $C$ is as given in Section \ref{con111}.   
If we let $h=1+3m$ for $m \geq 0$ in the triangulation 
$T_1=(2,3,0,45+16h),$ then we can get a triangulation 
$(2,3,0,61+48m)$ for $m \geq 0,$ 
which is the same as the triangulation $T=(2,3,0,13+48m)$ for $m\geq 1$.   
The triangulation $T$ gives all the values of the case $13(48),$ except the 
initial value of $N=13.$ 
For $N=13,$ the triangulation $(2,3,0,13)$ is given in Table \ref{tab230}. 

\subsubsection{Case 5}
There exists a triangulation $T=(2,3,0,31+60m)$ for all $m \geq 1.$
This can be constructed as follows:\\
Glueing of the patches $(1,1,1,25+7h)_{9+2h}$ and $(1,1,1,35+11h)_{9+2h},$ by 
using the method $A$ (as explained in Section \ref{glupat}) 
gives a triangulation $T_1=(2,3,0,71+20h)$ for all $h \geq 0$.   
The patch $(1,1,1,25+7h)_{9+2h}$ is of type B with $k=9$ and 
$(1,1,1,35+11h)_{9+2h}$ is of type C with $k=9$ (types $B$ and $C$ 
are as given   
in Section \ref{con111}).  
If we let $h=1+3m$ for $m \geq 0$ in the triangulation
$T_1=(2,3,0,71+20h),$ then we can get a triangulation of type
$(2,3,0,91+60m)$ for  $m \geq 0$,   
which is the same as the triangulation $T=(2,3,0,31+60m)$ for $m\geq 1$.   
The triangulation $T$ gives all the values of the case $31(60),$ except the 
initial value of $N=31.$ As for the case $N=7,$ we have not been able to
construct a triangulation of the type $(2,3,0,31).$

\subsubsection{Case 6}
There exists a triangulation $T=(2,3,0,43+60m)$ for all $m \geq 1.$
This can be constructed as follows:\\
Glueing of the patches $(1,1,1,14+3h)_{11+2h}$ and $(1,1,1,56+15h)_{11+2h},$ by 
the method $A$ (as explained in Section \ref{glupat}) gives a triangulation 
$T_1=(2,3,0,83+20h)$ for all $h \geq 0.$
The patch $(1,1,1,14+3h)_{11+2h}$ can be constructed 
from the patch $(1,1,1,5+3h')_{5+2h'}$ with the substitution of $h'=3+h.$
The patches $(1,1,1,5+3h')_{5+2h'}$ and $(1,1,1,56+15h)_{11+2h}$
are of type $C$ with $k=5$ and $k=11$ respectively
(type $C$ is as given in Section  \ref{con111}).  
If we let $h=1+3m$ for $m \geq 0$ in the triangulation
$T_1=(2,3,0,83+20h),$ then we can get a triangulation  
$(2,3,0,103+60m)$ for $m \geq 0,$ 
which is the same as the triangulation 
$T=(2,3,0,43+60m)$ for $m \geq 1$.  
The triangulation $T$ gives all the values of the case $43(60),$ except the 
initial value of $N=43.$ 
For $N=43,$ the triangulation $(2,3,0,43)$ 
can be constructed by using the same method as described 
in case $3$ above,
with the substitution of $m=1$.   

\vspace{0.5cm}
From the six cases above, we have considered all possible values of $N$
where $N \equiv 1(6),$ except the values $55(240)$, $79(240)$,
$85(240)$, $133(240)$, $175(240)$, $181(240)$, $199(240)$, $229(240)$.   
To determine whether triangulations exist for these other values,
we will consider triangulations of the form
 $(2,3,0,X+14m)$ and $(2,3,0,Y+28m)$ for some values of $X$ and $Y.$
The least common multiple of $14,28,240$ is $1680.$ 
Therefore the values of $N$ in this range are
in the set 
$N'=\{a+240m\mid a+240m\le 1680;\,a=55,79,85,133,175,181,199,229;\,
 m\ge 0\}$.  

\subsubsection{Case i}
There exists a triangulation $T=(2,3,0,99+14m)$ for all $m \geq 0.$
This can be constructed by glueing of the patches $(1,1,1,29+3m)_{21+2m}$ and 
 $(1,1,1,68+11m)_{21+2m}$ by 
the method $C,$ as explained in Section \ref{glupat}.
The patch $(1,1,1,29+3m)_{21+2m}$ can be constructed 
from the patch $(1,1,1,5+3h)_{5+2h}$ with the substitution of $h=8+m.$
The patch $(1,1,1,5+3h)_{5+2h}$ is of type $B$ with $k=5$
and $(1,1,1,68+11m)_{21+2m}$ is of type $A$ with $k=11$ 
(types $A$ and $B$ are as given in Section \ref{con111}).  
From the triangulation $T=(2,3,0,99+14m)$ for all $m \geq 0,$
we can get the following triangulations:

\begin{enumerate}
\item If $m=14+120k$ for all $k \geq 0$, then
we can get the triangulation $(2,3,0,295+1680k)$ for all $k \geq 0.$
This gives some values of $N$ in the $55(240)$ case.

\item If $m=23+120k$ for all $k \geq 0$, then
we can get the triangulation $(2,3,0,421+1680k)$ for all $k \geq 0.$
This gives some values of $N$ in the $181(240)$ case.

\item If $m=71+120k$ for all $k \geq 0$, then
we can get the triangulation $(2,3,0,1093+1680k)$ for all $k \geq 0.$
This gives some values of $N$ in the $133(240)$ case.

\item If $m=95+120k$ for all $k \geq 0$, then
we can get the triangulation $(2,3,0,1429+1680k)$ for all $k \geq 0.$
This gives some values of $N$ in the $229(240)$ caes.

\item If $m=119+120k$ for all $k \geq 0$, then
we can get the triangulation $(2,3,0,1765+1680k)$ for all $k \geq 0,$
which is the same as the triangulation $(2,3,0,85+1680k)$ for all $k \geq 1.$
This gives some values of $N$ in the $85(240)$ case.

\item If $m=50+120k$ for all $k \geq 0$, then
we can get the triangulation $(2,3,0,799+1680k)$ for all $k \geq 0.$
This gives some values of $N$ in the $79(240)$ case.

\item If $m=74+120k$ for all $k \geq 0$, then
we can get the triangulation $(2,3,0,1135+1680k)$ for all $k \geq 0.$
This gives some values of $N$ in the $175(240)$ case.

\item If $m=110+120k$ for all $k \geq 0$, then
we can get the triangulation $(2,3,0,1639+1680k)$ for all $k \geq 0.$
This gives some values of $N$ in the $199(240)$ case.
\end{enumerate}

\subsubsection{Case ii}
There exists a triangulation $T=(2,3,0,257+28m)$ for all $m \geq 0.$
This can be constructed by glueing of the patches $(1,1,1,161+15m)_{25+2m}$ and  $(1,1,1,94+13m)_{25+2m}$ by the method $C,$ 
as explained in Section \ref{glupat}.
The patch  $(1,1,1,161+15m)_{25+2m}$ can be constructed 
from the patch $(1,1,1,71+15h)_{13+2h}$ with the substitution of $h=6+m.$
The patch $(1,1,1,71+15h)_{13+2h}$ is of type $C$ with $k=13$ 
and $(1,1,1,94+13m)_{25+2m}$ is of type $A$ with $k=13$
(types $A$ and $C$ are as given in Section \ref{con111}).  
From the triangulation $T=(2,3,0,257+28m)$ for all $m \geq 0,$
we can get the following triangulations:

\begin{enumerate}
\item If $m=11+60k$ for all $k \geq 0$, then
we can get the triangulation $(2,3,0,565+1680k)$ for all $k \geq 0.$
This gives some values of $N$ in the $85(240)$ case.

\item If $m=23+60k$ for all $k \geq 0$, then
we can get the triangulation $(2,3,0,901+1680k)$ for all $k \geq 0.$
This gives some values of $N$ in the $181(240)$ case.

\item If $m=47+60k$ for all $k \geq 0$, then
we can get the triangulation $(2,3,0,1573+1680k)$ for all $k \geq 0.$
This gives some values of $N$ in the $133(240)$ case.

\item If $m=59+60k$ for all $k \geq 0$, then
we can get the triangulation $(2,3,0,1909+1680k)$ for all $k \geq 0,$
which is the same as the triangulation $(2,3,0,229+1680k)$
for all $k \geq 1.$
This gives some values of $N$ in the $229(240)$ case.
\end{enumerate}

\subsubsection{Case iii}
There exists a triangulation $T=(2,3,0,143+28m)$ for all values of $m \geq 0.$
This can be constructed by glueing of the patches $(1,1,1,57+11m)_{13+2m}$ and $(1,1,1,71+15m)_{13+2m}$ by the method $A,$ as explained in Section \ref{glupat}.
The patch $(1,1,1,57+11m)_{13+2m}$ can be constructed 
from the patch $(1,1,1,24+11h)_{7+2h}$ with the substitution of $h=3+m.$
The patch $(1,1,1,24+11h)_{7+2h}$ is of type $D$ with $k=7$
and $(1,1,1,71+15m)_{13+2m}$ is of type $C$ with $k=13$
(types $C$ and $D$ are as give in Section \ref{con111}).  
From the triangulation $T=(2,3,0,143+28m)$ for all $m \geq 0,$
we can get the following triangulations:

\begin{enumerate}
\item If $m=2+60k$ for all $k \geq 0$ , then
we can get the triangulation $(2,3,0,199+1680k)$ for all $k \geq 0.$
This gives some values of $N$ in the $199(240)$ case.

\item If $m=14+60k$ for all $k \geq 0$ , then
we can get the triangulation $(2,3,0,535+1680k)$ for all $k \geq 0.$
This gives some values of $N$ in the $55(240)$ case.

\item If $m=32+60k$ for all $k \geq 0$ , then
we can get the triangulation $(2,3,0,1039+1680k)$ for all $k \geq 0.$
This gives some values of $N$ in the $79(240)$ case.

\item If $m=44+60k$ for all $k \geq 0$ , then
we can get the triangulation $(2,3,0,1375+1680k)$ for all $k \geq 0.$
This gives some values of $N$ in the $175(240)$ case.
\end{enumerate}

\subsubsection{Case iv}
There exists a triangulation $T=(2,3,0,185+28m)$ for all values of $m \geq 0.$
This can be constructed by glueing of the patches $(1,1,1,53+7m)_{17+2m}$ and $(1,1,1,113+19m)_{17+2m}$ by the method $A,$ as explained in Section \ref{glupat}.
The patch $(1,1,1,53+7m)_{17+2m}$ can be constructed 
from the patch $(1,1,1,11+7h)_{5+2h}$ with the substitution of $h=6+m.$
The patch $(1,1,1,11+7h)_{5+2h}$ is of type $D$ with $k=5$
and $(1,1,1,113+19m)_{17+2m}$ is of type $E$ with $k=17$
(types $D$ and $E$ are as give in Section \ref{con111}).  
From the triangulation $T=(2,3,0,185+28m)$ for all $m \geq 0,$
we can get the following triangulations:

\begin{enumerate}

\item If $m=5+60k$ for all $k \geq 0$ , then
we can get the triangulation $(2,3,0,325+1680k)$ for all $k \geq 0.$
This gives some values of $N$ in the $85(240)$ case.

\item If $m=17+60k$ for all $k \geq 0$ , then
we can get the triangulation $(2,3,0,661+1680k)$ for all $k \geq 0.$
This gives some values of $N$ in the $181(240)$ case.

\item If $m=41+60k$ for all $k \geq 0$ , then
we can get the triangulation $(2,3,0,1333+1680k)$ for all $k \geq 0.$
This gives some values of $N$ in the $133(240)$ case.

\item If $m=53+60k$ for all $k \geq 0$ , then
we can get the triangulation $(2,3,0,1669+1680k)$ for all $k \geq 0.$
This gives some values of $N$ in the $229(240)$ case.

\end{enumerate}

\subsubsection{Case v}
There exists a triangulation $T=(2,3,0,167+28m)$ for all $m \geq 0.$
This can be constructed by glueing of the patches $(1,1,1,46+9m)_{17+2m}$ and 
$(1,1,1,119+19m)_{17+2m}$ by the method $C,$ 
as explained in Section \ref{glupat}.
The patch $(1,1,1$,$119$ $+19m)_{17+2m}$ can be constructed 
from the patch $(1,1,1,100+19h)_{15+2h}$ with the substitution of $h=1+m.$
The patch $(1,1,1,100+19h)_{15+2h}$ is of type $C$ with $k=15$
and the patch  $(1,1,1,46+9m)_{17+2m}$ is of type $A$ with $k=9$
(types $A$ and $C$ are as given in Section \ref{con111}).  
From the triangulation $T=(2,3,0,167+28m)$ for all $m \geq 0,$
we can get the following triangulations:

\begin{enumerate}
\item If $m=14+60k$ for all $k \geq 0$, then
we can get the triangulation $(2,3,0,559+1680k)$ for all $k \geq 0.$
This gives some values of $N$ in the $79(240)$ case.

\item If $m=26+60k$ for all $k \geq 0$, then
we can get the triangulation $(2,3,0,895+1680k)$ for all $k \geq 0.$
This gives some values of $N$ in the $175(240)$ case.

\item If $m=44+60k$ for all $k \geq 0$, then
we can get the triangulation $(2,3,0,1399+1680k)$ for all $k \geq 0.$
This gives some values of $N$ in the $199(240)$ case.

\item If $m=56+60k$ for all $k \geq 0$, then
we can get the triangulation $(2,3,0,1735+1680k)$ for all $k \geq 0,$
which is the same as the triangulation $(2,3,0,55+1680k)$
for all $k \geq 1.$
This gives some values of $N$ in the $55(240)$ case.
\end{enumerate}

\subsubsection{Case vi}
There exists a triangulation $T=(2,3,0,263+28m)$ for all $m \geq 0.$
This can be constructed by glueing of the patches $(1,1,1,35+3m)_{25+2m}$ and 
$(1,1,1,201+23m)_{25+2m}$ by the method $A,$ 
as explained in Section \ref{glupat}.
The patch  $(1,1,1,35+3m)_{25+2m}$ can be constructed 
from the patch $(1,1,1,2+3h)_{3+2h}$ with the substitution of $h=11+m.$
The patch $(1,1,1,2+3h)_{3+2h}$ is of type $D$ with $k=3$ 
and $(1,1,1,201+23m)_{25+2m}$ is of type $F$ with $k=25$ 
(types $D$ and $F$ are as given in Section \ref{con111}). 
From the triangulation $T=(2,3,0,263+28m)$ for all $m \geq 0,$
we can get the following triangulations:

\begin{enumerate}
\item If $m=2+60k$ for all $k \geq 0$, then
we can get the triangulation $(2,3,0,319+1680k)$ for all $k \geq 0.$
This gives some values of $N$ in the $79(240)$ case.

\item If $m=14+60k$ for all $k \geq 0$, then
we can get the triangulation $(2,3,0,655+1680k)$ for all $k \geq 0.$
This gives some values of $N$ in the $175(240)$ case.

\item If $m=32+60k$ for all $k \geq 0$, then
we can get the triangulation $(2,3,0,1159+1680k)$ for all $k \geq 0.$
This gives some values of $N$ in the $199(240)$ case.

\item If $m=44+60k$ for all $k \geq 0$, then
we can get the triangulation $(2,3,0,1495+1680k)$ for all $k \geq 0.$
This gives some values of $N$ in the $55(240)$ case.
\end{enumerate}

\subsubsection{Case vii}
There exists a triangulation $T=(2,3,0,147+28m)$ for all $m \geq 0.$
This can be constructed by glueing of the patches $(1,1,1,20+3m)_{15+2m}$ and 
$(1,1,1,110+23m)_{15+2m}$ 
by the method $A,$ as explained in Section \ref{glupat}.
The patch  $(1,1,1,20+3m)_{15+2m}$ can be constructed 
from the patch $(1,1,1,5+3h)_{5+2h}$ with the substitution of $h=5+m.$
The patches $(1,1,1,5+3h)_{5+2h}$ and $(1,1,1,110+23m)_{15+2m}$ are 
of type $B$ with $k=5$ and $k=15$ respectively (type $B$ is as given in Section \ref{con111}). 
From the triangulation $T=(2,3,0,147+28m)$ for all $m \geq 0,$
we can get the following triangulations:

\begin{enumerate}
\item If $m=1+60k$ for all $k \geq 0$, then
we can get the triangulation  $(2,3,0,175+1680k)$ for all $k \geq 0.$
This gives some values of $N$ in the $175(240)$ case.

\item If $m=19+60k$ for all $k \geq 0$, then
we can get the triangulation $(2,3,0,679+1680k)$ for all $k \geq 0.$
This gives some values of $N$ in the $199(240)$ case.

\item If $m=31+60k$ for all $k \geq 0$, then
we can get the triangulation  $(2,3,0,1015+1680k)$ for all $k \geq 0.$
This gives some values of $N$ in the $55(240)$ case.

\item If $m=49+60k$ for all $k \geq 0$, then
we can get the triangulation $(2,3,0,1519+1680k)$ for all $k \geq 0.$
This gives some values of $N$ in the $79(240)$ case.
\end{enumerate}

\subsubsection{Case viii}
There exists a triangulation $T=(2,3,0,135+28m)$ for all $m \geq 0.$
This can be constructed by glueing of the patches $(1,1,1,39+7m)_{13+2m}$
and $(1,1,1,81+19m)_{13+2m}$ 
by the method $A,$ as explained in Section \ref{glupat}.
The patch  $(1,1,1,39+7m)_{13+2m}$ can be constructed 
from the patch $(1,1,1,11+7h)_{5+2h}$ with the substitution of $h=4+m.$
The patch $(1,1,1,11+7h)_{5+2h}$ is of type $D$ with $k=5$ 
and $(1,1,1,81+19m)_{13+2m}$ is of type $C$ with $k=13$
(types $C$ and $D$ are as given in Section \ref{con111}). 
From the triangulation $T=(2,3,0,135+28m)$ for all $m \geq 0,$
we can get the following triangulations:

\begin{enumerate}
\item If $m=10+60k$ for all $k \geq 0$, then
we can get the triangulation $(2,3,0,415+1680k)$ for all $k \geq 0.$
This gives some values of $N$ in the $175(240)$ case.

\item If $m=28+60k$ for all $k \geq 0$, then
we can get the triangulation $(2,3,0,919+1680k)$ for all $k \geq 0.$
This gives some values of $N$ in the $199(240)$ case.

\item If $m=40+60k$ for all $k \geq 0$, then
we can get the triangulation $(2,3,0,1255+1680k)$ for all $k \geq 0.$
This gives some values of $N$ in the $55(240)$ case.

\item If $m=58+60k$ for all $k \geq 0$, then
we can get the triangulation 
$(2,3,0,1759+1680k)$ for all $k \geq 0,$
which is the same as the triangulation 
$(2,3,0,79+1680k)$ for all $k \geq 1.$
This gives some values of $N$ in the $79(240)$ case.
\end{enumerate}

\subsubsection{Case ix}
There exists a triangulation $T=(2,3,0,159+28m)$ for all $m \geq 0.$
This can be constructed by glueing of the patches $(1,1,1,62+11m)_{15+2m}$ 
and 
$(1,1,1,80+15m)_{15+2m}$ by the method $A,$ 
as explained in Section \ref{glupat}.
The patch $(1,1,1,62+11m)_{15+2m}$ can be constructed 
from the patch $(1,1,1,51+11h)_{13+2h}$ with the substitution of $h=1+m$.  
The patches $(1,1,1,51+11h)_{13+2h}$ and
$(1,1,1,80+15m)_{15+2m}$ are of type $E$ with $k=13$ and
$k=15$ respectively      
(type $E$ is as given in Section \ref{con111}). 
From the triangulation $T=(2,3,0,159+28m)$ for all $m \geq 0,$
we can get the following triangulations:

\begin{enumerate}
\item If $m=40+60k$ for all $k \geq 0$, then
we can get the triangulation  $(2,3,0,1279+1680k)$ for all $k \geq 0.$
This gives some values of $N$ in the $79(240)$ case.

\item If $m=52+60k$ for all $k \geq 0$, then
we can get the triangulation  $(2,3,0,1615+1680k)$ for all $k \geq 0.$
This gives some values of $N$ in the $175(240)$ case.

\item If $m=10+60k$ for all $k \geq 0$, then
we can get the triangulation  $(2,3,0,439+1680k)$ for all $k \geq 0.$
This gives some values of $N$ in the $199(240)$ case.

\item If $m=22+60k$ for all $k \geq 0$, then
we can get the triangulation  $(2,3,0,775+1680k)$ for all $k \geq 0.$
This gives some values of $N$ in the $55(240)$ case.
\end{enumerate}

\subsubsection{Case x}
There exists a triangulation $T=(2,3,0,377+28m)$ for all $m \geq 0.$
This can be constructed by glueing of the patches $(1,1,1,238+21m)_{41+2m}$ 
and $(1,1,1,137+7m)_{41+2m}$ by the method $C,$ 
as explained in Section \ref{glupat}.
The patch $(1,1,1,137+7m)_{41+2m}$ can be constructed 
from the patch $(1,1,1,18+7h)_{7+2h}$ with the substitution of $h=17+m$.   
The patch $(1,1,1,18+7h)_{7+2h}$ is of type $B$ with $k=7$  
and $(1,1,1,238+21m)_{41+2m}$ is of type $A$ with $k=21$
(types $A$ and $B$ are as given in Section \ref{con111}).
From the triangulation $T=(2,3,0,377+28m)$ for all $m \geq 0,$
we can get the following triangulations:

\begin{enumerate}
\item If $m=17+60k$ for all $k \geq 0$, then
we can get the triangulation  $(2,3,0,853+1680k)$ for all $k \geq 0.$
This gives some values of $N$ in the $133(240)$ case.

\item If $m=29+60k$ for all $k \geq 0$, then
we can get the triangulation  $(2,3,0,1189+1680k)$ for all $k \geq 0.$
This gives some values of $N$ in the $229(240)$ case. 

\item If $m=41+60k$ for all $k \geq 0$, then
we can get the triangulation  $(2,3,0,1525+1680k)$ for all $k \geq 0.$
This gives some values of $N$ in the $85(240)$ case. 

\item If $m=53+60k$ for all $k \geq 0$, then
we can get the triangulation  $(2,3,0,1861+1680k)$ for all $k \geq 0,$
which is the same as the triangulation 
$(2,3,0,181+1680k)$ for all $k \geq 1.$
This gives some values of $N$ in the $181(240)$ case. 
\end{enumerate}

\subsubsection{Case xi}
There exists a triangulation $T=(2,3,0,137+28m)$ for all $m \geq 0.$
This can be constructed by glueing of the patches $(1,1,1,51+11m)_{13+2m}$ 
and $(1,1,1,71+15m)_{13+2m}$ by the method $A,$ 
as explained in Section \ref{glupat}.
The patch $(1,1,1,51+11m)_{13+2m}$ is of type $E$ with $k=13$
and $(1,1,1,71+15m)_{13+2m}$ is of type $C$ with $k=13$  
(types $C$ and $E$ are as given in Section \ref{con111}).
From the triangulation $T=(2,3,0,137+28m)$ for all $m \geq 0,$
we can get the following triangulations:

\begin{enumerate}
\item If $m=17+60k$ for all $k \geq 0$, then
we can get the triangulation  $(2,3,0,613+1680k)$ for all $k \geq 0.$
This gives some values of $N$ in the $133(240)$ case.

\item If $m=29+60k$ for all $k \geq 0$, then
we can get the triangulation  $(2,3,0,949+1680k)$ for all $k \geq 0.$
This gives some values of $N$ in the $229(240)$ case.

\item If $m=41+60k$ for all $k \geq 0$, then
we can get the triangulation  $(2,3,0,1285+1680k)$ for all $k \geq 0.$
This gives some values of $N$ in the $85(240)$ case.

\item If $m=53+60k$ for all $k \geq 0$, then
we can get the triangulation  
$(2,3,0,1621+1680k)$ for all $k \geq 0.$
This gives some values of $N$ in the $181(240)$ case.
\end{enumerate}

\subsubsection{Case xii}
There exists a triangulation $T=(2,3,0,161+28m)$ for all $m \geq 0.$
This can be constructed by glueing of the patches $(1,1,1,46+9m)_{17+2m}$ 
and $(1,1,1,113+19m)_{17+2m}$ by the method $C,$ 
as explained in Section \ref{glupat}.
The patch $(1,1,1,46+9m)_{17+2m}$ is of type $A$ with $k=9$
and $(1,1,1,113+19m)_{17+2m}$ is of type $E$ with $k=17$
(types $A$ and $E$ are as given in Section \ref{con111}).
From the triangulation $T=(2,3,0,161+28m)$ for all $m \geq 0,$
we can get the following triangulations:

\begin{enumerate}
\item If $m=59+60k$ for all $k \geq 0$, then
we can get the triangulation $(2,3,0,1813+1680k)$ for all $k \geq 0,$
which is the same as the triangulation 
$(2,3,0,133+1680k)$ for all $k \geq 1$.   
This gives some values of $N$ in the $133(240)$ case.

\item If $m=23+60k$ for all $k \geq 0$, then
we can get the triangulation $(2,3,0,805+1680k)$ for all $k \geq 0.$
This gives some values of $N$ in the $85(240)$ case.

\item If $m=11+60k$ for all $k \geq 0$, then
we can get the triangulation $(2,3,0,469+1680k)$ for all $k \geq 0.$
This gives some values of $N$ in the $229(240)$ case.

\item If $m=35+60k$ for all $k \geq 0$, then
we can get the triangulation $(2,3,0,1141+1680k)$ for all $k \geq 0.$
This gives some values of $N$ in the $181(240)$ case.
\end{enumerate}

\subsubsection{Case xiii}
There exists a triangulation $T=(2,3,0,233+28m)$ for all $m \geq 0.$
This can be constructed by glueing of the patches $(1,1,1,94+13m)_{25+2m}$ 
and $(1,1,1,137+15m)_{25+2m}$ by the method $C,$ as explained in Section \ref{glupat}.
The patch $(1,1,1,137+15m)_{25+2m}$ can be constructed 
from the patch $(1,1,1,107+15h)_{21+2h}$ with the substitution of $h=2+m.$
The patch $(1,1,1,107+15h)_{21+2h}$ is of type $F$ with $k=21$
and $(1,1,1,94+13m)_{25+2m}$ is of type $A$ with $k=13$
(types $A$ and $F$ are as given in Section \ref{con111}). 
From the triangulation $T=(2,3,0,233+28m)$ for all $m \geq 0,$
we can get the following triangulations:

\begin{enumerate}
\item If $m=29+60k$ for all $k \geq 0$, then
we can get the triangulation $(2,3,0,1045+1680k)$ for all $k \geq 0.$
This gives some values of $N$ in the $85(240)$ case.

\item If $m=5+60k$ for all $k \geq 0$, then
we can get the triangulation $(2,3,0,373+1680k)$ for all $k \geq 0.$
This gives some values of $N$ in the $133(240)$ case.

\item If $m=41+60k$ for all $k \geq 0$, then
we can get the triangulation $(2,3,0,1381+1680k)$ for all $k \geq 0.$
This gives some values of $N$ in the $181(240)$ case.

\item If $m=17+60k$ for all $k \geq 0$, then
we can get the triangulation $(2,3,0,709+1680k)$ for all $k \geq 0.$
This gives some values of $N$ in the $229(240)$ case.
\end{enumerate}

\vspace{0.5cm}
From the above $13$ cases we found triangulations for all values of  
$N$ in $N'$  
except for the values $55,79,85,181,229$.   
We can construct the triangulations for these values as follows:

\begin{enumerate}
\item The triangulation $(2,3,0,55)$ exists. This can be
constructed by glueing of the patches $(1,1,1,19)_9$ and 
$(1,1,1,25)_9$ by the method $A,$ as explained in Section \ref{glupat}.
The patch $(1,1,1,19)_9$ is of type $M$ with $k=9$
and $(1,1,1,25)_9$ is of type $L$ with $k=9$
(types $M$ and $L$ are as given in Section \ref{con111}). 

\item The triangulation $(2,3,0,79)$ exists. This can be
constructed by glueing of the patches $(1,1,1,26)_{11}$ and 
$(1,1,1,40)_{11}$ by the method $A,$ as explained in Section \ref{glupat}.
The patch $(1,1,1,26)_{11}$ is of type $O$ with $k=11$
and $(1,1,1,40)_{11}$ is of type $M$ with $k=11$
(types $O$ and $M$ are as given in Section \ref{con111}). 

\item The triangulation $(2,3,0,85)$ exists. This can be
constructed by glueing of the patch $(1,1,1,37)_9$ with
itself (that is, $(1,1,1,37)_9$ and $(1,1,1,37)_9$)   
by the method $A,$ as explained in Section \ref{glupat}.
The patch $(1,1,1,37)_9$ is of type $I$ with $k=9$
(type $I$ is as given in Section \ref{con111}). 

\item The triangulation $(2,3,0,181)$ exists. This can be
constructed by glueing of the patch $(1,1,1,83)_{13}$ with itself   
(that is, $(1,1,1,83)_{13}$ and $(1,1,1$, $83)_{13}$) by the method $A,$
as explained in Section \ref{glupat}.
The patch $(1,1,1$, $83)_{13}$ is of type $I$ with $k=13$
(type $I$ is as given in Section \ref{con111}). 

\item The triangulation $(2,3,0,229)$ exists. This can be
constructed by glueing of the patches $(1,1,1,31)_{23}$ and 
$(1,1,1,196)_{23}$ by the method $B$, as explained in Section \ref{glupat}.
The patch $(1,1,1,31)_{23}$ is of type $A$ with $k=3,m=9$
and $(1,1,1,196)_{23}$ is of type $E$ with $k=19,m=2$
(types $A$ and $E$ are as given in Section \ref{con111}). 
\end{enumerate}

\subsection{$3(6)$ type} 
To construct this type, we will consider the triangulations of the form   
$(2,3,0,$ $X+Ym)$ with $Y=6,12,18,24,30,36$ and all possible values of $X.$ 
The least common multiple of $6,12,18,24,30,36$ is $360$. 
Therefore the values
of $N$ to consider in this $3(6)$ type are  
$\{ 3+Ym\mid 3+Ym\le 360,\,Y=6,12,18,24,30,36;\,m\ge 0\}$. 
For the constructions we will consider the following cases:

\subsubsection{Case 1}
The triangulation $(2,3,0,21+12m)$ exists for all $m \geq  0,$
by glueing of the patches $(0,3,0,16)_{12}$ and $(2,0,0,5)_{12}.$
The patch $(0,3,0,16)_{12}$ can be constructed from the patch 
given in Theorem \ref{T1}(iii)
with
$h=4,k=2$ and the patch $(2,0,0,5)_{12}$ 
is a special case of patches of the type $(2,0,0)$
with $k=6,r=0$ (as discussed in Section \ref{patch200}).

\subsubsection{Case 2}
The triangulation $(2,3,0,33+18m)$ exists for all $m \geq  0,$
by glueing of the patches $(0,3,0,25)_{18}$ and $(2,0,0,8)_{18}.$
The patch $(0,3,0,25)_{18}$ can be constructed from the patch  
given in Theorem \ref{T1}(ii) with   
$h=6$ and the patch $(2,0,0,8)_{18}$ 
is a special case of patches of the type $(2,0,0)$
with $k=9,r=0$ (as discussed in Section \ref{patch200}).

\subsubsection{Case 3}
The triangulation $(2,3,0,57+18m)$ exists for all $m \geq  0,$
by glueing of the patches $(0,3,0,49)_{18}$ and $(2,0,0,8)_{18}.$
The patch $(0,3,0,49)_{18}$ can be constructed from the patch  
given in Theorem \ref{T1}(iv) with 
$h=6,k=4,l=2$ and the patch $(2,0,0,8)_{18}$ 
is a special case of patches of the type $(2,0,0)$
with $k=9,r=0$ (as discussed in Section \ref{patch200}).

\subsubsection{Case 4}
The triangulation $(2,3,0,93+30m)$ exists for all $m \geq  0,$
by glueing of the patches $(0,3,0,79)_{30}$ and $(2,0,0,14)_{30}.$
The patch $(0,3,0,79)_{30}$ can be constructed from the patch  
given in Theorem \ref{T1}(iii) with 
$h=10,k=2$ and the patch $(2,0,0,14)_{30}$ 
is a special case of patches of the type $(2,0,0)$
with $k=15,r=0$ (as discussed in Section \ref{patch200}).

\subsubsection{Case 5}
The triangulation $(2,3,0,141+30m)$ exists for all $m \geq  0,$
by glueing of the patches $(0,3,0,127)_{30}$ and $(2,0,0,14)_{30}.$
The patch $(0,3,0,127)_{30}$ can be constructed from the patch  
given in Theorem \ref{T1}(iv) with 
$h=10$, $k=4$, $l=4$ and the patch $(2,0,0,14)_{30}$ 
is a special case of patches of the type $(2,0,0)$
with $k=15,r=0$ (as discussed in Section \ref{patch200}).

\subsubsection{Case 6}
The triangulation $(2,3,0,129+30m)$ exists for all $m \geq  0,$
by glueing of the patches $(0,3,0,115)_{30}$ and $(2,0,0,14)_{30}.$
The patch $(0,3,0,115)_{30}$ can be constructed from the patch  
given in Theorem \ref{T1}(iv) with   
$h=10$, $k=6$, $l=2$ and the patch $(2,0,0,14)_{30}$ 
is a special case of patches of the type $(2,0,0)$
with $k=15,r=0$ (as discussed in Section \ref{patch200}).

\subsubsection{Case 7}
The triangulation $(2,3,0,165+30m)$ exists for all $m \geq  0,$
by glueing of the patches $(0,3,0,151)_{30}$ and $(2,0,0,14)_{30}.$
The patch $(0,3,0,151)_{30}$ can be constructed from the patch  
given in Theorem \ref{T1}(iv) with    
$h=10$, $k=8$, $l=6$ and the patch $(2,0,0,14)_{30}$ 
is a special case of patches of the type $(2,0,0)$
with $k=15,r=0$ (as discussed in Section \ref{patch200}).

\vspace{0.5cm}
From the above seven cases, we found triangulations of the type 
$(2,3,0,N)$ for all $N\equiv 3(6)$ 
except for the case $N\equiv 27(180)$ and the values 
$N=3,9,15,39,63,99,135$.  
We next find triangulations for the case $N\equiv 27(180)$.    
To do this, we will consider the following cases:  
\medskip

\subsubsection {Case (i)}
The triangulation $T=(2,3,0,27+1260m)$ exists for all $m \geq 1.$
This can be constructed as follows:\\
Glueing of the patches $(0,3,0,301)_{42}$ and $(2,0,0,20)_{42}$
gives a triangulation $T_1=(2,3,0,321+42h)$ 
for $h \geq 0$.   
The patch $(2,0,0,20)_{42}$ is a special case of patches of the type $(2,0,0)$
with $k=21,r=0$ (as discussed in Section \ref{patch200})
and the patch $(0,3,0,301)_{42}$ is
as given in Theorem \ref{T1}(iv) 
with $h=14,k=12,l=10.$ 
If we substitute $h=23+30m$ in the triangulation $T_1=(2,3,0,321+42h),$  
then we can get the triangulation $(2,3,0,1287+1260m)$ for   
$m \geq 0$ which is the same as the triangulation $T=(2,3,0,27+1260m)$ 
for $m \geq 1$.   
The triangulation $T$ gives all the values of the $27(1260)$
case, except the initial value $27.$
For $N=27$ the triangulation $(2,3,0,27)$ can be constructed as follows:\\
Glueing of the patches $(1,1,1,4+3m)_{5+2m}$ 
and $(1,1,1,11+7m)_{5+2m}$ by the method $C,$ 
as explained in Section \ref{glupat},
gives a triangulation $(2,3,0,17+10m)$ for $m \geq 0$.  
If $m=1,$ then we can get the triangulation $(2,3,0,27).$
The patch $(1,1,1,4+3m)_{5+2m}$ is of type $A$ with $k=3$
and $(1,1,1,11+7m)_{5+2m}$ is of type $D$ with $k=5$
(types $A$ and $D$ are as given in Section \ref{con111}).

\subsubsection {Case (ii)}
The triangulation $T=(2,3,0,207+1260m)$ exists for all $m \geq 1.$
This can be constructed as follows:\\
Glueing of the patches $(0,3,0,229)_{42}$ and $(2,0,0,20)_{42}$
gives a triangulation $T_1=(2,3,0,249+42h)$ for   
$h \geq 0.$
The patch $(2,0,0,20)_{42}$ is a special case of patches of the type $(2,0,0)$
with $k=21,r=0$ (as discussed in Section \ref{patch200})
and the patch $(0,3,0,229)_{42}$ is
as given in Theorem \ref{T1}(iv) 
with $h=14,k=6,l=4.$ 
If we substitute $h=29+30m$ in the triangulation 
$T_1=(2,3,0,249+42h),$ then we can get the triangulation 
$(2,3,0,1467+1260m)$ for $m \geq 0$,  
which is the same as the triangulation $T=(2,3,0,207+1260m)$ 
for $m \geq 1$.   
The triangulation $T$ gives all the values of the $207(1260)$
case, except the initial value $207.$
For $N=207$ the triangulation $(2,3,0,207)$ can be constructed as follows:\\
Glueing of the patches $(1,1,1,4+3m)_{5+2m}$ 
and $(1,1,1,11+7m)_{5+2m}$ by the method $C,$ as explained 
in Section \ref{glupat},
gives a triangulation $(2,3,0,17+10m)$ for $m \geq 0$.  
If $m=19,$ then we can get the triangulation $(2,3,0,207).$
The patch $(1,1,1,4+3m)_{5+2m}$ is of type $A$ with $k=3$
and $(1,1,1,11+7m)_{5+2m}$ is of type $D$ with $k=5$
(types $A$ and $D$ are as given in Section \ref{con111}).

\subsubsection {Case (iii)}
The triangulation $T=(2,3,0,387+1260m)$ exists for all $m \geq 0.$
This can be constructed as follows:\\
Glueing of the patches $(0,3,0,157)_{42}$ and $(2,0,0,20)_{42}$
gives a triangulation $T_1=(2,3,0,177+42h)$ for    
$h \geq 0.$
The patch $(2,0,0,20)_{42}$ is a special case of patches of the type $(2,0,0)$
with $k=21,r=0$ (as discussed in Section \ref{patch200})
and the patch $(0,3,0,157)_{42}$ is
as given in Theorem \ref{T1}(iii) 
with $h=14,k=4$.    
If we substitute $h=5+30m$ in the triangulation $T_1=(2,3,0,177+42h),$
then we can get the triangulation $T=(2,3,0,387+1260m)$ for $m \geq 0.$ 
The triangulation $T$ gives all the values of the $387(1260)$ case. 

\subsubsection {Case (iv)}
The triangulation $T=(2,3,0,567+1260m)$ exists for all $m \geq 0.$
This can be constructed as follows:\\
Glueing of the patches $(0,3,0,169)_{42}$ and $(2,0,0,20)_{42}$
gives a triangulation $T_1=(2,3,0,189+42h)$ for $h\geq 0$.   
The patch $(2,0,0,20)_{42}$ is a special case of patches of the type $(2,0,0)$
with $k=21,r=0$ (as discussed in Section \ref{patch200})
and the patch $(0,3,0,169)_{42}$ is
as given in Theorem \ref{T1}(iv)   
with $h=14,k=2,l=2$.   
If we substitute $h=9+30m$ in the triangulation $T_1=(2,3,0,189+42h),$ 
then we can get the triangulation $T=(2,3,0,567+1260m)$ for $m\geq 0.$ 
The triangulation $T$ gives all the values of the $567(1260)$ case.

\subsubsection {Case (v)}
The triangulation $T=(2,3,0,747+1260m)$ exists for $m\geq 0.$
This can be constructed as follows:\\
Glueing of the patches $(0,3,0,265)_{42}$ and $(2,0,0,20)_{42}$
gives a triangulation  $T_1=(2,3,0,285+42h)$ 
for $h\geq 0.$
The patch $(2,0,0,20)_{42}$
is a special case of patches of the type $(2,0,0)$
with $k=21,r=0$ (as discussed in Section \ref{patch200})
and the patch $(0,3,0,265)_{42}$ is
from the patch given in Theorem \ref{T1}(iv) 
with $h=14,k=10,l=6$.    
If we substitute $h=11+30m$ in the triangulation $T_1=(2,3,0,285+42h),$  
then we can get the triangulation $T=(2,3,0,747+1260m)$ for $m\geq 0$.  
The triangulation $T$ gives all the values of the $747(1260)$ case.

\subsubsection {Case (vi)}
The triangulation $T=(2,3,0,927+1260m)$ exists for all $m\geq 0.$
This can be constructed as follows:\\
Glueing of the patches $(0,3,0,277)_{42}$ and $(2,0,0,20)_{42}$
gives a triangulation  $T_1=(2,3,0,297+42h)$ for 
$h\geq 0$. The patch $(2,0,0,20)_{42}$ 
is a special case of patches of the type $(2,0,0)$
with $k=21,r=0$ (as discussed in Section \ref{patch200})
and the patch $(0,3,0,277)_{42}$ is
from the patch given in Theorem \ref{T1}(iv) 
with $h=14,k=8,l=8$.   
If we substitute $h=15+30m$ in the triangulation $T_1=(2,3,0,297+42h),$  
then we can get the triangulation $T=(2,3,0,927+1260m)$ for $m\geq 0.$
The triangulation $T$ gives all the values of the $927(1260)$ case.

\subsubsection {Case (vii)}
The triangulation $T=(2,3,0,1107+1260m)$ exists for all $m \geq 0.$
This can be constructed as follows:\\
Glueing of the patches $(0,3,0,205)_{42}$ and $(2,0,0,20)_{42}$
gives a triangulation $T_1=(2,3,0,225+42h)$ for $h\geq 0.$
The patch $(2,0,0,20)_{42}$ 
is a special case of patches of the type $(2,0,0)$
with $k=21,r=0$ (as discussed in Section \ref{patch200})
and the patch $(0,3,0,205)_{42}$ is
as given in Theorem \ref{T1}(iv) 
with $h=14,k=8,l=2.$ 
If we substitute $h=21+30m$ in the triangulation $T_1=(2,3,0,225+42h),$  
then we can get the triangulation $T=(2,3,0,1107+1260m)$ for $m\geq 0$.   
The triangulation $T$ gives all the values of the $1107(1260)$ case.

\vspace{0.3cm}
Case (i) to (vii) above give triangulations for all the values
of the $27(180)$ case. The remaining values of $N$ to consider
are $N=3,9,15,39,63,99,135.$

\begin{enumerate}
\item
The triangulation $(2,3,0,9)$ can be constructed as follows: 
Glueing of the patches $(1,1,1,2+3m)_{3+2m}$ and $(1,1,1,2+3m)_{3+2m}$ by
the method $A,$ as explained in Section \ref{glupat},    
gives a triangulation $(2,3,0,9+8m)$ for $m \geq 0$.  
If $m=0,$ then we can get the triangulation $(2,3,0,9).$
The patch $(1,1,1,2+3m)_{3+2m}$ is of type $D$ with $k=3$ 
(type $D$ is as given in Section \ref{con111}). 

\item
The triangulation $(2,3,0,39)$ can be constructed by glueing 
of the patches $(1,1,1,12)_7$ and $(1,1,1,18)_7$ by
the method $A,$ as explained in Section \ref{glupat}.
The patch $(1,1,1,12)_7$ is of type $K$ with $k=7$ and 
$(1,1,1,18)_7$ is of type $J$ with $k=7$
(types $J$ and $K$ are as given in Section \ref{con111}). 

\item 
The triangulation $(2,3,0,63)$ can be constructed by glueing
of the patches $(1,1,1,11)_9$ and $(1,1,1,41)_9$ by the method $A,$
as explained in Section \ref{glupat}.
The patch $(1,1,1,11)_9$ is given in Table $2.6$
and $(1,1,1,41)_9$ is of type $H$ with $k=9$
(type $H$ is as given in Section \ref{con111}). 

\item
The triangulation $(2,3,0,99)$ can be constructed by glueing 
of the patches $(1,1,1,40)_{11}$ and $(1,1,1,46)_{11}$ by
the method $A,$ as explained in Section \ref{glupat}.
The patch $(1,1,1,40)_{11}$ is of type $M$ with $k=11$ and 
$(1,1,1,46)_{11}$ is of type $L$ with $k=11$
(types $K$ and $L$ are as given in Section \ref{con111}). 

\item
The triangulation $(2,3,0,135)$ can be constructed by glueing 
of the patches $(1,1,1,33)_{13}$ and $(1,1,1,87)_{13}$ by
the method $A,$ as explained in Section \ref{glupat}.
The patch $(1,1,1,33)_{13}$ is of type $E$ with $k=11,\,m=1$ and 
$(1,1,1,87)_{13}$ is of type $H$ with $k=13$
(types $E$ and $H$ are as given in Section \ref{con111}). 
\end{enumerate}

We have not been able to construct triangulations for the case
$N=3,15.$

\subsection {$5(6)$ type}
The triangulation $T=(2,3,0,5+6m)$ for all $m\geq 0$
can be constructed by glueing of the patches $(2,0,0,2)_6$ and $(0,3,0,3)_6.$ The patch $(2,0,0,2)_6$ is a special case of patches of 
the type $(2,0,0)$
with $k=3,r=0$ (as discussed in Section \ref{patch200})
and the patch $(0,3,0,3)_6$ is as given in 
Theorem \ref{T1}(ii) with $h=2$.   
The triangulation $T$ gives all the values of the $5(6)$ type.

\section{Triangulations of type $(3,0,3)$} 
In this section we show that the triangulation 
$(3,0,3,N)$ exists for all values of $N$ except for $N=0,2$
and possibly $N=4,12$.   
 
\begin{table}[h]
\begin{center}
\begin{tabular}{|c|l|} \hline 
$N=a_6 $  &  Triangulation \\ \hline
1  & $~ \{ 125,~127,~135,~136,~146,~147,~247,~235,~346,~234 ~\} $ \\ \hline 

3  & $~ \{ 123,~129,~138,~148,~145,~237,~256,~267,~259,~348,~346,~367, $ \\
   & $ ~~~159,~456 ~\} $ \\ \hline 
 
6  & $~ \{ 123,~134,~145,~15B,~12B,~237,~27C,~2AC,~2AB,~345,~356, $ \\
   & $ ~~~367,~569,~689,~678,~789,~ABC,~59B,79C,9BC ~ \}  $  \\ \hline 
    
8  & $ ~\{ 12B,~128,~15B,~15C,~189,~19C,~237,~23E,~27B,~28E,~348, $  \\
   & $ ~~~34A,~37A,~38E,~456,~459,~46A,~489,~56B,~59C,~67B,~67D, $ \\
   & $ ~~~6AD,~7AD ~ \} $ \\ \hline
   
\end{tabular}  
\end{center}
\caption{Some triangulations of type $(3,0,3)$}
\label{tab303}
\end{table}

To construct triangulations of this type we will consider two cases: 
(1) $N$ odd, ~(2)~$N$ even.

\subsection{$N$ odd } \label{odd303}
There exists a triangulation $T:(3,0,3,2m+1)$ for all $m \geq 0.$
This can be constructed as follows: \\

In Section \ref{Tri400} it will be shown that the   
triangulation $T_1:(4,0,0,2n)$ exists for each $n \geq 2$.
This triangulation has 
a point of degree $3$ adjacent to three points, say $x,y,z$ of degree $6$.  
If we remove this point of degree three, then the degrees of the points $x,y,z$
will become $5$, that is, we will get three points of degree $5$.   
Hence, we will get the triangulation $T:(3,0,3,2m+1)$ for $m\geq 0$.  

\subsection*{$N$ even }
In this case we will consider three types: $0(6),2(6),4(6).$

\subsection{$0(6)$ type}
To construct triangulations of this type, we will consider two 
sub-types: $0(12),$ $6(12).$

\subsubsection{Case 1:$~0(12)$ type}
There exists a triangulation $T:(3,0,3,0+12m)$ for all $m \geq 3.$
This can be constructed as follows:\\
Glueing of the patches $(1,0,3,31)_{12}$ and $(2,0,0,5)_{12}$ gives 
a triangulation $(3,0,3,36+12m)$ for $m \geq 0,$ 
which is the same as the triangulation $T:(3,0,3,0+12m)$ for $m \geq 3.$ 
The patch $(1,0,3,31)_{12}$ can be constructed by using 
the edge-fullering construction
on the patch $(1,0,3,7)_6,$ that is, $EF(1,0,3,7)_6= (1,0,3,31)_{12}$ where 
$(1,0,3,7)_6$ is given in Table $2.4.$
The patch $(2,0,0,5)_{12}$ is a special case of patches of the type $(2,0,0)$
with $k=6,r=0$ (as discussed in Section \ref{patch200}).   

The triangulation $T$ gives all the values of the $0(12)$ type, 
except the initial
values $N=0,12,24.$  For $N=0,$ the triangulation cannot be
constructed as has been shown by Eberhard \cite{Eber} 
and Br\"uckner \cite{bruck} 
(see also Gr\"unbaum \cite{GR}). 
We have not been able to construct a triangulation
for the case $N=12.$ 
For $N=24,$ the triangulation 
$(3,0,3,24)$ can be constructed by using the edge-fullering 
construction on the triangulation $(3,0,3,3)$, that is, 
$ EF(3,0,3,3)=(3,0,3,24).$ 
The construction of the triangulation $(3,0,3,3)$ is given in 
Table \ref{tab303}.

\subsubsection{Case2:$~6(12)$ type}
We will divide this type into two more sub-types:~$6(24)$ and $18(24).$

\subsubsection{Case 2.1: $6(24)$ type}
There exists a triangulation $T:(3,0,3,6+24m)$ for all $m \geq 1.$
This can be constructed as follows:\\
Glueing of the patches $(1,0,3,27)_8$ and $(2,0,0,3)_8$ gives 
a triangulation $T_1:(3,0,3,30+8h)$ for $h \geq 0$.   
If $h=3m$, then we can get the triangulation 
$(3,0,3,30+24m)$ for $m \geq 0$,   
which is the same as the triangulation $T:(3,0,3,6+24m)$ for $m \geq 1$.   
The patch $(1,0,3,27)_8$ can be constructed by using the generic
construction method for $\beta_5=1$
three times on the patch $(1,0,3,6)_5$.    
The patch $(1,0,3,6)_5$ is given in Table 
\ref{penbd} whereas  
the patch $(2,0,0,3)_8$ is a special case of patches of the type $(2,0,0)$
with $k=4,r=0$ (as discussed in Section \ref{patch200}).
The triangulation $T$ gives all the values of the $6(24)$ type, 
except the initial
value of $N=6$. 
For $N=6,$ the triangulation is given in Table \ref{tab303} above.

\subsubsection{Case 2.2: $18(24)$ type}
There exists a triangulation $T:(3,0,3,18+24m)$ for all $m \geq 0.$
This can be constructed as follows:\\
Glueing of the patches $(1,0,3,15)_8$ and $(2,0,0,3)_8$ gives 
a triangulation $T_1:(3,0,3,18+8h)$ for $h \geq 0.$
If $h=3m,$ then we can get the triangulation  
$T:(3,0,3,18+24m)$ for $m \geq 0.$ 
The construction of the patch $(1,0,3,15)_8$ is given in Table $2.6$ and
$(2,0,0,3)_8$ is a special case of patches of the type $(2,0,0)$
with $k=4,r=0$ (as discussed in Section \ref{patch200}). 
The triangulation $T$ gives all the values of the $18(24)$ type.

\subsection{$2(6)$ type}
There exists a triangulation $T:(3,0,3,2+6m)$ for all $m \geq 1.$
This can be constructed as follows: \\
Glueing of the patches $(1,0,3,12)_6$ and $(2,0,0,2)_6$ gives a 
triangulation $T_1:(3,0,3,14+6m)$ for $m \geq 0,$ which is the same as   
the triangulation $T:(3,0,3,2+6m)$ for $m \geq 2$.   
The patch $(1,0,3,12)_6$ can be constructed from the patch
$(1,0,3,6)_5$ by using the generic construction method
with $\beta_5=1$.  
The patch $(1,0,3,6)_5$ is given in Table 2.3   
whereas  
$(2,0,0,2)_6$ is a special case of patches of the type $(2,0,0)$
with $k=3,r=0$ (as discussed in Section \ref{patch200}). 
The triangulation $T$ gives all the values of the $2(6)$ type, except the two initial
values of $N=2,8.$ For $N=8,$ the triangulation is given in Table \ref{tab303} above.
For $N=2,$ the triangulation $(3,0,3,2)$ cannot be constructed. 
We can argue this as follows: \\ 

Let $X=\{abcdefg \}$ with the triangles $abc,acd,ade,aef,afg,abg$ and let the degree of the point $a$ be $d(a)=6.$ 
\\

\noindent
Case (I): Suppose that the points $b,\,c,\,d,\,e,\,f,\,g$ all have degree  
$<6$. If a point $h$ of degree $6$ is included in this construction, then
we may assume (without loss of generality) that $h$ is adjacent to $e$ and $d$,
and there must be $5$ more triangles $hfg,\,hef,\,hcd,\,hbc,\,hbg$.
But then the points $b,c,d,e,f,g$ will have degree $4$. Thus the
triangulation $(3,0,3,2)$ cannot be constructed in this case.
\\

\noindent
Case (II): Suppose without loss of generality
that the point $e$ has degree $6$. Let $h$ be another
point to be added to this construction (so that $h$ has degree
$3$ or $5$). Suppose first that $h$ is adjacent
to $e$. Without loss of generality, we may assume 
that $hed$ forms a triangle. Since $h$ has degree $3$ or $5$ 
and we do not want to add any new points to this construction, the point 
$e$ will end up having degree $3$ or $4$; a contradiction.
Suppose now that $h$ is not adjacent to $e$. 
Again, since $h$ has degree $3$ or $5$ and we do not want to add any new     
points to this construction, the point $e$ will 
end up having degree $2,3$ or $4$ 
which is also a contradiction.  
\\

Hence the triangulation $(3,0,3,2)$ cannot be constructed.

\subsection{$4(6)$ type}
There exists a triangulation $T:(3,0,3,4+6m)$ for all $m \geq 1.$
This can be constructed as follows:\\
Glueing of the patches $(1,0,3,8)_6$ and $(2,0,0,2)_6$ gives 
a triangulation $(3,0,3,$ $10+6h)$ for $h \geq 0,$ which 
is the same as the triangulation $T:(3,0,3,4+6m)$ for $m \geq 1$.   
The construction of the patch $(1,0,3,8)_6$ is given in Table $2.4$
and $(2,0,0,2)_6$ is a special case of patches of the type $(2,0,0)$
with $k=3,r=0$ (as discussed in Section \ref{patch200}).   

The triangulation $T$ gives all the values of the $4(6)$ type,
except the initial value of $N=4.$ We have not been able to construct 
a triangulation
for the case $N=4.$

\section{Triangulations of type $(3,1,1)$} 
It is known by Gr\"unbaum and Motzkin \cite{GM} (see also \cite{GR})
that the triangulation $(3,1,1,2m)$ does not exist for all $m \geq 0.$
In this section we show that the triangulation $(3,1,1,N)$ exists for 
all odd values of $N$ except for $N=1$ and possibly $N=17$.   
We will consider two types: $1(4)$ and $3(4).$

\begin{table}[h]
\begin{center}
\begin{tabular}{|c|l|} \hline 
$N=a_6 $  &  Triangulation \\ \hline
3  & $~ \{ 124,~126,~135,~145,~138,~168,~235,~245,~237,~267,~367,~368 ~\} $ \\ \hline 

5  & $~ \{ 123,~126,~137,~158,~178,~156,~239,~249,~24A,~26A,~349, $ \\
   & $ ~~~345,~357,~456,~46A,~578 ~\} $ \\ \hline 
 
9  & $~ \{ 124,~12C,~134,~13B,~1AC,~1AB,~234,~235,~25C,~357,~37B, $ \\
   & $ ~~~568,~56C,~587,~678,~679,~69D,~6CD,~79B,~9AE,~9BE,~9AD, $ \\
   & $ ~~~AEB,~ACD ~ \}  $  \\ \hline 
     
\end{tabular} 
\end{center}
\caption{Some triangulations of type $(3,1,1)$}
\label{tab311}
\end{table}

\subsection{$1(4)$ type}
To construct triangulations of this type, we will consider two 
sub-types: $1(8)$ and $5(8).$

\subsubsection{$1(8)$ type}
We will divide this type into three sub-types, namely, $1(24),9(24),17(24).$

\subsubsection{Case 1: $1(24)$}
There exists a triangulation $T:(3,1,1,1+24m)$ for all $m \geq 1.$
This can be constructed as follows: \\
Glueing of the patches $(1,1,1,5)_6$ and $(2,0,0,2)_6$ gives a triangulation
$T_1 :(3,1,1,7+6h)$ for $h \geq 0$.    
If $h=3+4m$ for $m \geq 0,$
then we can get the triangulation $(3,1,1,25+24m)$ for $m \geq 0,$
which is the same as the triangulation 
$T :(3,1,1,1+24m)$ for $m \geq 1$.   
The construction of the patch $(1,1,1,5)_6$ 
is given in Table $2.4$ and 
$(2,0,0,2)_6$ is a special case of patches of the type $(2,0,0)$
with $k=3,r=0$ (as discussed in Section \ref{patch200}).  

The triangulation $T$ gives all the values of the $1(24)$ sub-type, 
except the initial value of $N=1.$ For $N=1,$ the triangulation cannot be
constructed as has been shown by Eberhard \cite{Eber} 
and Br\"uckner \cite{bruck}  
(see also Gr\"unbaum \cite{GR}).

\subsubsection{Case 2: $9(24)$}
There exists a triangulation $T:(3,1,1,9+24m)$ for all $m \geq 1.$
This can be constructed as follows: \\
Glueing of the patches $(1,1,1,13)_6$ and $(2,0,0,2)_6$ gives a triangulation
$T_1 :(3,1,1,15+6h)$ for $h \geq 0.$ 
If $h=3+4m$ for $m \geq 0$, then
we can get the triangulation $(3,1,1,33+24m)$ for $m \geq 0,$
which is the same as the triangulation $T:(3,1,1,9+24m)$ 
for $m \geq 1.$ 
The construction of the patch $(1,1,1,13)_6$ is given in Table $2.4$
and $(2,0,0,2)_6$ is a special case of patches of the type $(2,0,0)$ with $k=3,r=0$
(as discussed in Section \ref{patch200}).  

The triangulation $T$ gives all the values of the $9(24)$ sub-type, except
the initial value of $N=9.$ For $N=9,$ the triangulation is given in Table \ref{tab311} above.

\subsubsection{Case 3: $17(24)$}
To construct this $17(24)$ sub-type, we will consider $5$ 
cases: $17(120),41(120),$ $65(120),89(120),113(120).$ 

\subsubsection{Case 3.1: $41(120)$}
There exists a triangulation $T:(3,1,1,41+120m)$ for all $m \geq 1.$
This can be constructed as follows: \\
Glueing of the patches $(1,1,1,72)_{20}$ and $(2,0,0,9)_{20}$ gives a 
triangulation $T_1 :(3,1,1,81+20h)$ for  
$h \geq 0$. If $h=4+6m$ for $m \geq 0,$ then
we can get the triangulation $T:(3,1,1,161+120m)$ for $m \geq 0$
which is the same as the triangulation $T:(3,1,1,41+120m)$ for $m \geq 1.$
The patch $(1,1,1,72)_{20}$ 
can be constructed from the patch $(1,1,1,10)_8$
by using the generic construction method with 
$\beta_4=1$ repeatedly for four times.   
The construction of the patch $(1,1,1,10)_8$
is given in Table $2.6.$
The patch $(2,0,0,9)_{20}$ is a special case of patches of the type $(2,0,0)$
with $k=10,r=0$ (as discussed in Section \ref{patch200}). 
The triangulation $T$ gives all the values of the $41(120)$ case, except the initial
value $N=41.$ For $N=41,$ the triangulation $(3,1,1,41)$ can be
constructed as follows:\\
By glueing the patches $(1,1,1,35)_{14}$ and $(2,0,0,6)_{14}$,    
we can get the triangulation $(3,1,1,41)$.  
The patch $(1,1,1,35)_{14}$ 
can be constructed from the patch $(1,1,1,$ $10)_8$
by using the generic construction method with $\beta_4=1$ twice.
The construction of the patch $(1,1,1,10)_8$
is given in Table $2.6.$ 
The patch $(2,0,0,6)_{14}$ is a special case of patches of the type $(2,0,0)$
with $k=7,r=0$ (as discussed in Section \ref{patch200}).\\

\subsubsection{Case 3.2: $89(120)$}
There exists a triangulation $T:(3,1,1,89+120m)$ for all $m \geq 0.$
This can be constructed as follows: \\
Glueing of the patches $(1,1,1,15)_{10}$ and $(2,0,0,4)_{10}$ 
gives a triangulation $T_1 :(3,1,1,19+10h)$ for $h \geq 0.$ 
If $h=7+12m$ for $m \geq 0,$ then
we can get the triangulation $T:(3,1,1,89+120m)$ for $m \geq 0$.   
The construction of the patch $(1,1,1,15)_{10}$ is given in Table $2.6$ and 
$(2,0,0,4)_{10}$ is a special case of patches of the type $(2,0,0)$ with $k=5,r=0$ 
(as discussed in Section \ref{patch200}).
The triangulation $T$ gives all the values of the $89(120)$ case.

\subsubsection{Case 3.3: $113(120)$}
There exists a triangulation $T:(3,1,1,113+120m)$ for all $m \geq 0.$
This can be constructed as follows: \\
Glueing of the patches $(1,1,1,29)_{10}$ and $(2,0,0,4)_{10}$ 
gives a triangulation $T_1 :(3,1,1,33+10h)$ for $h \geq 0$.  
If $h=8+12m$ for $m \geq 0,$ then
we can get the triangulation $T:(3,1,1,113+120m)$ for $m \geq 0.$
The patch $(1,1,1,29)_{10}$ 
can be constructed from the patch $(1,1,1,10)_8$
by using the generic construction method with $\beta_5=1$ twice. 
The construction of the patch $(1,1,1,10)_8$
is given in Table $2.6.$
The patch $(2,0,0,4)_{10}$ is a special case of patches of the type $(2,0,0)$
with $k=5,r=0$ (as discussed in Section \ref{patch200}).
The triangulation $T$ gives all the values of the $113(120)$ case.

\vspace{0.5cm}
For the construction of the $17(120)$ and $65(120)$ cases,
we will consider triangulations of the form 
$(3,1,1,X+14N)$ and $(3,1,1,Y+20N)$ for some values of $X$ and $Y.$
The least common multiple of the numbers $14,20,120$  is $840.$ 
Therefore, we have to consider the following $14$ sub-cases: 
$17(840)$, $137(840)$, $257(840)$, $377(840)$, $497(840)$,  
$617(840)$, $737(840)$, $65(840)$, $185(840)$, $305(840)$,    
$425(840)$, $545(840)$, $665(840)$, $785(840)$.   

\subsubsection{Case (i): $17(840)$}
There exists a triangulation $T:(3,1,1,17+840m)$ for all $m \geq 1.$
This can be constructed as follows: \\
Glueing of the patches $(1,1,1,25)_{14}$ and $(2,0,0,6)_{14}$ 
gives a triangulation $T_1 :(3,1,1,31+14h)$ for $h \geq 0.$ 
If $h=59+60m$ for $m \geq 0$, then
we can get the triangulation $(3,1,1,857+840m)$ for $m \geq 0$, which
is the same as the triangulation $T :(3,1,1,17+840m)$ for  
$m \geq 1.$ The patch $(1,1,1,25)_{14}$ can be constructed 
from the patch $(1,1,1,11)_{10}$ by using the 
generic construction method with 
$\beta_4=1$ and $\beta_5=1$.  
The construction of the patch $(1,1,1,11)_{10}$
is given in Table $2.6.$
The patch $(2,0,0,6)_{14}$ is a special case of patches of the type $(2,0,0)$
with $k=7,r=0$ (as discussed in Section \ref{patch200}).
The triangulation $T$ gives all the values of the $17(840)$ case, except 
for the initial value $N=17.$ 
We have not been able to construct a triangulation
for the case $N=17.$

\subsubsection{Case (ii): $185(840)$}
There exists a triangulation $T:(3,1,1,185+840m)$ for all $m \geq 0.$
This can be constructed by substituting $h=11+60m$ for $m \geq 0$
in the triangulation $T_1 :(3,1,1,31+14h)$,     
which is as given in the construction of case (i) above.
The triangulation $T$ gives all the values of the $185(840)$ case.

\subsubsection{Case (iii): $377(840)$}
There exists a triangulation $T:(3,1,1,377+840m)$ for all $m \geq 0.$
This can be constructed as follows: \\
Glueing of the patches $(1,1,1,35)_{14}$ and $(2,0,0,6)_{14}$ gives 
a triangulation $T_1 :(3,1,1,41+14h)$ for  
$h \geq 0.$ If $h=24+60m$ for $m \geq 0,$ then
we can get the triangulation $T :(3,1,1,377+840m)$ for $m \geq 0.$
The patch $(1,1,1,35)_{14}$ can be constructed from 
the patch $(1,1,1,10)_8$ by using 
the generic construction method with $\beta_4=1$ twice.
The construction of the patch $(1,1,1,10)_8$ is given in Table $2.6$.   
The patch $(2,0,0,6)_{14}$ is a special case of patches of the type $(2,0,0)$
with $k=7,r=0$ (as discussed in Section \ref{patch200}).
The triangulation $T$ gives all the values of the $377(840)$ case.

\subsubsection{Case (iv): $545(840)$}
There exists a triangulation $T:(3,1,1,545+840m)$ for all $m \geq 0.$
This can be constructed by substituting $h=36+60m$ for $m \geq 0$
in the triangulation $T_1 :(3,1,1,41+14h)$, which is as 
given in the construction of case (iii) above.
The triangulation $T$ gives all the values of the $545(840)$ case.

\subsubsection{Case (v): $257(840)$}
There exists a triangulation $T:(3,1,1,257+840m)$ for all $m \geq 0.$
This can be constructed as follows: \\
Glueing of the patches $(1,1,1,41)_{14}$ and $(2,0,0,6)_{14}$ 
gives a triangulation $T_1 :(3,1,1,47+14h)$ for $h \geq 0$.    
If $h=15+60m$ for $m \geq 0$, then
we can get the triangulation $T :(3,1,1,257+840m)$ for $m \geq 0.$
The patch $(1,1,1,41)_{14}$ can be constructed from the 
patch $(1,1,1,14)_{12}$ by using  
the generic construction method with $\beta_5=1$ twice. 
The construction of the patch $(1,1,1,14)_{12}$ is given in Table $2.6$.   
The patch $(2,0,0,6)_{14}$ is a special case of patches of the type $(2,0,0)$
with $k=7,r=0$ (as discussed in Section \ref{patch200}).
The triangulation $T$ gives all the values of the $257(840)$ case. 

\subsubsection{Case (vi): $425(840)$}
There exists a triangulation $T:(3,1,1,425+840m)$ for all $m \geq 0.$
This can be constructed by substituting $h=27+60m$ for $m \geq 0$
to the triangulation $T_1 :(3,1,1,47+14h)$,     
which is as given in the construction of case (v) above.
The triangulation $T$ gives all the values of the $425(840)$ case. 

\subsubsection{Case (vii): $305(840)$}
There exists a triangulation $T:(3,1,1,305+840m)$ for all $m \geq 0.$
This can be constructed as follows: \\
Glueing of the patches $(1,1,1,61)_{14}$ and $(2,0,0,6)_{14}$ gives 
a triangulation 
$T_1 :(3,1,1,67+14h)$ for $h \geq 0$. If $h=17+60m$ 
for $m \geq 0$, then
we can get the triangulation $T:(3,1,1,305+840m)$.     
The patch $(1,1,1,61)_{14}$ can be constructed from 
the patch $(1,1,1,11)_{10}$ by using 
the generic construction method
with $\beta_5=1$ repeatedly for four times.   
The construction of the patch $(1,1,1,11)_{10}$
is given in Table $2.6.$ 
The patch $(2,0,0,6)_{14}$ is a special case of patches of the type $(2,0,0)$
with $k=7,r=0$ (as discussed in Section \ref{patch200}).
The triangulation $T$ gives all the values of the $305(840)$ case. 

\subsubsection{Case (viii): $737(840)$}
There exists a triangulation $T:(3,1,1,737+840m)$ for all $m \geq 0.$
This can be constructed as follows: \\
Glueing of the patches $(1,1,1,101)_{14}$ and $(2,0,0,6)_{14}$ 
gives a triangulation $T_1 :(3,1,1,107+14h)$ for  
$h \geq 0.$ If $h=45+60m$ for $m \geq 0,$ then
we can get the triangulation $T:(3,1,1,737+840m)$ for $m \geq 0.$ 
The patch $(1,1,1,101)_{14}$ can be constructed from 
the patch $(1,1,1,2)_3$ by using 
the generic construction method with $\beta_5=1$ repeatedly for 
$11$ times.    
The patch $(1,1,1,2)_3$
is given in Table 2.1 whereas   
the patch $(2,0,0,6)_{14}$ is a special case 
of patches of the type $(2,0,0)$
with $k=7,r=0$ (as discussed in Section \ref{patch200}).
The triangulation $T$ gives all the values of the $737(840)$ case. 

\subsubsection{Case (ix): $65(840)$}
There exists a triangulation $T:(3,1,1,65+840m)$ for all $m \geq 1.$
This can be constructed as follows:\\
If we substitute $h=57+60m$ for $m \geq 0$
in the triangulation $T_1 :(3,1,1,107+14h)$,   
which is as given in the construction of case (viii) above, then
we can get the triangulation $(3,1,1,905+840m)$ for $m \geq 0$ which
is the same as the triangulation $T :(3,1,1,65+840m)$ for $m \geq 1.$ 
The triangulation $T$ gives all the values of the $65(840)$ case, 
except the initial value of $N=65.$ For $N=65,$ the 
triangulation $T:(3,1,1,65)$ can be constructed as follows: \\
Glueing of the patches $(1,1,1,11)_{10}$ and $(2,0,0,4)_{10}$
gives a triangulation $(3,1,1,15+10h)$ for $h \geq 0.$ 
If $h=5,$ then we can get the triangulation $(3,1,1,65).$
The construction of the patch $(1,1,1,11)_{10}$ is given in Table
$2.6$ whereas $(2,0,0,4)_{10}$ 
is a special case of patches of the type $(2,0,0)$
with $k=5,r=0$ (as discussed in Section \ref{patch200}).\\

\subsubsection{Case (x): $617(840)$}
There exists a triangulation $T:(3,1,1,617+840m)$ for all $m \geq 0.$
This can be constructed as follows: \\
Glueing of the patches $(1,1,1,65)_{14}$ and $(2,0,0,6)_{14}$ 
gives a triangulation $T_1 :(3,1,1,71+14h)$ for $h \geq 0.$ 
If $h=39+60m$ for $m \geq 0,$ then
we can get the triangulation $T:(3,1,1,617+840m)$ for $m \geq 0.$ 
The patch $(1,1,1,65)_{14}$ can be constructed from 
the patch $(1,1,1,15)_{10}$ by using  
the generic construction method with $\beta_5=1$ repeatedly for
four times.   
The construction of the patch $(1,1,1,15)_{10}$
is given in Table $2.6$.   
The patch $(2,0,0,6)_{14}$ is a special case of patches of the type $(2,0,0)$
with $k=7,r=0$ (as discussed in Section \ref{patch200}).
The triangulation $T$ gives all the values of the $617(840)$ case. 

\subsubsection{Case (xi): $785(840)$}
There exists a triangulation $T:(3,1,1,785+840m)$ for all $m \geq 0.$
This can be constructed by substituting $h=51+60m$ for $m \geq 0$
in the triangulation $T_1 :(3,1,1,71+14h)$, which is as 
given in the construction of case (x) above.
The triangulation $T$ gives all the values of the $785(840)$ case.

\subsubsection{Case (xii): $137(840)$}
There exists a triangulation $T:(3,1,1,137+840m)$ for all $m \geq 0.$
This can be constructed as follows: \\
Glueing of the patches $(1,1,1,61)_{14}$ and $(2,0,0,6)_{14}$ 
gives a triangulation $T_1 :(3,1,1,67+14h)$ for $h \geq 0.$ 
If $h=5+60m$ for $m \geq 0,$ then
we can get the triangulation $T:(3,1,1,137+840m)$ for $m \geq 0$.   
The patch $(1,1,1,61)_{14}$ can be constructed from the patch 
$(1,1,1,11)_{10}$ by using  
the generic construction method with $\beta_5=1$ 
repeatedly for four times.   
The construction of the patch $(1,1,1,11)_{10}$ is given in Table $2.6$.   
The patch $(2,0,0,6)_{14}$ is a special case of patches of the type $(2,0,0)$
with $k=7,r=0$ (as discussed in Section \ref{patch200}).
The triangulation $T$ gives all the values of the $137(840)$ case.

\subsubsection{Case (xiii): $497(840)$}
There exists a triangulation $T:(3,1,1,497+840m)$ for all $m \geq 0.$
This can be constructed as follows: \\
Glueing of the patches $(1,1,1,29)_{14}$ and $(2,0,0,6)_{14}$ gives a triangulation $T_1 :(3,1,1,35+14h)$ for all $h \geq 0.$ 
If $h=33+60m$ for $m \geq 0,$ then
we can get the triangulation $T:(3,1,1,497+840m)$.    
The patch $(1,1,1,29)_{14}$ can be constructed from 
the patch $(1,1,1,15)_{10}$ by using the generic construction method with 
$\beta_4=\beta_5=1$. 
The construction of the patch $(1,1,1,15)_{10}$ is given in Table $2.6$.  
The patch $(2,0,0,6)_{14}$ is a special case of patches of the type $(2,0,0)$
with $k=7,r=0$ (as discussed in Section \ref{patch200}).
The triangulation $T$ gives all the values of the $497(840)$ case. 

\subsubsection{Case (xiv): $665(840)$}
There exists a triangulation $T:(3,1,1,665+840m)$ for all $m \geq 0.$
This can be constructed by substituting $h=45+60m$ for $m \geq 0$
in the triangulation $T_1 :(3,1,1,35+14h)$, where $T_1$ is as      
given in the construction of case (xiii) above.
The triangulation $T$ gives all the values of the $665(840)$ case. 

\vspace{0.5cm}
From the $14$ cases above, we have shown the existence of triangulations 
of the 
type $(3,1,1,N)$ where $N \equiv 17(24)$ and $N \equiv 65(24)$, 
except for $N=17$.  
We have not been able to show the existence (or non-existence) 
of a triangulation
of the type $(3,1,1,17).$
\medskip

\subsubsection{$5(8)$ type}

There exists a triangulation $T:(3,1,1,5+8m)$ for all $m \geq 1.$
This can be constructed as follows: \\
Glueing of the patches $(1,1,1,10)_8$ and $(2,0,0,3)_8$
gives a triangulation  $(3,1,1$, $13+8m)$ for  
$m \geq 0$ which 
is the same as the triangulation $T:(3,1,1,5+8m)$ for $m \geq 1$.   
The construction of the patch $(1,1,1,10)_8$ is given in Table
$2.6$ and $(2,0,0,3)_8$ 
is a special case of patches of the type $(2,0,0)$
with $k=4,r=0$ (as discussed in Section \ref{patch200}). 
The triangulation $T$ gives all the values of the $5(8)$ type, except
the initial value $N=5$. For $N=5$, the triangulation 
$(3,1,1,5)$ is given
in Table \ref{tab311}.

\subsection{$3(4)$ type}
The triangulation $(3,1,1,3+4m)$ exists for all $m \geq 0.$
This can be constructed by glueing of the patches 
$(1,1,1,2)_4$ and $(2,0,0,1)_4.$ The construction of the patch
$(1,1,1,2)_4$ is given in Table \ref{rectbd} and
$(2,0,0,1)_4$ is a special case of patches of the type $(2,0,0)$
with $k=2,r=0$ (as discussed in Section \ref{patch200}).
The triangulation $T$ gives all the values of the $3(4)$ type.

\section{Triangulations of type $(4,0,0)$}  \label{Tri400}
In this section we show that triangulations of the type $(4,0,0,N)$
exist for all even values of $N$ except for $N=2.$ In the case $N=2$
or $N$ is odd, Gr\"unbaum and Motzkin \cite{GM} (see also \cite{GR})
have shown that triangulations of the type $(4,0,0,N)$ do not exist.
A triangulation $T:(4,0,0,2m)$ for $m \geq 3$ can be constructed as 
follows:\\
By taking two copies of the patch $(2,0,0,k-1)_{2k}$ 
and glueing them with $b>0$ belts, we can get the triangulation  
$T_1=(4,0,0,2(k-1)+b \cdot 2k).$ 
Here $k \geq 2$ is the edge distance 
between the points of degree $3$ of the patch.   

The triangulation $T_1$ is the same as the triangulation
$T=(4,0,0,2m)$ for all $m\geq 3$.   
The patch $(2,0,0,k-1)_{2k}$ is 
a special case of patches of the type $(2,0,0)$
with $r=0$ (refer to Section \ref{patch200}). 
For $N=4$, the triangulation $(4,0,0,4)$ is as given in Table 3.19
below.    
\medskip

\begin{table}[h]
\begin{center}
\begin{tabular}{|c|l|} \hline 
$N=a_6 $  &  Triangulation \\ \hline
0  & $~ \{ 123,~124,~134,~234 ~\} $ \\ \hline 

4  & $~ \{ 124,~126,~145,~135,~137,~167,~246,~348,~345,~367,~368,~468 ~\} $ \\ \hline 

6  & $~ \{ 126,~127,~137,~138,~146,~148,~235,~237,~259,~269,~345,~348, $ \\
   & $ ~~~45A,~46A,~569,~56A ~\} $ \\ \hline 
 
10  & $~ \{ 125,~12B,~15B,~256,~236,~23C,~2BC,~367,~347,~34D,~3CD, $ \\
   & $ ~~~47D,~5AB,~58A,~568,~689,~679,~789,~78D,~8AD,~ABE, $ \\
   & $ ~~~AEC,~ACD,~BCE ~\} $ \\ \hline

\end{tabular} 
\end{center}
\caption{Some triangulations of type $(4,0,0)$}
\end{table}

\newpage

\addcontentsline{toc}{4}{Appendix \hspace{11.0cm}}
\Large 
\centerline{\bf{Appendix}}
\normalsize
\bigskip

\begin{table}[h]
\begin{center}
\begin{tabular}{|c|l|c|c|} \hline
\multicolumn{1}{|c|}{Triangulation}&
\multicolumn{3}{|c|}{Values of $N$ where the type $T$}\\ 
 type $(T)$ & (i) exists & (ii) does not exist &(iii) is not known \\
 \hspace{2.5cm} &  		 	&  & yet to exist \\ \hline

$(0,0,12,N)$  &$0,2,3,4,\dots$ & $1$  & $-$  \\ \hline 
$(0,1,10,N)$  & $2,3,\dots $ & $0,1$  & $-$  \\ \hline 
$(0,2,8,N)$  & $0,1,2,\dots $ & $-$   & $-$  \\ \hline 
$(0,3,6,N)$  & $0,1,2,\dots $ & $-$   & $-$  \\ \hline 
$(0,4,4,N)$  & $0,2,3,4,\dots $ & $1$   & $-$  \\ \hline 
$(0,5,2,N)$  & $0,2,3,4,\dots $ & $1$   & $-$  \\ \hline 
$(0,6,0,N)$  & $0,2,3,4,\dots $ & $1$   & $-$  \\ \hline 
$(1,0,9,N)$  & $3,5,6,7,\dots $ & $0,1,2$   & $4$  \\ \hline 
$(1,1,7,N)$  & $2,3,\dots $ & $0,1$   & $-$  \\ \hline 
$(1,2,5,N)$  & $1,2,\dots $ & $0$   & $-$  \\ \hline 
$(1,3,3,N)$  & $0,1,2,\dots $ & $-$   & $-$  \\ \hline 
$(1,4,1,N)$  & $2,3,\dots $ & $0,1$   & $-$  \\ \hline 
$(2,0,6,N)$  & $0,2,3,\dots $ & $1$   & $-$  \\ \hline 
$(2,1,4,N)$  & $1,2,\dots $ & $0$   & $-$  \\ \hline 
$(2,2,2,N)$  & $0,1,2,\dots $ & $-$   & $-$  \\ \hline 
$(2,3,0,N)$  & $N \neq 1,3,7,15,31$ & $1$   & $3,7,15,31$  \\ \hline 
$(3,0,3,N)$  & $N \neq 0,2,4,12$ & $0,2$   & $4,12$  \\ \hline 
$(3,1,1,N)$  & $N$ odd, $N \neq 1,17$ & $N$ even, $N=1$   &$17$  \\ \hline 
$(4,0,0,N)$  & $N$ even, $N \neq 2$ & $N$ odd, $N=2$   & $-$  \\ \hline 

\end{tabular} 
\end{center}
\caption{Elliptic triangulations of spheres} 
\end{table}

\newpage
\addcontentsline{toc}{4}{Bibliography \hspace{11.0cm}}

\end{document}